\documentclass{amsart}
\usepackage[utf8]{inputenc}
\usepackage[T2A,T1,OT1]{fontenc}
\usepackage[russian,british]{babel}
\usepackage{uniinput}
\usepackage[cmyk]{tudcolors}
\usepackage{tikz}[2010/10/13]
\usetikzlibrary{matrix,arrows,decorations.pathmorphing,decorations.pathreplacing,decorations.footprints}
\usepackage{tsabk}
\usepackage{tabularx}
\usepackage{subfigure}
\usepackage{hyperref}
\usepackage{amssymb}
\newtheorem{theorem}{Theorem}
\newtheorem{lemma}[theorem]{Lemma}
\newtheorem{corollary}[theorem]{Corollary}
\theoremstyle{definition}
\newtheorem{definition}{Definition}
\theoremstyle{remark}

\author[T. Schlemmer]{Tobias Schlemmer}
\address{Department of Mathematics\\Technische Universit\"at Dresden}
\email{Tobias.Schlemmer@gmx.de}
\urladdr{http://www.math.tu-dresden.de/~schlemme/}
\title[Linear extensiones on abelian po-groups]{Linear extensions of
  partial orders \protect\\on abelian groups}
\thanks{This work was supported by TU Dresden, Professur f. Angewandte
Algebra}
\DeclareMathOperator\projektionssymbol{P}
\newcommand\projektion[2]{\projektionssymbol_{#1,#2}}
\newcommand\alert{\emph}
\newcommand\defindex[2][]{\emph{#2}}
\newcommand{\svnid}[1]{}
\DeclareMathOperator\codim{codim}
\newcommand\translit[1]{{\fontencoding{T1}\selectfont{}#1}}
\newcommand*{\ordnungsideal}[2][\leq]{\mathop{\downarrow_{#1}}{#2}}
\DeclareMathOperator\Act{Act}
\DeclareMathOperator\vAct{\mathfrak{Act}}
\tikzstyle{alert1}=[]%
\tikzstyle{alert2}=[]%
\tikzstyle{alert3}=[]%
\tikzstyle{alert4}=[]%
\tikzstyle{alert5}=[]%
\tikzstyle{alert6}=[]%

\addto\captionsbritish{\def\ggTname{gcd}\def\kgVname{lcm}}

\begin{document}
\bibliographystyle{plain}
\maketitle
\begin{abstract}
  Partially ordered groups, also known as po-groups, are groups with a
  compatible partial order. Results from \translit{M.\,I. Zajceva} and
  H.-H.\ Teh are combined in order to provide a full characterisation
  of linear order extensions of a given order on a group. In contrast
  to Teh this approach provides a method to discuss linear orders of
  different abelian rank in a uniform manner. This will be achieved by
  modelling the linear orders using hyperplanes in a real vector
  spaces. Among some additional remarks a construction of an
  archimedian directed order is given for every torsion free abelian
  group.
\end{abstract}

\section*{Introduction}

Linear and lattice ordered groups are well-studied algebraic
structures. The structure of linearly ordered abelian groups has
mainly been discovered by the work of Otto Hölder \cite{Holder:1901} and
Hans Hahn \cite{Hahn:1907}. Later, Friedrich Levi
\cite{springerlink:10.1007/BF03015602} provided a first
characterisation of lattice ordered groups and showed that an abelian group
can be linearly ordered iff it is torsion free. It must be noticed
that at the time of  these early works the theory of modules and vector spaces
had not settled down. First works on these topics from the
middle of the 19th century have rarely been noticed in the
community. It took nearly 100 years until todays definitions for
vector spaces and modules had been fixed.

\textru{Анатолий Иванович Мальцев}~\cite{Malcev:1949}
(\translit{A.\,I. Maľcev}) investigated necessary and sufficient
conditions for linear orders on abelian groups. Later, 
\textru{М.\,И. Зайцева} \cite{M.I.Zajceva:1953}
(\translit{M.\,I.~Zajceva}) published a
characterisation of finitely generated archimedian linearly ordered abelian groups.
With a corresponding decomposition into archimedian subgroups, in that work she
discusses also the main properties of a description of
linear orders on finitely generated groups. The cardinality of such
orders on a given group has been determined by Shin-Ichi \cite{Matsusita:1953} to
$\aleph$. Eventually, Hoon-Heng Teh \cite{Teh:1961} provided another
characterisation of linearly ordered groups based on Hahn's theorem.

The current work provides a characterisation of linear order
extensions of abelian partially ordered groups, which has been
developed from scratch, based on the usage and an example from Charles
Holland and other authors
(cf. e.\,g., \cite{springerlink:10.1007/BF01190757}).  Levi discussed
the characteristics of lattice ordered groups starting with abelian
groups embedded into one and two dimensional vector spaces and
analysed them with rising dimension (bottom-up approach). He provides
tools to describe the image of a positive cone of a lattice ordered
group in a vector space by means of rays and edges. We will discuss
the same phenomenon by construction of hyperplanes in a vector space
that are constructed using independent sets of a given group (top-down
approach). Though, this method considers the positive cone in less granularity,
it provides a better bridge to methods from convex geometry than
Levi's approach.

\translit{Zajceva}'s equations for archimedian orders on $n$-generated
abelian groups can be interpreted as certain
hyperplanes in $n$-dimensional real vector spaces. Thus, this work can also be considered an extension of her results
onto arbitrary abelian groups and allows to discuss existing orders
on them.  A further advantage of this view consists in the possibility
to investigate archimedian and several non-archimedian orders together
based on a common mathematical structure. This facet and the discussion of
linear order extensions of given order relations are also an enhancement to
Teh's work, who discusses the archimedian rank of linear orders.

The constructive nature of the given article is dedicated to be more
accessible for non-mathematicians and mathematicians with a background
different from group theory. This is achieved by using injective
abelian groups, which need different tools and provide slightly
different insights than the usual upproach using free abelian
$\ell$-groups or free vector lattices and the duality theorems
provided by W.\,M. Beynon~\cite{Beynon01071975}.

To achieve this, after a short clearing of notions, mappings between
orders on different structures will be investigated. This discussion
has the aim to enable us to represent a given order on a torsion free
abelian group in a convenient way in the vector space of the direct
sum $\direktsumme{E}\reellezahlen$, where the set of indices $E$
corresponds to a maximal independent set in the given group. There we
will use a characterisation of linearly ordered groups that is a
little bit different from those given by \translit{Zajceva} and Teh
(cf. Theorem
\ref{satz:charakterisierung-linearer-ordnungen}). Finally, with
Theorem \ref{satz:ordnungserweiterungscharakterisierung} we will
provide a characterisation of linear order extensions on partially
ordered abelian groups, with the use of half spaces and linearly
ordered hyperplanes. Thus, we can use methods from convex geometry to
discuss linearly ordered groups. This will be demonstrated in
Section \ref{sec:further-results}, where some known results have been
resembled. This section also shows how to construct an archimedian directed
order on an arbitrary torsion free abelian group\iffalse, which has been
reported as open problem in 2009 by
V.\,V.~Bludov et.al.~\cite{Bludov:2009}. This discussion also includes an
example of an torsion free abelian group that does not permit an
archimedian linear order\fi.

\section{Preliminaries}

In this section we repeat the basics used in this article. The
main facts can be found in the usual textbooks about lattice ordered
groups (e.\,g. \cite{marlow1988latticeordered,Birkhoff93,r1995theory,fuchs:1963,w1989latticeordered,Kopytov:1994lr,Kokorin:1972ru,KOP84}).

A group $\gruppe = \Gruppe G\cdot\multinverses e$ is called 
\defindex{partially ordered} iff there exists a partial order
${\leq}\subseteq G\times G$ such that the following condition holds
for all  elements $x,y,a,b\in G$:\footnote{The given condition is
  sufficient to assure that the partial order $\leq$ is also compatible with
  the other group operations (inverse element and neutral element).}
\begin{equation}
	a\leq b \Rightarrow xay \leq xby.
\end{equation}
It is called \defindex{linearly ordered} iff $\leq$ is a linear order.
A group homomorphism $φ$ from a partially ordered group \gruppe{} into
another one $\gruppe'$ is called \defindex{o-homomorphism} iff it is
also an order homomorphism. Furthermore it is called an
\defindex{o-isomorphism} if it is additionally an order isomorphism.

A very fundamental theorem allows us to characterise each
partial order $\leq$ on a group $\gruppe$ by means of its \defindex{positive
cone} $\poskegel{\gruppe}:=\Menge{g\in G}{0\leq g}$:
\begin{theorem} \label{satz:poskegelbedingungen}%
The positive cone $\poskegel \gruppe$ of a partially ordered group $\OrdGruppe G
\cdot\multinverses e\leq$ fulfils the following conditions:
\begin{align}
	&\poskegel \gruppe \cdot \poskegel \gruppe \subseteq \poskegel
        \gruppe&&\text{$\poskegel \gruppe$ is a semigroup}\label{eq:poskegelbedingung_halbgruppe}\\
	&\poskegel \gruppe\cap \multnegpot{(\poskegel \gruppe)} = \menge e&&
        \text{$\poskegel \gruppe$ is a pure subset}\label{eq:poskegelbedingung_echte_halbgrupe}\\ %pure subset
        &\forall g\in G: g^{-1}\poskegel \gruppe g \subseteq \poskegel
         \gruppe&&\text{$\poskegel \gruppe$ is invariant subsemigroup}.\label{eq:poskegelbedingung_invariant}
\intertext{If $\gruppe$ is a directed set, it further fulfils,}
&G=\poskegel \gruppe\cdot \multnegpot{(\poskegel \gruppe)}.\label{eq:poskegelbedingung_gerichtet}
\intertext{and if $\gruppe$ is  linearly ordered}
&G=\poskegel \gruppe\cup  \multnegpot{(\poskegel \gruppe)}.\label{eq:poskegelbedingung_linear}
\end{align}
Conversely, if in \gruppe\ a set $P\subseteq G$ exists,
that fulfils \eqref{eq:poskegelbedingung_halbgruppe} to
\eqref{eq:poskegelbedingung_invariant},  there exists an order  on
$\gruppe$ such that $P$ is the set of positive elements of
\gruppe{}. If furthermore the Condition
\eqref{eq:poskegelbedingung_gerichtet} is met, the order is directed
and, if \eqref{eq:poskegelbedingung_linear} is true, it is a linear order.
\end{theorem}

\begin{proof} Cf. \cite{Kopytov:1994lr}, Theorem 2.1.1 or
  \cite{KOP84}, II.2.1.
\end{proof}

An Element $a∈\poskegel{\gruppe}$ is called \defindex{infinitesimal}
with respect to another element $b$ iff the relation $a^n ≤ b$ holds
for every integer $n∈\ganzzahlen$. Consequently $0$ is infinitesimal
with respect to every other element. If there is no other element
infinitesimal with respect to any other element then the group is
called \defindex{archimedian}.

It is a well-known fact that each finite group can
only be discretely ordered. This implies that only torsion free
groups can be linearly ordered. Thus, throughout this paper we will
refer to \textbf{torsion free groups} when ever we talk about
\textbf{groups}, unlike otherwise stated.

An
\defindex{independent set} $E\subseteq G$ of an abelian group
$\gruppe=\Gruppe G+-0$ is defined as a set of its elements, which
fulfils the equation
\begin{equation}
  \label{eq:independent set}
  \sum_{a\in E'}\vektor a(a)a = 0,
 \end{equation}
for each finite subset $E'\subseteq E$ and a mapping $\vektor{a}\in\direktsumme{E}\ganzzahlen$ iff  for all  $a\in E'$ the condition $\vektor a(a)
= 0$ is met.  Zorn's lemma assures the existence of a maximal
independent set. A subset $E\subseteq G$ is independent,
iff the subgroup $\erzgruppe E$, which is generated by $E$, is
isomorphic to the
direct sum of the cyclic groups $\erzgruppe a$ of all of its
elements $a\in E$. A subgroup $S\subseteq G$ is called \defindex{essential subgroup} iff
for each subgroup $S'\subseteq G$ the intersection $S\cap
S'\neq\menge{0}$ is non-trivial. An independent set $E$ is maximal iff its
generated subgroup $\erzgruppe E$ is an essential subgroup of
\gruppe.  For further information we refer to \cite{Fuchs:1970}
or any other text book on the theory of abelian groups. 

Given
an arbitrary  set $E$ we define the direct sum $\direktsumme E\gruppe$ according
to 
\begin{equation}
\direktsumme E\gruppe :=
\Menge[\Big]{\vektor x}{\vektor x:E\to\gruppe, |\supp\vektor
  x|<\infty},
\end{equation}
where $\supp\vektor{x}:=\Menge{a\in E}{\vektor{x}(a)\neq 0}$ is the
support of the mapping $\vektor{x}$. The unit vectors will be denoted by
\begin{equation}
  \label{eq:unit-vector}
  \vektor{e}_a(b):=
  \begin{cases}
    1&\text{iff }a=b,\\
    0&\text{otherwise}.
  \end{cases}
\end{equation}
If the generating set $E$ is not uniquely determined by the context
(e.g. when we are using different bases of a vector space) we will
write $\vektor x_E(e)$ instead of $\vektor x(e)$ for any element
$e∈E$.

\begin{figure}\noindent\parbox[t]{0.4\linewidth}{\caption[Positive
    cone in a vector space.]
{\linebreak[2]\raggedright Positive
    cone in a vector space}\label{fig:poskegel-im-Vektorraum}}%
\parbox[t]{0.6\linewidth}{\centerline{\raisebox{-\height}{% Sketch output, version 0.3 (build 7d, Fri Nov 22 14:14:09 2013)
% Output language: PGF/TikZ,LaTeX
\begin{tikzpicture}[line join=round]
\tikzstyle{ann} = [fill=white,font=\footnotesize,inner sep=1pt]
\providecommand\vektor{}\filldraw[draw=HKS36K80,draw opacity=0.7,line width=0.0001,fill=HKS36K80,fill opacity=0.7](0,0)--(.548,2.207)--(.624,2.247)--(0,0)--cycle;
\filldraw[draw=HKS36K80,draw opacity=0.7,line width=0.0001,fill=HKS36K80,fill opacity=0.7](0,0)--(.549,2.123)--(.548,2.207)--(0,0)--cycle;
\filldraw[draw=HKS36K80,draw opacity=0.7,line width=0.0001,fill=HKS36K80,fill opacity=0.7](0,0)--(.624,2.247)--(.765,2.239)--(0,0)--cycle;
\filldraw[draw=HKS36K80,draw opacity=0.7,line width=0.0001,fill=HKS36K80,fill opacity=0.7](0,0)--(.627,2.01)--(.549,2.123)--(0,0)--cycle;
\filldraw[draw=HKS36K80,draw opacity=0.7,line width=0.0001,fill=HKS36K80,fill opacity=0.7](0,0)--(.765,2.239)--(.951,2.183)--(0,0)--cycle;
\filldraw[draw=HKS36K80,draw opacity=0.7,line width=0.0001,fill=HKS36K80,fill opacity=0.7](0,0)--(.771,1.884)--(.627,2.01)--(0,0)--cycle;
\filldraw[draw=HKS36K80,draw opacity=0.7,line width=0.0001,fill=HKS36K80,fill opacity=0.7](0,0)--(.951,2.183)--(1.152,2.087)--(0,0)--cycle;
\filldraw[draw=HKS36K80,draw opacity=0.7,line width=0.0001,fill=HKS36K80,fill opacity=0.7](0,0)--(.957,1.764)--(.771,1.884)--(0,0)--cycle;
\filldraw[draw=HKS36K80,draw opacity=0.7,line width=0.0001,fill=HKS36K80,fill opacity=0.7](0,0)--(1.339,1.968)--(1.482,1.842)--(0,0)--cycle;
\filldraw[draw=HKS36K80,draw opacity=0.7,line width=0.0001,fill=HKS36K80,fill opacity=0.7](0,0)--(1.152,2.087)--(1.339,1.968)--(0,0)--cycle;
\filldraw[draw=HKS36K80,draw opacity=0.7,line width=0.0001,fill=HKS36K80,fill opacity=0.7](0,0)--(1.159,1.669)--(.957,1.764)--(0,0)--cycle;
\filldraw[draw=HKS36K80,draw opacity=0.7,line width=0.0001,fill=HKS36K80,fill opacity=0.7](0,0)--(1.482,1.842)--(1.561,1.728)--(0,0)--cycle;
\filldraw[draw=HKS36K80,draw opacity=0.7,line width=0.0001,fill=HKS36K80,fill opacity=0.7](0,0)--(1.344,1.613)--(1.159,1.669)--(0,0)--cycle;
\filldraw[draw=HKS36K80,draw opacity=0.7,line width=0.0001,fill=HKS36K80,fill opacity=0.7](0,0)--(1.561,1.728)--(1.562,1.645)--(0,0)--cycle;
\filldraw[draw=HKS36K80,draw opacity=0.7,line width=0.0001,fill=HKS36K80,fill opacity=0.7](0,0)--(1.486,1.605)--(1.344,1.613)--(0,0)--cycle;
\filldraw[draw=HKS36K80,draw opacity=0.7,line width=0.0001,fill=HKS36K80,fill opacity=0.7](0,0)--(1.562,1.645)--(1.486,1.605)--(0,0)--cycle;
\draw[arrows=<->,line width=.4pt](.42,-.767)--(0,0)--(0,2.167);
\draw[arrows=->,line width=.4pt](0,0)--(-2.1,-.153);
\path (1.068,1.95) node {$P$};\path (.42,-.767) node[below] {$\vektor a_1$}
     (0,2.167) node[above] {$\vektor a_2$}
     (-2.1,-.153) node[left] {$\vektor a_3$};\end{tikzpicture}% End sketch output
}}}%
\end{figure}
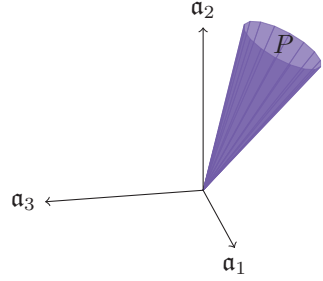%

A vector space \vektorraum{} over an ordered field\footnote{The
  additive group $\Gruppe K+-0$ is an ordered group such that $0\leq
  1$, as well as the multiplicative group of the positive elements
  $\Gruppe{\poskegel K\setminus \menge{0}}\cdot\multinverses1$.}
$\koerper=\tupel{K,+,-,0,\cdot,\multinverses,1}$ is called
\defindex{ordered vector space} iff \vektorraum\ is an ordered group
w.r.t.\ addition and for all vectors $\vektor v \in
\poskegel\vektorraum$, and all positive elements
$\alpha\in\poskegel\koerper$ of the field the condition $\alpha\vektor
v\in\poskegel\vektorraum$ is met. Consequently, the set
$\poskegel\vektorraum$ is also a positive cone in the geometrical
sense (cf.\ Fig.\ \ref{fig:poskegel-im-Vektorraum}).\footnote{The
  complete argumentation be found in the proof of Theorem
  \ref{lemma:Z-vektorraum}.}  If $A$ is a set of vectors of a vector
space, $\linhuelle A$ denotes their linear closure. The set of linear
combinations of vectors from $A$ with exclusively non-negative
coefficients will be symbolised by $\linhuelle^+A$.  As usual a subset
$A\subseteq\vektorraum$ is called linear subspace $A\leq\vektorraum$
if $\linhuelle A= A$ is true.  A subspace $\hyperebene$ of a vector
space \vektorraum\ is called \defindex{hyperplane} if there exists a
vector $\vektor a\in\vektorraum$ such that the conditions
$\linhuelle\hyperebene=\hyperebene\lneq \vektorraum$ and
$\linhuelle(\hyperebene\cup\menge{\vektor a}) =\vektorraum$ hold.  For
some base $\basis$ of the vector space $\vektorraum$ let $S⊆\basis$
be a subbase. Then for every vector $\vektor x$ the sum
\begin{equation}
  \projektion\basis S  := \sum_{\vektor s∈S}\vektor x_{\basis}(\vektor s)\vektor s
\end{equation}
is called \defindex{(coordinate) projection of $\vektor x$ into
  $\linhuelle S$ along $\basis$}.

% Als Hyperebene $\hyperebene$ bezeichnen wir
% lineare Unterräume mit Kodimension $1$. Das heißt es existiert ein
% linearer Unterraum $\reellezahlen\vektor x$ der Dimension $1$ besitzt,
% so dass $\direktsumme M\reellezahlen=\hyperebene \vektorraumplus
% \reellezahlen\vektor x$ gilt. 

In Section \ref{sec:linear-order-in-real-vector-space} we need the
notion of a neighbourhood. Since vectors in the direct sum have finite
support, we can use the standard scalar product in the resulting
vector space for any ordered field $\koerper$. To be consistent with
the following definitions, we use a generalised version: Let $\basis$
be a base of the vector space $\direktsumme{E}\koerper$, $\vektor
x,\vektor y\in\direktsumme E\koerper$ two vectors, and
$\abbildung{r}{E}{\poskegel{\koerper}\setminus\menge\nullvektor}$ a
positive mapping. Then
\begin{equation}
\skalprod{\vektor x}{\vektor y}_{r,\basis} := \sum_{\vektor b∈\basis}r(\vektor b)\cdot\vektor x_{\basis}(\vektor b)\vektor y_{\basis}(\vektor b)
\end{equation}
denotes a scalar product, as $\vektor x$ and $\vektor y$ have a finite
support, and have values that differ from zero only at finitely many
entries. For the same reason given any basis \basis\ and any positive integer $n∈\natzahlen$
\begin{equation}
  \label{eq:norm}
  \|\vektor x\|_{n,r,\basis}:=\sqrt[n]{\sum_{\vektor b\in \basis}r(\vektor b)\cdot|\vektor x_{\basis}(\vektor b)|^n}\quad\text{and}\quad
  \|\vektor x\|_{\infty,r,\basis}:=\max_{\vektor b\in \basis}r(\vektor b)\cdot|\vektor x_{\basis}(\vektor b)|
\end{equation}
are norms, and 
\begin{equation}
  \label{eq:metrik}
  \Delta_{n,r,\basis}(\vektor x,\vektor y):=\|\vektor y-\vektor x\|_{n,r,\basis}\quad\text{and}\quad
  \Delta_{\infty,r,\basis}(\vektor x,\vektor y):=\|\vektor y-\vektor x\|_{\infty,r,\basis}
\end{equation}
are metrics on $\direktsumme E\koerper$. For the special case $n=2$ we
have the usual identity $\sqrt{\skalprod {\vektor x}{\vektor
    x}}_{r,\basis}=\|\vektor x\|_{2,r,\basis}$.

For $\vektor x\in\direktsumme E\koerper$ and
$r\in\koerper$ we use the following definition for an
\defindex{open ball} $B_{r,n,\basis}(\vektor x)$ (including $n=\infty$):
\begin{equation}
B_{r,n,\basis}(\vektor x) := \Menge{\vektor y\in \direktsumme
  E\koerper}{\Delta_{n,\frac{1}{r},\basis}(\vektor x,\vektor y)<1}.
\end{equation}
%Mit Hilfe dieses Mengensystems wird $\direktsumme M\koerper$ zu
%einem topologischen Raum.
As the generated topologies of these open balls are not necessarily
the same, we define a finer topology that includes all of them: 

\begin{definition}
  A subset $M⊆\vektorraum$ of a vector space $\vektorraum$ is called
  \defindex{open}, iff for each pair of vectors $\vektor a ∈ M,
  \vektor b∈\vektorraum$ the set
  \[
  \Menge{λ∈\reellezahlen}{λ\vektor a+(1-λ)\vektor b∈M}
  \]
  is open in the set of real numbers $\reellezahlen$ with the standard
  topology.
\end{definition}

As all vectors have finite support, this defines a topology on
$\vektorraum$ where all open balls are open sets. The
\defindex{boundary} of a set $M$ with respect to this topology is
denoted by $\rand M$, and its \defindex{interior} by $\inneres M$.
%
\iffalse
 of a subset
$N\subseteq\direktsumme{E}\koerper$ according to the formula
\begin{align}
\rand N &:= \Menge[\big]{\vektor x \in\direktsumme{E}\koerper}{\forall
  r\in\poskegel\koerper\setminus\menge 0:B_r(\vektor x)\cap N
  \neq \emptyset \text{ and } B_r(\vektor x)\setminus N\neq \emptyset}
\intertext{and the \defindex{interior} of $N$ by}
\inneres N&:=N\setminus \rand N.
\end{align}
\fi

For two subsets 
$A,B\subseteq\gruppe$ of a group \gruppe\ we use the usual definition of
a complex sum
\begin{equation}
A+B:=\Menge {a+b}{a\in A, b\in B}.
\end{equation}
If $A$ and $B$ are linear subspaces of a vector space whose
intersection $A\cap B=\menge\nullvektor$ is the singleton containing the zero
vector $\nullvektor$, then we denote this by 
$A\vektorraumplus B$. Thus, for a hyperplane
$\hyperebene\leq\vektorraum$ there exists 
a vector $\vektor a$ such that
$\vektorraum = \hyperebene
\vektorraumplus\koerper\vektor a$.

\section{Preparation of injective groups}

In this section we mainly resemble Walkers theorem \cite{0085.25502}
for ordered groups providing a construction of the embedding into the
injective abelian Groups with the formalism used throughout the
current article.

We will characterise all linear order extensions of an abelian
group \gruppe{}, in this article using the vector space
$\direktsumme E\reellezahlen$ for a set $E$, that is constructed depending on
\gruppe{}. In order to achieve this we will use properties from the
integer module $\direktsumme{E}\ganzzahlen$ and its divisible
extension, the rational vector space $\direktsumme E\rationalezahlen$
as the latter is the smallest injective group containing
$\direktsumme{E}\ganzzahlen$. While the first one is isomorphic to a
subgroup of $\gruppe$, the latter one is isomorphic to a supergroup of
the divisible hull of \gruppe. As it is useful to have a well-defined
base in the considered vector spaces, we shortly discuss the necessary
steps for the construction of the mentioned algebraic structures on
the basis of the given torsion free abelian group. Furthermore, in
this section we discuss a method to transfer orders between these
structures.\footnote{It proved to be useful, to favour monomorphisms
  over subgroup relations, here.} Let's start with the following
well-known fact:

\begin{lemma}[Folklore]\label{lemma:monomorphismusordnung}%
    Let $\gruppe = \OrdGruppe G\cdot\multinverses e\leq $ be a po-group, $\gruppe'=\Gruppe
    {G'}{\cdot'}{{\multinverses}'}{e'}$ a group and $\varphi:\gruppe'\to \gruppe$ an
    injective homomorphism. Then $\gruppe'$ forms a po-group together
    with the relation
    $\sqsubseteq$ defined by 
    \begin{equation}
      g\sqsubseteq h:\Leftrightarrow\varphi(g)\leq\varphi(h). \label{eq:monomorphismusordnung}
    \end{equation}
 \end{lemma}
\begin{proof}
  The mapping $φ$ is a group isomorphism from  $G'$ onto $φ[G']$ and
  according to equation \eqref{eq:monomorphismusordnung} it is also an
  order isomorphism. So it remains to show that the group $\gruppe'$
  is a partially ordered group together with the relation
  $\sqsubseteq$. Let $a,b,x,y\in G'$ elements of the group
  $\gruppe'$. Then the following equivalences hold:
  \begin{align*}
     a\sqsubseteq b&\Leftrightarrow φ(a)\leq φ(b)\Leftrightarrow
     \underbrace{φ(x)\cdot φ(a)\cdot φ(y)}_{{}=φ(x\cdot'a\cdot'y)}\leq\underbrace{φ(x) \cdot φ(b) \cdot φ(y)}_{{}=φ(x \cdot'b\cdot'y)}\\
     &\Leftrightarrowφ(x\cdot'a\cdot'y) \leqφ(x\cdot'b\cdot'y)\\
&\Leftrightarrow x\cdot'a\cdot'y\sqsubseteq x\cdot'b\cdot'y.
   \end{align*}
   So $φ$ is an o-homomorphism from $\gruppe'$ into \gruppe.
\end{proof}

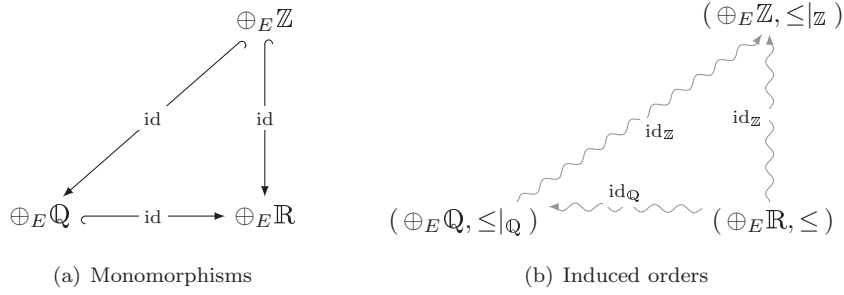
\begin{figure}\par\noindent%
  \parbox{0.45\linewidth}{\subfigure[][Monomorphisms]{\label{subfig:zqr-monomorphismen}%
      \begin{tikzpicture}[description/.style={fill=white,
          inner sep=2pt},>=latex,x=0.9em,y=1em]%
        \useasboundingbox (-8.5,-5) rectangle (8.5,5);%
        \matrix (m) [matrix of math nodes, row sep = 6em, 
        column sep=5.5em,text height=1.5ex,text depth=0.25ex,ampersand replacement=\&]%
        {\& \direktsumme E\ganzzahlen\\
          \direktsumme E\rationalezahlen \&\direktsumme
          E\reellezahlen \\};%
        \path[left hook->,font=\scriptsize] (m-2-1) edge node[description] {$\id$} (m-2-2);%
        \path[right hook->,font=\scriptsize] (m-1-2) edge node[description] {$\id$} (m-2-2);%
        \path[right hook->,font=\scriptsize] (m-1-2) edge node[description] {$\id$} (m-2-1);%
      \end{tikzpicture}%
    }%
  }%
  \parbox{0.55\linewidth}{%
  \subfigure[][Induced orders]{\label{subfig:zqr-ordnungen-einseitig}%
    \begin{tikzpicture}[x=1em,y=1em,description/.style={inner
         sep=2pt,black,auto,fill=white,swap},>=latex,decoration={snake,
         amplitude=0.15em,
         foot of=human,
         foot length=1.5em,
         stride length=4em,foot sep=-0.7em}]%
       \useasboundingbox (-9,-5) rectangle (9,5);%
       \matrix (m) [matrix of math nodes, row sep = 6.4em, 
       column sep=5.5em,text height=1.5ex,text depth=0.25ex,ampersand replacement=\&]%
      { \&\paar{\direktsumme E\ganzzahlen}{{\leq}{|_\ganzzahlen}}\\
        \paar{\direktsumme E\rationalezahlen}{{\leq}|_\rationalezahlen} \&  \paar{\direktsumme E\reellezahlen}{\leq} \\};
      \path[->,font=\scriptsize,gray]    (m-2-2.170) edge[decorate] node[description] {$\id_\rationalezahlen$} (m-2-1.10);%
      \path[->,font=\scriptsize,gray,swap=false]    (m-2-2.100) edge[decorate] node[description] {$\id_\ganzzahlen$} (m-1-2.-100);%
      \path[->,font=\scriptsize,gray]    (m-2-1.20) edge[decorate] node[description] {$\id_\ganzzahlen$} (m-1-2.-120);%
     \end{tikzpicture}%
   }%
  }%
  \caption{Monomorphisms and induced orders according to Lemma
    \ref{lemma:monomorphismusordnung}}%
  \label{fig:monomorphismusordnung-Z-Q-G}
\end{figure}%

This lemma allows us to relate possible orders on the three
structures $\direktsumme E\ganzzahlen$, $\direktsumme
E\rationalezahlen$ and $\direktsumme
E\reellezahlen$. Figure~\ref{fig:monomorphismusordnung-Z-Q-G} shows
these relationships with respect to the identical embeddings. In the
left hand Figure \ref{subfig:zqr-monomorphismen} the diagram of the
embeddings is shown, while the arrows on the right hand side
\ref{subfig:zqr-ordnungen-einseitig} show the possible ways of
induction or transfer of the orders between the structures.

In order to add an arbitrary torsion free abelian group into this
system of monomorphisms and induced orders, suppose $E$ is an
independent set of the group \gruppe.  Thus, considering cyclic groups, a
positive integer $n\in\natzahlen\setminus\menge 0$ exists for each
element $g\in G$ such that $ng\in\erzgruppe E$. Thereby, we found a
more or less unique representation for each Element of \gruppe{} by
elements of a maximal independent set $E$:

\begin{lemma}\label{lemma:defnatzahl}\label{lemma:defnatzahl-eindeutig}
  Let\/ \gruppe\ be a torsion free abelian group and $E$ a maximal
 independent set in\/ \gruppe{}. Then for each element $g\in G$ a positive integer 
$q_g\in\natzahlen\setminus\menge 0$ and a mapping $\abbildung{\vektor
  p_g}{E}{\ganzzahlen}$ 
 exist such that the following equation holds:
 \begin{equation}
    \label{eq:defnatzahl}
    q_gg=\sum_{a\in E}\vektor p_g(a)a.
  \end{equation}
  The factors $q_g$ and $\vektor p_g(a)$ are uniquely defined up to
  multiplication by a common rational number.
\end{lemma}

\begin{proof}As $\erzgruppe{E}$ is an essential subgroup, for each
  element $g\in G$ the
  intersection $\erzgruppe{E}\cap\erzgruppe{g}\neq\emptyset$ is
  non-empty. This implies that a
  positive integer $q_g$ exists such that $q_gg\in\erzgruppe{E}$. This proves
  the existence of a representation of the form \eqref{eq:defnatzahl}.

 Let
\[ q_gg=\sum_{a\in E}\vektor p_g(a)a\text{ and } q_gg =\sum _{a\in
  E}\vektor p'_g(a)a\]
two such representations according to \eqref{eq:defnatzahl}. Considering
their difference we get
\[0=(q_g-q_g)g=\sum_{a\in E}(\vektor p_g-\vektor p'_g)(a) a.\]
Since $E$ is independent, we deduce from this equation that for each
$a\in E$ the condition $\vektor p_g(a) - \vektor p'_g(a) = 0$ holds. So we
proved the identity $\vektor p_g =\vektor p'_g$ for any fixed number $q_g$.

Let's assume that we have different representations of the form
\eqref{eq:defnatzahl}. Since the mapping $\vektor p_g$ is unique for
each $q_g$, there must be different numbers on the left hand side
too. Let $q_g\in\natzahlen\setminus\menge0$ the smallest integer
for which a representation in the form \eqref{eq:defnatzahl} exists.
Furthermore, let us assume there exists a number 
$q'_g\in\natzahlen\setminus\menge{0,1}$, which is coprime to $q_g$
such that:
\[q'_gg=\sum_{a\in E}\vektor p'_g(a)a.\]
Then the following two equations hold:
\begin{align*}
  (q_g-q'_g)g&=\sum_{a\in E} (\vektor p_g-\vektor p'_g)(a)a\text{ and}\\
  (q_g+q'_g)g&=\sum_{a\in E} (\vektor p_g+\vektor p'_g)(a)a.
\end{align*}
This allows us to use the Euclidian algorithm to find two integers $c$
and $d$ such that $cq_g+dq'_g=\ggT\paar {q_g}{q'_g}$ holds. Using
these the equation
\[
\ggT\paar {q_g}{q'_g}g = (cq_g+dq'_g)g = \sum_{a\in E}(c\vektor
p_g+d\vektor p'_g)(a)a. (\widerspruch)
\]
holds. Since $0<\ggT\paar {q_g}{q'_g}\leq q_g $ and $q_g$ was chosen
minimal and $q'_g$ coprime to $q_g$, there exists a contradiction. So
$q'_g$ is a multiple of $q_g$. 

Let $q'_g = bq_g$ for some positive integer
$b\in\natzahlen\setminus\menge{0}$. Then we have
\[ 
\sum_{a\in E}b\vektor p_g (a)a =b q_gg=q'_gg = \sum_{a\in E}\vektor p'_g (a)a.
\]
This proves that both $q'_g$ and $\vektor{p'_g}$ are multiples with
the same factor of $q_g$ and $\vektor{p_g}$, respectively.
\end{proof}

\noindent Using this Lemma~\ref{lemma:defnatzahl-eindeutig} we can
find an embedding of $\gruppe$ with a maximal independent set $E$ into the direct sum
$\direktsumme{E}\rationalezahlen$ (cf. Walker's theorem \cite{0085.25502}): 

\begin{theorem}\label{lemma:G-Q-Monomorphismus}\label{satz:Monomorphismus-G-Q}
  Let $\gruppe = \Gruppe G+-0$ be a torsion free abelian group and
  $E\subseteq G$ a maximal independent set in \gruppe. 
  Let for each $g\in G$ a mapping  $\abbildung{p_g}{E}{\ganzzahlen}$ and
   a positive integer $q_g\in\natzahlen\setminus\menge 0$ defined according to
  equation \eqref{eq:defnatzahl}. Then the mapping 
  \begin{equation}
    \label{eq:homomorphismus-G-Q}
    \wabbildung{\abbildung{φ}{G}{\direktsumme
      E\rationalezahlen}}g{\frac{\vektor p_g}{q_g}}
\end{equation}
is a monomorphism from 
    \gruppe\ into $\direktsumme E\rationalezahlen$.
\end{theorem}

\begin{proof}
  Lemma \ref{lemma:defnatzahl} ensures  for each $g\in G$ the
  existence of integers
  $q_g, \vektor p_g(a)$ (for all $a\in E$) that fulfil the following equation:
  \[
  q_gg=\sum_{a\in E}\vektor p_g(a)a
  \]
  Let 
  $\wabbildung{\abbildung\varphi G{\direktsumme
      E\rationalezahlen}}g{\frac {\vektor p_g}{q_g}}$ be the mapping
  as defined in
  equation \eqref{eq:homomorphismus-G-Q}. Then for each
  $c\in\ganzzahlen$ the equation
  $\varphi(g)=\frac{c\vektor p_g}{cq_g}=\frac{\vektor p_g}{q_g}$
  holds. From this we infer
  together with Lemma \ref{lemma:defnatzahl-eindeutig} that $\varphi$
  is well-defined.

  Furthermore, having
  $\varphi(cg)=\frac {c\vektor p_g}{q_g}=c\varphi(g)$ we can combine the two representations
 \[
  q_gg=\sum_{a\in E}\vektor p_g(a)a, \qquad q_hh=\sum_{a\in E}\vektor p_h(a)a
  \]
  of elements $g,h\in G$, into one common formula: 
 \[
  q_gq_h(g+h)=q_h\sum_{a\in E}\vektor p_g(a)a + q_g\sum_{a\in E}\vektor p_h(a)a.
  \]
  Thus, the mapping $\varphi$ is a group homomorphism from \gruppe\
  into $\direktsumme E\rationalezahlen$ with respect to the addition.
  Consequently,
  \[
  \varphi(g+h)=\frac{q_h\vektor p_g}{q_gq_h} + \frac{q_g\vektor p_h}{q_gq_h} =
  \varphi(g) + \varphi(h)
  \]

Let $\varphi(g)=\nullvektor$. Then for each $a\in E$ the element
  $\vektor p_g(a)=0$ is zero, because $E$ is an independent set. Thus, for
  any $q_g\in\natzahlen\setminus\menge{0}$ the equality $q_gg=0$ holds, which
  implies $g=0$. So $φ$ is injective, which means it is a monomorphism
  from $\gruppe$ into $\direktsumme E\rationalezahlen$.
\end{proof}

\begin{corollary}
  The homomorphism  $φ$ maps each element $a\in E$ of the independent
  set  onto a unit vector, i.e., 
  $\varphi(a) = \vektor e_a$.
\end{corollary}

With Theorem \ref{satz:Monomorphismus-G-Q} each order in the real (rational)
vector space $\direktsumme{E}\reellezahlen$
($\direktsumme{E}\rationalezahlen$) defines an order on the group $\gruppe$.

\begin{theorem}\label{lemma:z-g-monomorphismus}\label{lemma:Monomorphismus-Z-G}%
  Let $\gruppe$ be a torsion free abelian group and $E$ a maximal independent
  set in \gruppe{}. Furthermore, let
  $\abbildung{φ}{\gruppe}{\direktsumme{E}\rationalezahlen}$
    be defined as in equation \eqref{eq:homomorphismus-G-Q}. Then the
    mapping 
    \begin{equation}
      \label{eq:z-g-mornomorphismus}
      \wabbildung{\abbildung{ψ}{\direktsumme{E}\ganzzahlen}{\gruppe}}{\vektor{x}}{\sum_{a\in
        E}\vektor{x}(a)a}
    \end{equation}
    is a monomorphism such that for all $g\in \erzgruppe E$ and all
    $\vektor{x}\in\direktsumme{E}\ganzzahlen$ the
    conditions
  \begin{equation}\label{eq:phi-und-psi-kommutieren-auf-Z}
    ψ\bigl(φ(g)\bigr)=g\qquad\text{and}\qquadφ\bigl(ψ(\vektor x)\bigr) = \vektor x
  \end{equation}
  hold.
\end{theorem}

\begin{proof}
  Firstly we show that $ψ$ is a homomorphism. Since \gruppe\ is
  commutative, we can rewrite the sum in the
  following way:
 \begin{equation*}
    ψ(\vektor x + \vektor y) = \sum _{a\in E}(\vektor x+\vektor
    y)(a)a = \sum_{a\in E} \vektor x(a) a + \sum_{a\in E}\vektor y(a)a
    = ψ(\vektor x) + ψ(\vektor y).
  \end{equation*}
  The last identity holds, since \gruppe{} is torsion free.
  Furthermore, $E$ is an independent set in \gruppe{} and the equation
  $\vektor e_a(a) = 1$ holds. Thus, $\sum_{a\in E} \vektor x(a)a = 0$
  implies $\vektor x = \nullvektor$, which proves that $ψ$ is an
  injective homomorphism.

  Since $φ$ and $ψ$ are monomorphisms and $\erzgruppe{E}$ is a
  subgroup in \gruppe{}, it suffices to show the equations
  \eqref{eq:phi-und-psi-kommutieren-auf-Z} for the elements of the
  generating sets $E$ and $\Menge {\vektor e_a}{a\in
    E}$. For those it is easy to show the following identities:
  \begin{align*}
    ψ\bigl(φ(a)\bigr)&=ψ(\vektor e_a) = a,\text{ and}\\
    φ\bigl(ψ(\vektor e_a)\bigr)&=φ(a) = \vektor e_a.
  \end{align*}
  Since these identities hold for any element of the corresponding set
  they are satisfied for any element of the corresponding (sub)group, too.
\end{proof}

\begin{figure}\par\noindent%
  \parbox{0.45\linewidth}{\subfigure[][Monomorphisms]{\label{subfig:zqrg-monomorphismen}%
 \begin{tikzpicture}[description/.style={fill=white,
     inner sep=2pt},>=latex,x=0.9em,y=1em]%
  \useasboundingbox (-9,-5) rectangle (9,5);%
  \matrix (m) [matrix of math nodes, row sep = 6em, 
    column sep=5.5em,text height=1.5ex,text depth=0.25ex,ampersand replacement=\&]%
    {\gruppe\&\direktsumme E\ganzzahlen\\
      \direktsumme E\rationalezahlen   \& \direktsumme E\reellezahlen \\};%
  %
%  \only<1,6,8->{
  \path[right hook->,font=\scriptsize,alert6] (m-1-2) edge node[auto] {$\psi$}  (m-1-1);%
 % \only<1,5,6,7->{
  \path[right hook->,font=\scriptsize,alert5] (m-1-1) edge node[auto] {$\varphi$} (m-2-1);%
%  \only<1,3->{
  \path[left hook->,font=\scriptsize,alert3] (m-2-1) edge node[description] {$\id$} (m-2-2);%
 % \only<1,3->{
  \path[right hook->,font=\scriptsize,alert3] (m-1-2) edge node[description] {$\id$} (m-2-2);%
 % \only<1,3->{
  \path[right hook->,font=\scriptsize,alert3] (m-1-2) edge node[description] {$\id$} (m-2-1);%
 \end{tikzpicture}%
%\caption{Monomorphisms according to Lemma \ref{lemma:z-gmonomorphismus}}%
}%
}%
\parbox{0.55\linewidth}{%
  \subfigure[][Induced orders]{\label{subfig:zqrg-ordnungen-einseitig}%
  \begin{tikzpicture}[x=1em,y=1em,description/.style={inner
       sep=2pt,black,auto,fill=white,swap},>=latex,decoration={snake,
       amplitude=0.15em,
       foot of=human,
       foot length=1.5em,
       stride length=4em,foot sep=-0.7em},ampersand replacement=\&]%
     \useasboundingbox (-9,-5) rectangle (9,5);%
     \matrix (m) [matrix of math nodes, row sep = 6.4em, 
    column sep=5.5em,text height=1.5ex,text depth=0.25ex]%
    {  \paar\gruppe\sqsubseteq \&\paar{\direktsumme E\ganzzahlen}{{\leq}|_\ganzzahlen}\\
      \paar{\direktsumme E\rationalezahlen}{{\leq}|_\rationalezahlen} \&  \paar{\direktsumme E\reellezahlen}{\leq} \\};
    % \only<1,6,8->{
    \path[<-,font=\scriptsize,gray,alert6]    (m-1-2) edge[decorate] node[description] {$\psi$}  (m-1-1);%
    % \only<1,5,6,7->{
    \path[<-,font=\scriptsize,gray,alert5]   (m-1-1) edge[decorate] node[description] {$\varphi$} (m-2-1) ;%
    % \only<1,4,5,6,9->{\path[->,font=\scriptsize,gray,alert4]    (m-2-1.-10) edge[decorate] node[description] {$\conv_\reellezahlen$} (m-2-2.-170);}%
    % \only<1,3->{
    \path[->,font=\scriptsize,gray,alert3]    (m-2-2.170) edge[decorate] node[description] {$\id_\rationalezahlen$} (m-2-1.10);%
    % \only<1,4,5,6,9->{\path[->,font=\scriptsize,gray,alert4,swap=false]    (m-1-2.-40) edge[decorate] node[description] {$\conv_\reellezahlen$} (m-2-2.40);}%
    % \only<1,3->{
    \path[->,font=\scriptsize,gray,alert3,swap=false]    (m-2-2.100) edge[decorate] node[description] {$\id_\ganzzahlen$} (m-1-2.-100);%
    % \only<1,4,5,6,9->{\path[->,font=\scriptsize,gray,alert4]   (m-1-2.-161) edge[decorate] node[description] {$\conv_\rationalezahlen$} (m-2-1.51);}%
    % \only<1,3->{
    \path[->,font=\scriptsize,gray,alert3]    (m-2-1.20) edge[decorate] node[description] {$\id_\ganzzahlen$} (m-1-2.-120);%
  \end{tikzpicture}%
%  \caption{Induced orders according to Lemma \ref{lemma:z-g-monomorphismus}}%
  }%
  }%
  \caption{Monomorphisms and induced orders according to Lemmata
    \ref{lemma:monomorphismusordnung}, \ref{satz:Monomorphismus-G-Q}
    and \ref{lemma:Monomorphismus-Z-G}}%
  \label{fig:monomorphismusordnung-Z-Q-R-G}
\end{figure}
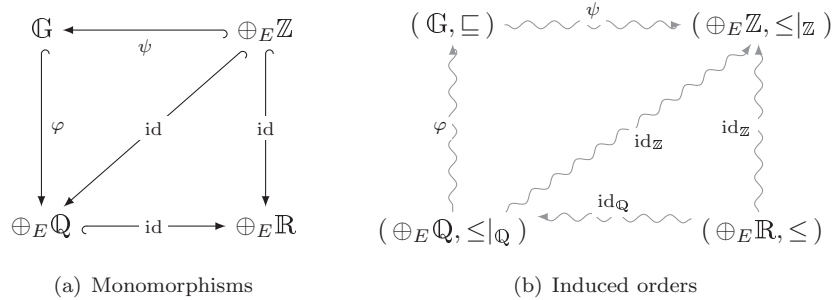%

The preceding two theorems prove that the diagram in Figure
\ref{subfig:zqrg-monomorphismen} commutes. Thus, we can transfer
orders along the arrows in Figure
\ref{subfig:zqrg-ordnungen-einseitig} between the different
groups. Together with Theorem \ref{satz:charakterisierung-linearer-ordnungen},
this would be sufficient to rework the theorem of {\ru М.\,И. Зайцева}
(M.\,I. Zajceva) \cite{M.I.Zajceva:1953} in the language of vector
spaces. Since we want to discuss order extensions, we must firstly
transfer the order from a given group into the corresponding real
vector space. The following section will provide this link.

\section{Copying the order into the vector space}

Any order in the vector space constructed in the last section implies
an order in our group
(cf. Figure~\ref{subfig:zqr-ordnungen-einseitig}). So far the vector
space has no idea about the existing order in the group
$\gruppe$. This section closes this gap. There are some well-known
theorems such as Hahn's theorem \cite{Hahn:1907} or the
Conrad-Harvey-Holland-Theorem (cf. \cite{glass1999partially},
chap. 4.6) which cover the embedding of ordered abelian groups into
real vector spaces. These theorems provide embeddings for each order,
but they do not assure that it is possible to use the same embedding
for all orders.

\iffalse
In order to characterise linear order extensions of a partially
ordered torsion free abelian group using a real vector space, we need a
notion of the order on the group in the vector space.
\fi
We can transfer the order from the group onto an integer module
(see Figure \ref{fig:monomorphismusordnung-Z-Q-R-G}). However, going this
path we loose information about group elements. On the other hand, the
integer module $\direktsumme{E}\ganzzahlen$ is a set of size nearly zero in the
real vector space $\direktsumme{E}\reellezahlen$, so it is not easy to tell which linear vector space
orders correspond to one linear order in the integer module.

As the
notion convex is already order-theoretically defined in the language
of partially ordered groups, we choose a different wording, here:

\begin{definition}[cf.\ {\cite{r1995theory}, Def.~3.2}]
  Let $\gruppe=\OrdGruppe G+-0\leq$ a partially ordered abelian
  group. The order $\leq$ is called \alert{semiclosed} iff for all
  $a\in G$ and $n\in\natzahlen\setminus\menge 0$ the following
  implication holds:
  \begin{equation}
    0\leq na\Rightarrow 0\leq a.
  \end{equation}
\end{definition}

We show in this section that semiclosed orders on any abelian
torsion free group can be easily extended to vector space
orders. This will be done by forming the convex hull. Doing so, we
prove that in this particular case the induced order on the group will
be the same as the original one. Fortunately, the class of semiclosed
ordered torsion free abelian groups is the relevant class of groups to
consider, here. Each partially ordered abelian group has a canonical semiclosed
order extension, which itself is a suborder of any linear order
extension of the group. So the restriction to semiclosed groups does
not influence the generality of this construction. 

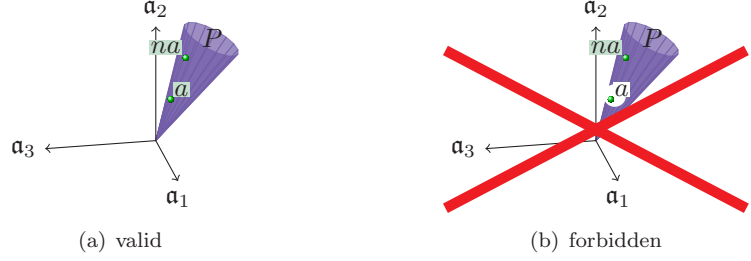
\begin{figure}
  \parbox{0.5\linewidth}{\subfigure[][valid]{\label{subfig:zqr-halb-abgeschlossen-erlaubt}%
     \centerline{% Sketch output, version 0.3 (build 7d, Fri Nov 22 14:14:09 2013)
% Output language: PGF/TikZ,LaTeX
\begin{tikzpicture}[line join=round]
\tikzstyle{ann} = [fill=white,font=\footnotesize,inner sep=1pt]
\providecommand\vektor{}\tikzstyle{plane}=[shading=axis,top color=white,bottom
color=HKS41K80,shading angle=-50]\tikzstyle{plane2}=[shading=axis,top color=white,bottom
color=HKS41K80,shading angle=-30]\filldraw[draw=HKS36K80,draw opacity=0.7,line width=0.0001,fill=HKS36K80,fill opacity=0.7](0,0)--(.383,1.545)--(.437,1.573)--(0,0)--cycle;
\filldraw[draw=HKS36K80,draw opacity=0.7,line width=0.0001,fill=HKS36K80,fill opacity=0.7](0,0)--(.384,1.486)--(.383,1.545)--(0,0)--cycle;
\filldraw[draw=HKS36K80,draw opacity=0.7,line width=0.0001,fill=HKS36K80,fill opacity=0.7](0,0)--(.437,1.573)--(.536,1.567)--(0,0)--cycle;
\filldraw[draw=HKS36K80,draw opacity=0.7,line width=0.0001,fill=HKS36K80,fill opacity=0.7](0,0)--(.439,1.407)--(.384,1.486)--(0,0)--cycle;
\filldraw[draw=HKS36K80,draw opacity=0.7,line width=0.0001,fill=HKS36K80,fill opacity=0.7](0,0)--(.536,1.567)--(.666,1.528)--(0,0)--cycle;
\filldraw[draw=HKS36K80,draw opacity=0.7,line width=0.0001,fill=HKS36K80,fill opacity=0.7](0,0)--(.539,1.319)--(.439,1.407)--(0,0)--cycle;
\filldraw[draw=HKS36K80,draw opacity=0.7,line width=0.0001,fill=HKS36K80,fill opacity=0.7](0,0)--(.666,1.528)--(.807,1.461)--(0,0)--cycle;
\filldraw[draw=HKS36K80,draw opacity=0.7,line width=0.0001,fill=HKS36K80,fill opacity=0.7](0,0)--(.67,1.235)--(.539,1.319)--(0,0)--cycle;
\filldraw[draw=HKS36K80,draw opacity=0.7,line width=0.0001,fill=HKS36K80,fill opacity=0.7](0,0)--(.937,1.377)--(1.038,1.289)--(0,0)--cycle;
\filldraw[draw=HKS36K80,draw opacity=0.7,line width=0.0001,fill=HKS36K80,fill opacity=0.7](0,0)--(.807,1.461)--(.937,1.377)--(0,0)--cycle;
\filldraw[draw=HKS36K80,draw opacity=0.7,line width=0.0001,fill=HKS36K80,fill opacity=0.7](0,0)--(.811,1.168)--(.67,1.235)--(0,0)--cycle;
\filldraw[draw=HKS36K80,draw opacity=0.7,line width=0.0001,fill=HKS36K80,fill opacity=0.7](0,0)--(1.038,1.289)--(1.093,1.21)--(0,0)--cycle;
\filldraw[draw=HKS36K80,draw opacity=0.7,line width=0.0001,fill=HKS36K80,fill opacity=0.7](0,0)--(.941,1.129)--(.811,1.168)--(0,0)--cycle;
\filldraw[draw=HKS36K80,draw opacity=0.7,line width=0.0001,fill=HKS36K80,fill opacity=0.7](0,0)--(1.093,1.21)--(1.093,1.152)--(0,0)--cycle;
\filldraw[draw=HKS36K80,draw opacity=0.7,line width=0.0001,fill=HKS36K80,fill opacity=0.7](0,0)--(1.04,1.123)--(.941,1.129)--(0,0)--cycle;
\filldraw[draw=HKS36K80,draw opacity=0.7,line width=0.0001,fill=HKS36K80,fill opacity=0.7](0,0)--(1.093,1.152)--(1.04,1.123)--(0,0)--cycle;
\draw[arrows=<->,line width=.4pt](.294,-.537)--(0,0)--(0,1.517);
\draw[arrows=->,line width=.4pt](0,0)--(-1.47,-.107);
\path (.748,1.365) node {$P$};\path (.294,-.537) node[below] {$\vektor a_1$}
     (0,1.517) node[above] {$\vektor a_2$}
     (-1.47,-.107) node[left] {$\vektor a_3$};\fill[fill=white] (.25,.6) node (n0){};
\shade[shading=ball,ball color=HKS57K100] (n0) +(-0.05,-0.05) circle (0.04);
\draw (n0) +(-0.05,-0.05) node[anchor=south west,inner sep=00.02cm,outer sep=0.03cm,fill=HKS57K20] {$a$};\fill[fill=white] (.5,1.2) node (n1){};
\shade[shading=ball,ball color=HKS57K100] (n1) +(-0.1,-0.1) circle (0.04);
\draw (n1) +(-0.1,-0.1) node[anchor=south east,inner sep=0.02cm,outer sep=00.05cm,fill=HKS57K20] {$na$};\end{tikzpicture}% End sketch output
}}}%
  \parbox{0.5\linewidth}{\subfigure[][forbidden]{\label{subfig:zqr-halb-abgeschlossen-verboten}%
      \centerline{% Sketch output, version 0.3 (build 7d, Fri Nov 22 14:14:09 2013)
% Output language: PGF/TikZ,LaTeX
\begin{tikzpicture}[line join=round]
\tikzstyle{ann} = [fill=white,font=\footnotesize,inner sep=1pt]
\providecommand\vektor{}\tikzstyle{plane}=[shading=axis,top color=white,bottom
color=HKS41K80,shading angle=-50]\tikzstyle{plane2}=[shading=axis,top color=white,bottom
color=HKS41K80,shading angle=-30]\filldraw[draw=HKS36K80,draw opacity=0.7,line width=0.0001,fill=HKS36K80,fill opacity=0.7](5,0)--(5.383,1.545)--(5.437,1.573)--(5,0)--cycle;
\filldraw[draw=HKS36K80,draw opacity=0.7,line width=0.0001,fill=HKS36K80,fill opacity=0.7](5,0)--(5.384,1.486)--(5.383,1.545)--(5,0)--cycle;
\filldraw[draw=HKS36K80,draw opacity=0.7,line width=0.0001,fill=HKS36K80,fill opacity=0.7](5,0)--(5.437,1.573)--(5.536,1.567)--(5,0)--cycle;
\filldraw[draw=HKS36K80,draw opacity=0.7,line width=0.0001,fill=HKS36K80,fill opacity=0.7](5,0)--(5.439,1.407)--(5.384,1.486)--(5,0)--cycle;
\filldraw[draw=HKS36K80,draw opacity=0.7,line width=0.0001,fill=HKS36K80,fill opacity=0.7](5,0)--(5.536,1.567)--(5.666,1.528)--(5,0)--cycle;
\filldraw[draw=HKS36K80,draw opacity=0.7,line width=0.0001,fill=HKS36K80,fill opacity=0.7](5,0)--(5.539,1.319)--(5.439,1.407)--(5,0)--cycle;
\filldraw[draw=HKS36K80,draw opacity=0.7,line width=0.0001,fill=HKS36K80,fill opacity=0.7](5,0)--(5.666,1.528)--(5.807,1.461)--(5,0)--cycle;
\filldraw[draw=HKS36K80,draw opacity=0.7,line width=0.0001,fill=HKS36K80,fill opacity=0.7](5,0)--(5.67,1.235)--(5.539,1.319)--(5,0)--cycle;
\filldraw[draw=HKS36K80,draw opacity=0.7,line width=0.0001,fill=HKS36K80,fill opacity=0.7](5,0)--(5.937,1.377)--(6.038,1.289)--(5,0)--cycle;
\filldraw[draw=HKS36K80,draw opacity=0.7,line width=0.0001,fill=HKS36K80,fill opacity=0.7](5,0)--(5.807,1.461)--(5.937,1.377)--(5,0)--cycle;
\filldraw[draw=HKS36K80,draw opacity=0.7,line width=0.0001,fill=HKS36K80,fill opacity=0.7](5,0)--(5.811,1.168)--(5.67,1.235)--(5,0)--cycle;
\filldraw[draw=HKS36K80,draw opacity=0.7,line width=0.0001,fill=HKS36K80,fill opacity=0.7](5,0)--(6.038,1.289)--(6.093,1.21)--(5,0)--cycle;
\filldraw[draw=HKS36K80,draw opacity=0.7,line width=0.0001,fill=HKS36K80,fill opacity=0.7](5,0)--(5.941,1.129)--(5.811,1.168)--(5,0)--cycle;
\filldraw[draw=HKS36K80,draw opacity=0.7,line width=0.0001,fill=HKS36K80,fill opacity=0.7](5,0)--(6.093,1.21)--(6.093,1.152)--(5,0)--cycle;
\filldraw[draw=HKS36K80,draw opacity=0.7,line width=0.0001,fill=HKS36K80,fill opacity=0.7](5,0)--(6.04,1.123)--(5.941,1.129)--(5,0)--cycle;
\filldraw[draw=HKS36K80,draw opacity=0.7,line width=0.0001,fill=HKS36K80,fill opacity=0.7](5,0)--(6.093,1.152)--(6.04,1.123)--(5,0)--cycle;
\draw[arrows=<->,line width=.4pt](5.294,-.537)--(5,0)--(5,1.517);
\draw[arrows=->,line width=.4pt](5,0)--(3.53,-.107);
\draw[draw=red,line width=4](3,-.9)--(7,1.2);
\draw[draw=red,line width=4](3,1.2)--(7,-.9);
\path (5.748,1.365) node {$P$};\path (5.294,-.537) node[below] {$\vektor a_1$}
     (5,1.517) node[above] {$\vektor a_2$}
     (3.53,-.107) node[left] {$\vektor a_3$};\fill[fill=white] (5.25,.6) node (n0){} circle (0.15);
\shade[shading=ball,ball color=HKS57K100] (n0) +(-0.05,-0.05) circle (0.04);
\draw (n0) +(-0.05,-0.05) node[anchor=south west,inner sep=00.02cm,outer sep=0.03cm,fill=white] {$a$};\fill[fill=white] (5.5,1.2) node (n1){};
\shade[shading=ball,ball color=HKS57K100] (n1) +(-0.1,-0.1) circle (0.04);
\draw (n1) +(-0.1,-0.1) node[anchor=south east,inner sep=0.02cm,outer sep=00.05cm,fill=HKS57K20] {$na$};\end{tikzpicture}% End sketch output
}}}%
\caption{Semiclosed order in an ordered group}
\label{fig:semiclosed}
\end{figure}
 
As shown in Figure \ref{fig:semiclosed} the concept of a semiclosed
order can be considered as an order on the group, whose positive cone
is convex in a geometrical sense (we will discuss this later). 
In geometric  discussions about convex sets the convex hull plays an
important role. As a semiclosed order has a “convex“ positive cone,
we can ask, whether the positive cone of any partial order has  a well
defined “convex hull“. The following lemma discusses this fact:

%\frametitle{Semiclosed extension of an order}
\begin{lemma}[cf.\ {\cite{r1995theory}, Cor.~29.10}]\label{lemma:halbabschluss}
  Let $\gruppe=\OrdGruppe G+-0\leq$ be a partially ordered abelian group
  and the order relation $\leq'$ defined by
\begin{equation}
    \label{eq:ordnungserweiterung}
    0\leq' g:\Leftrightarrow \exists m \in\natzahlen\setminus\menge 0: 0\leq mg.
  \end{equation}
  Then $≤'$ is a semiclosed order relation and \OrdGruppe G+-0{\leq'}\ is a
  partially ordered group. Each linear or lattice order on \gruppe\ is
  an extension of $\leq'$.
\end{lemma}

\begin{proof}
  Firstly, we have to show that $P:=\Menge{g\in G}{0\leq' g}$ fulfils
  the conditions \eqref{eq:poskegelbedingung_halbgruppe} to
  \eqref{eq:poskegelbedingung_invariant}. Let $g,h\in P$ be two elements
  of this set. Then there exist two positive integers
  $m,n\in\natzahlen\setminus\menge 0$ such that $0\leq mg$ and $0\leq
  nh$. This implies $0\leq mng$ and $0\leq
  mnh$. Consequently, we have $0\leq mng+mnh=mn(g+h)$ and  $0\leq'
  g+h$. Since \gruppe{}
  is a commutative group, the set $P$ is a semigroup and invariant. If for any nonzero element $g\neq 0$ the inequality
  $0\leq mg$ holds, there exists no positive number $n\in \natzahlen$
  such that $ng\leq 0$.
  Otherwise from $ng\leq 0$ follows $mng \leq 0$ and from $0\leq mg$
  follows $0 \leq mng$. This would imply $g=0$, which has been excluded.\widerspruch

  Thus, with $P\cap -P=\menge 0$ the set $P$ is a pure subset of \gruppe{}. This
  shows that $\OrdGruppe
  G+-0{\leq'}$ is a partially ordered group. Since the condition
  \eqref{eq:ordnungserweiterung} is always true in lattice ordered
  groups (cf. \cite{r1995theory}, Prop.~3.6) and the positive cones of
  $\leq$ and $\leq'$ differ by exactly those elements, which
  contradict this condition, the order $\leq'$ is a suborder of any
  lattice order extension of $\leq$.  
\end{proof}

Up to now, the transitive closure of the arrows in Figure
\ref{subfig:zqrg-ordnungen-einseitig} is a partial order. Our aim is
to consider the orders on the different groups to be more or less
equivalent. The following theorem extends a given order on a subgroup
to a semiclosed order on the containing vector space. As a result of
this theorem, we will be able to invert the arrows in the lower right
triangle of this figure.

\begin{theorem}\label{lemma:Z-vektorraum}%
  Let $\direktsumme E\ganzzahlen$ an integer module and
  $\direktsumme E\reellezahlen$ a vector space of equal
  dimension. Then there exists an injective mapping, which maps each
  semiclosed order on $\direktsumme E\ganzzahlen$ to a vector space
  order on $\direktsumme E\reellezahlen$. The same is true for the
  combinations $\direktsumme E\ganzzahlen/\direktsumme
  E\rationalezahlen$ and $\direktsumme E\rationalezahlen/\direktsumme
  E\reellezahlen$.
\end{theorem}

\begin{proof} It is sufficient to prove that the positive cone $\poskegel{(\direktsumme
    E\ganzzahlen)}$ can be embedded into a convex invariant subset $P$ of the
  vector space $\direktsumme{E}\reellezahlen$ that fulfils the
  Conditions \eqref{eq:poskegelbedingung_halbgruppe} to
  \eqref{eq:poskegelbedingung_invariant}. We will show the injectivity
  afterwards. As candidate for the set $P$ we choose the convex hull
  of $φ\bigl[ψ[\poskegel{(\direktsumme
    E\ganzzahlen)}]\bigr]$ in the real vector space
  $\direktsumme{E}\reellezahlen$ using the monomorphisms $φ$ and $ψ$ defined in
  \eqref{eq:homomorphismus-G-Q} and \eqref{eq:z-g-mornomorphismus}. 
  Let 
  \begin{align*}
    P&:=\conv
    \bigl(φ\bigl[ψ[\poskegel{(\direktsumme{E}\ganzzahlen)}]\bigr]\bigr)\\
    P'&:=\linhuelle^+\bigl(φ\bigl[ψ[\poskegel{(\direktsumme{E}\ganzzahlen)}]\bigr]\bigr).
  \end{align*}
  As $P'$contains any positive linear combinations of its elements it
  also contains the special case of convex combinations. This implies
  $P\subseteq P'$.
  Since the zero vector $\nullvektor\in\direktsumme{E}\ganzzahlen$
  lies in the positive cone and for each positive vector $\vektor
  x\in\poskegel{(\direktsumme E\ganzzahlen)}$ its multiple
  $n\vektor{x}\in\poskegel{(\direktsumme E\ganzzahlen)}$ is positive
  for every
  $n\in\natzahlen$, for each $α\in\poskegel{\reellezahlen}$ we also
  have $α\vektor{x}\in P$. This can be deduced from the equation
  \[
  \alpha\sum_{i=0}^n \beta_{i}\vektor
  x_{i}=\sum_{i=0}^n\alpha\beta_{i}\vektor x_{i}.
  \]
 Thus, for any two vectors $\vektor{x}, \vektor{y}\in P$ and any two
  positive numbers $α,β\in \poskegel{\reellezahlen}$ we
  get:\footnote{For $α=β=0$ we already proved $\nullvektor\in
    \ganzzahlen\subseteq P$.}
  \[
  \alpha\vektor x+\beta\vektor y =\frac
  \alpha{\alpha+\beta}\bigl((\alpha+\beta)\vektor x\bigr) +
  \left(1-\frac \alpha{\alpha+\beta}\right) (\alpha+\beta)\vektor y
  \in P,
  \]
  which proves $P'\subseteq P$ and together with the preceding results
  $P=P'$.
  Furthermore, the set $P$ is a subsemigroup of $\direktsumme E\reellezahlen$. As
  this vector space is commutative, $P$ is already an invariant semigroup.

The set $\mathfrak B:=\Menge{\vektor e_a}{a\in E}$ is a generating set
of the module  $\direktsumme
  E\ganzzahlen$. Thus, $\mathfrak B$ is a basis of
  $\direktsumme{E}\reellezahlen$. Since $P\subseteq \direktsumme{E}\reellezahlen$, the following
  equation holds:
\[\linhuelle{P}=\linhuelle\Bigl(\linhuelle^+\bigl(\poskegel{(\direktsumme{E}\ganzzahlen)}\bigr)\Bigr).\]

\iffalse
can be considered as a
generating set for $P$ in the following sense: each vector $\vektor
x\in P$ can be written as linear combination of vectors from
$\mathfrak B$ with non-negative coefficients.  As already shown the set
$\linhuelle^+(\mathfrak B)$ of linear combinations of elements from
$\mathfrak B$ with non-negative coefficients is a subset of the convex
hull 
$\conv\bigl(\poskegel{(\direktsumme E\ganzzahlen)}\bigr)$. On the
other hand it is convex with respect to the vector space $\direktsumme
E\reellezahlen$, since any convex combination uses positive factors.
So we have found a representation of the convex hull  of 
$\poskegel{(\direktsumme E\ganzzahlen)}$ in $\direktsumme
E\reellezahlen$.

The next step it to prove that the integer elements of $P$ are those
elements of $\poskegel{(\direktsumme E\ganzzahlen)}$, which fulfil the
condition $P\cap\ganzzahlen^M=\poskegel{(\direktsumme E\ganzzahlen)}$.
Since $\mathfrak B$ is a generating system of $\poskegel{(\direktsumme
  E\ganzzahlen)}$, the identities $\linhuelle {\mathfrak B} =
\linhuelle\bigl(\poskegel{(\direktsumme E\ganzzahlen)}\bigr) =
\linhuelle P$ and  $\linhuelle^+ {\mathfrak B} =
\linhuelle^+\bigl(\poskegel{(\direktsumme E\ganzzahlen)}\bigr) =
\linhuelle^+ P=P$ hold.
\fi

Thus, for any vector $\vektor x \in \direktsumme{E}\ganzzahlen\cap P$
there exists a linear combination of elements of a linear independent
set $B\subseteq \poskegel{(\direktsumme{E}\ganzzahlen)}$ and a positive integer $n\in \natzahlen$
such that $\vektor x\in\sum_{i=1}^n\alpha_i\vektor a_i$ where
$\alpha_i\in\poskegel\reellezahlen$ and $\vektor a_i\in B$. Since the range of
$\vektor x$ is a set of integers, the coefficients $\alpha_i$ cannot
be irrational, because they are solutions of the system of linear
equations \[\forall m\in E:\sum_{i=1}^n\alpha_i\vektor a_i(m)=\vektor
x(m).\] As $\vektor{a}_i(m)\in\ganzzahlen$ and
$\vektor{x}(m)\in\ganzzahlen$ are integers, the coefficients $α_i$ are
rational numbers. Let $k$ be the least common denominator of all those
numbers $α_i$. Then we can multiply the system with $k$ and get the
modified system of linear equations
 \[
 \forall m\in E:\sum_{i=1}^n(k\alpha_i)\vektor a_i(m)=\vektor x'(m),
\]
where $\vektor{x}'$ is defined as $\vektor{x}':=k\vektor{x}$. Thus, the
vector $\vektor x'$ is positive in $\direktsumme{E} \ganzzahlen$ and
so the vector $\vektor{x}$ is non-negative. Since the order on
$\direktsumme{E}\ganzzahlen$ is semiclosed, this also implies
$\vektor{x}\in\poskegel{(\direktsumme{E}\ganzzahlen)}$.  Consequently,
we showed that $\vektor x\in P\cap\direktsumme{E}\ganzzahlen$ holds iff $\vektor x\in\poskegel{(\direktsumme
  E\ganzzahlen)}$.

So we get $P\cap\direktsumme{E}\ganzzahlen=\poskegel{(\direktsumme
  E\ganzzahlen)}$ and a partial order ${\leq'}\subseteq{\direktsumme
E\reellezahlen\times\direktsumme E\reellezahlen}$ defined by
\[
\vektor x\leq'\vektor y,\text{ iff }\vektor x-\vektor y\in P.
\]
Finally, we mapped each order $\leq$ on $\direktsumme E\ganzzahlen$
onto an order $\leq'$ on $\direktsumme E\reellezahlen$ such that the
positive cone of $\leq'$ reproduces the positive cone of $\leq$ by the
intersection with $\direktsumme{E}\ganzzahlen$. From that follows that the mapping
from the set of partial orders on the group $\direktsumme
E\ganzzahlen$ into the set of vector space orders on $\direktsumme
E\reellezahlen$, is injective.

The proof of the other two combinations follows the same path.
\end{proof}

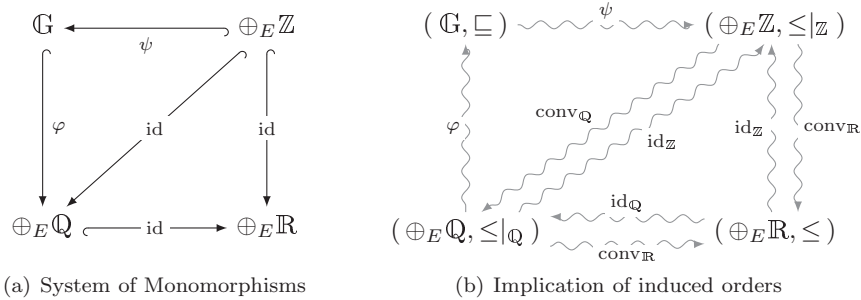
\begin{figure}[tbp]%
  \parbox{0.45\linewidth}{\subfigure[][System of Monomorphisms]{\label{subfig:gzqr-monomorphismen2}%
      \begin{tikzpicture}[description/.style={fill=white,
          inner sep=2pt},>=latex,x=0.9em,y=1em]%
        \useasboundingbox (-8.5,-5) rectangle (8.5,5);%
        \matrix (m) [matrix of math nodes, row sep = 6em, 
        column sep=5.5em,text height=1.5ex,text depth=0.25em,ampersand replacement=\&]%
        {\gruppe\&\direktsumme E\ganzzahlen\\
          \direktsumme E\rationalezahlen\& \direktsumme E\reellezahlen \\};%
        \path[right hook->,font=\scriptsize] (m-1-2) edge node[auto] {$\psi$}  (m-1-1);%
        \path[right hook->,font=\scriptsize] (m-1-1) edge node[auto] {$\varphi$} (m-2-1);%
        \path[left hook->,font=\scriptsize] (m-2-1) edge node[description] {$\id$} (m-2-2);%
        \path[right hook->,font=\scriptsize] (m-1-2) edge node[description] {$\id$} (m-2-2);%
        \path[right hook->,font=\scriptsize] (m-1-2) edge node[description] {$\id$} (m-2-1);%
      \end{tikzpicture}%
    }%
  }%
  \parbox{0.55\linewidth}{%
    \subfigure[][Implication of induced orders]{\label{subfig:gzqr-ordnungen-beide-richtungen}%
      \begin{tikzpicture}[x=1em,y=1em,description/.style={inner
          sep=2pt,black,auto,fill=white,swap},>=latex,decoration={snake,
          amplitude=0.15em,
          foot of=human,
          foot length=1.5em,
          stride length=4em,foot sep=-0.7em}]%
        \useasboundingbox (-9,-5) rectangle (9,5);%
        \matrix (m) [matrix of math nodes, row sep = 6.4em, 
        column sep=5.5em,text height=1.5ex,text depth=0.25ex, ampersand replacement=\&]%
        {  \paar\gruppe\sqsubseteq \&\paar{\direktsumme E\ganzzahlen}{{\leq}|_\ganzzahlen}\\
          \paar{\direktsumme E\rationalezahlen}{{\leq}|_\rationalezahlen} \&  \paar{\direktsumme E\reellezahlen}{\leq} \\};
        \path[<-,font=\scriptsize,gray]    (m-1-2) edge[decorate] node[description] {$\psi$}  (m-1-1);%
        \path[<-,font=\scriptsize,gray]   (m-1-1) edge[decorate] node[description] {$\varphi$} (m-2-1) ;%
        \path[->,font=\scriptsize,gray]    (m-2-1.-10) edge[decorate] node[description] {$\conv_\reellezahlen$} (m-2-2.-170);%
        \path[->,font=\scriptsize,gray]    (m-2-2.170) edge[decorate] node[description] {$\id_\rationalezahlen$} (m-2-1.10);%
        \path[->,font=\scriptsize,gray,swap=false]    (m-1-2.-40) edge[decorate] node[description] {$\conv_\reellezahlen$} (m-2-2.40);%
        \path[->,font=\scriptsize,gray,swap=false]    (m-2-2.100) edge[decorate] node[description] {$\id_\ganzzahlen$} (m-1-2.-100);%
        \path[->,font=\scriptsize,gray]   (m-1-2.-161) edge[decorate] node[description] {$\conv_\rationalezahlen$} (m-2-1.51);%
        \path[->,font=\scriptsize,gray]    (m-2-1.20) edge[decorate] node[description] {$\id_\ganzzahlen$} (m-1-2.-120);%
      \end{tikzpicture}%
    }%
  }%
  \caption{Transfer of orders on groups with respect to monomorphisms
    and convex hulls from an abelian Group \gruppe{} onto a real
    vector space and back}\label{fig:gzqr-ordnung-monomorphismus-beide-richtungen}
\end{figure}

\noindent As already stated, this theorem ensures that the paths in the lower right triangle
in Figure \ref{subfig:gzqr-ordnungen-beide-richtungen} are
equivalent. Due to the additional arrow from
$\paar{\direktsumme{E}\ganzzahlen}{{\leq}|_\ganzzahlen}$ to
$\paar{\direktsumme{E}\rationalezahlen}{{\leq}|_\rationalezahlen}$
this is not clear for the upper left triangle in this figure. The
following lemma closes this gap.

\begin{lemma}\label{lemma:ordungsuebertragung-G-Q}
  If the group \gruppe{} is semiclosed ordered, $E$ is a maximal
  independent set in $\gruppe$ and the monomorphisms $φ$ and $ψ$ are
  defined according to the Equations \eqref{eq:homomorphismus-G-Q} and
  \eqref{eq:z-g-mornomorphismus} the following equation holds:
  \begin{equation}\label{eq:conf_psi=phi}
    \conv_\rationalezahlen\bigl[\psi^{-1}[\poskegel\gruppe\cap\erzgruppe E]\bigr]\cap\varphi[G] =
    \varphi[\poskegel \gruppe]
 \end{equation}
\end{lemma}

\begin{proof}
  Let $\vektor{x}\in φ[\poskegel{\gruppe}]$. Then there exists an
  element $g\in \poskegel{\gruppe}$ such that $φ(g)=\vektor{x}$. Since
  $E$ is a maximal independent set there exists a natural number
  $n\in\natzahlen\setminus\menge{0}$ such that $ng \in\erzgruppe{E}$
  is an element of the subgroup which is generated by $E$. Thus $ng$
  is contained in the image of the group homomorphism $ψ$. Considering
  $\direktsumme{E}\rationalezahlen$ the vector $\frac{1}{n}ψ^{-1}(ng)$
  is an element of the convex hull
  $\conv_\rationalezahlen\bigl(\menge{\nullvektor
    ,ψ^{-1}(ng)}\bigr)$. Since the neutral element is always an
  element of the positive cone, we have shown one direction of
  \eqref{eq:conf_psi=phi}:
\[
    \conv_\rationalezahlen\bigl[\psi^{-1}[\poskegel\gruppe\cap\erzgruppe E]\bigr]\cap\varphi[G] \supseteq
    \varphi[\poskegel \gruppe].
\]

Let $\vektor{x}\in
\conv_\rationalezahlen\bigl[\psi^{-1}[\poskegel\gruppe\cap\erzgruppe
E]\bigr]\cap\varphi[G]$. Then there exists elements
$g\in G$ and $g_1,\ldots,g_n\in\poskegel{\gruppe}\cap\erzgruppe{E}$,
positive rational numbers $x_1,\ldots,x_n\in\poskegel{\rationalezahlen}$ and a
positive integer $n \in\natzahlen\setminus\menge{0}$ such that the
following condition holds:
\begin{equation*}
 x = φ(g) =\sum_{i=1}^n x_iψ^{-1}(g_i),\text{ where }
  \sum_{i=1}^n x_i=1.
\end{equation*}

Let $q\in\natzahlen$ be the smallest common denominator of
$\menge{x_1,\ldots,x_n}$. Then the products
$\menge{qx_1,\dotsc,qx_n}\subseteq\natzahlen$ are positive integers.
From Theorem \ref{lemma:Monomorphismus-Z-G} and Equation 
\eqref{eq:phi-und-psi-kommutieren-auf-Z} we infer that for each
$i\in\menge{1,\dotsc,n}$ we can use the identity
$ψ^{-1}(g_i)=φ(g_i)$. Using these two modifications we get the equation
\[
  qφ(g) =\sum_{i=1}^n qx_iφ(g_i).
\]
Since $φ$ is an injective group homomorphism, we can rewrite this
equation as 
\begin{equation*}
  φ(qg) = φ\left(\sum_{i=1}^n qx_ig_i\right)\text{ and get }
  qg = \sum_{i=1}^n qx_ig_i.
\end{equation*}
Thus, $qg\in\poskegel{\gruppe}$, which implies
$g\in\poskegel{\gruppe}$, since the order of $\gruppe$ is
semiclosed. This proves the other inclusion 
\[
    \conv_\rationalezahlen\bigl[\psi^{-1}[\poskegel\gruppe\cap\erzgruppe E]\bigr]\cap\varphi[G] \subseteq
    \varphi[\poskegel \gruppe]
\]
and, thus, the lemma is proved.
\end{proof}

\section{Linear orders in real vector spaces}\label{sec:linear-order-in-real-vector-space}

Here, we will discuss linear vector space orderings. As already
mentioned in the last section this touches the area of fundamental
theorems of the theory of $\ell$-groups with the additional insight
how different orders can be represented in the same vector space.

We will prove that the positive cone in the generated subspace of a
class of archimedian equivalent elements of a linearly ordered vector
space is defined by linear half spaces which are bounded by linearly
ordered hyperplanes. Using Zorn's lemma, it can be shown that this is
equivalent to the description by linearly ordered bases and an
extension to Teh's result \cite{Teh:1961}. In contrast to that work,
we use properties of normed vector spaces to prove this result for all
torsion free abelian groups.

%\setcapindent*{0pt}%
% \setcapwidth[l]{0.1\linewidth}%
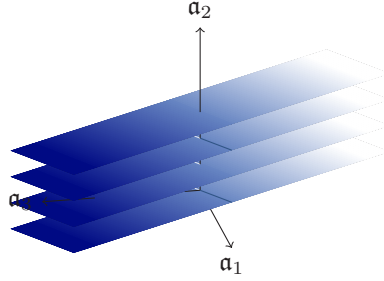
\begin{figure}\noindent\parbox[t]{0.4\linewidth}{\caption[Linear order in a real vector space]
{\linebreak[2]\raggedright Linear order in a real vector space}\label{fig:tikz-linordnung-im-vektorraum}}%
\parbox[t]{0.6\linewidth}{\centerline{\raisebox{-\height}{% Sketch output, version 0.3 (build 7d, Fri Nov 22 14:14:09 2013)
% Output language: PGF/TikZ,LaTeX
\begin{tikzpicture}[line join=round]
\tikzstyle{ann} = [fill=white,font=\footnotesize,inner sep=1pt]
\providecommand\vektor{}\tikzstyle{plane}=[shading=axis,top color=white,bottom
color=HKS41K80,shading angle=-50]\tikzstyle{plane2}=[shading=axis,top color=white,bottom
color=HKS41K80,shading angle=-30]\draw[arrows=<-,line width=.4pt](.42,-.767)--(0,0);
\filldraw[plane,draw=black,draw opacity=0.01,line width=0.1,fill opacity=1](-1.68,-.833)--(2.52,.513)--(1.68,.833)--(-2.52,-.513)--cycle;
\draw[draw opacity=0.8,draw=white](-.42,.16)--(0,0);
\draw[arrows=-,line width=.4pt](0,0)--(0,.347);
\draw[arrows=-,line width=.4pt](0,0)--(-1.4,-.102);
\draw[draw opacity=0.8,draw=HKS41K80](.42,.187)--(0,.347);
\filldraw[plane,draw=black,draw opacity=0.01,line width=0.1,fill opacity=1](-1.68,-.487)--(2.52,.86)--(1.68,1.18)--(-2.52,-.167)--cycle;
\draw[draw opacity=0.8,draw=white](-.42,.507)--(0,.347);
\draw[arrows=-,line width=.4pt](0,.347)--(0,.693);
\draw[arrows=->,line width=.4pt](-1.4,-.102)--(-2.1,-.153);
\filldraw[plane,draw=black,draw opacity=0.01,line width=0.1,fill opacity=1](-1.68,-.14)--(2.52,1.207)--(1.68,1.527)--(-2.52,.18)--cycle;
\draw[draw opacity=0.8,draw=white](-.42,.853)--(0,.693);
\draw[draw opacity=0.8,draw=white](-.42,1.2)--(0,1.04);
\draw[arrows=-,line width=.4pt](0,.693)--(0,1.04);
\draw[draw opacity=0.8,draw=HKS41K80](.42,.88)--(0,1.04);
\filldraw[plane,draw=black,draw opacity=0.01,line width=0.1,fill opacity=1](-1.68,.207)--(2.52,1.553)--(1.68,1.873)--(-2.52,.527)--cycle;
\draw[draw opacity=0.8,draw=HKS41K80](.42,-.16)--(0,0);
\draw[draw opacity=0.8,draw=HKS41K80](.42,.533)--(0,.693);
\draw[arrows=->,line width=.4pt](0,1.04)--(0,2.167);
\path (.42,-.767) node[below] {$\vektor a_1$}
     (0,2.167) node[above] {$\vektor a_2$}
     (-2.1,-.153) node[left] {$\vektor a_3$};\end{tikzpicture}% End sketch output
}}}%
\end{figure}%

At first we discuss possible linear orders on the vector space
$\direktsumme{E}\reellezahlen$ for a given set $E$. In each of these orders a hyperplane will
contain all infinitesimal elements, and the space will be filled with
copies of this hyperplane (c.\,f., Figure \ref{fig:tikz-linordnung-im-vektorraum}).

\begin{lemma}\label{lemma:hyperebene-definiert-ordnung}%
  Let \hyperebene\ be a linearly ordered hyperplane in the vector
  space $\direktsumme E\reellezahlen$ and $\vektor a\in\direktsumme
  E\reellezahlen\setminus\hyperebene$ a vector which is linearly
  independent from \hyperebene. Then the set
  \begin{equation}
    \poskegel\hyperebene\cup \bigl(\hyperebene + (\poskegel\reellezahlen\setminus\menge 0)\vektor a\bigr)
  \end{equation}
  defines a linear order on the vector space $\Gruppe{\direktsumme{E}\reellezahlen}+-\nullvektor$.
\end{lemma}
\begin{proof}
Let $P:=\poskegel\hyperebene\cup \bigl(\hyperebene +
(\poskegel\reellezahlen\setminus\menge 0)\vektor a\bigr)$ be the half
space, which is defined by the positive cone of the hyperplane $\poskegel{\hyperebene}$ and the vector
$\vektor{a}$. Since
$\direktsumme{E}\reellezahlen$ is a group considering the addition, We
have to prove the conditions
\eqref{eq:poskegelbedingung_halbgruppe} to
\eqref{eq:poskegelbedingung_invariant} considering $P$ as positive
cone of $\direktsumme{E}\reellezahlen$. This will prove that $P$ is an
invariant semigroup in $\direktsumme{E}\reellezahlen$ and for each
vector $\vektor x\in\direktsumme{E}\reellezahlen$ either  $\vektor x\in
P$ or $-\vektor x\in P$ holds.

Let $\vektor x\in P$ be a vector in $P$. Then either the condition
$\vektor x\in\poskegel\hyperebene$ or $\vektor x\in\hyperebene +
(\poskegel\reellezahlen\setminus\menge 0)\vektor a$ holds.
\begin{enumerate}
\item The condition $\vektor x\in\poskegel\hyperebene$ is
  true iff  $-\vektor x\in\negkegel\hyperebene$, since
  \hyperebene\ is a vector space itself and linearly ordered
  by precondition.
\item If $\vektor
  x\in\hyperebene+(\poskegel\reellezahlen\setminus\menge
  0)\vektor a$, then a unique vector $\vektor
  x'\in\hyperebene$ and a unique positive real number 
  $x\in\poskegel\reellezahlen\setminus\menge 0$ exist such that 
  $\vektor x = \vektor x' + x\vektor a$. This is equivalent to
  \begin{equation*}
    -\vektor x = -\vektor x' + (-x)\vektor a\in \hyperebene+(\negkegel\reellezahlen\setminus\menge 0)\vektor a
  \end{equation*}
\end{enumerate}
Furthermore, since
\hyperebene\ is as hyperplane in $\direktsumme{E}\reellezahlen$ and
$\vektor a\not\in\hyperebene$, the following holds:
\begin{align*}
	\direktsumme{E}\reellezahlen&=\hyperebene\vektorraumplus\reellezahlen\vektor a%\\&
		=\hyperebene+\reellezahlen\vektor a\\
		&=\hyperebene+(\poskegel\reellezahlen\setminus\menge 0)\vektor a
		\disjverein (\poskegel\hyperebene\setminus\menge \nullvektor)
		 \disjverein\menge{\nullvektor}\disjverein(\negkegel\hyperebene\setminus\menge \nullvektor)\disjverein 
			\hyperebene+(\negkegel\reellezahlen\setminus\menge 0)\vektor a\\
			&=P\setminus\menge \nullvektor\disjverein \menge \nullvektor \disjverein -P\setminus\menge \nullvektor
\end{align*}
In short: $\direktsumme{E}\reellezahlen = P\cup -P$, so the Condition
\eqref{eq:poskegelbedingung_linear} is fulfilled. From the two points
above follows $P\cap -P = \menge \nullvektor$, which is the Condition \eqref{eq:poskegelbedingung_echte_halbgrupe}. 

Let's consider the remaining Conditions \eqref{eq:poskegelbedingung_halbgruppe} and
\eqref{eq:poskegelbedingung_invariant}. Let $\vektor x,\vektor y\in P$
be two vectors in $P$. Then there exist unique decompositions
\begin{align*}
	\vektor x &= \vektor x' + x\vektor a\text{ and }\\
	\vektor y &= \vektor y' + y\vektor a
\end{align*}
where $\vektor x',\vektor y'\in\hyperebene$ and
$x,y\in\poskegel\reellezahlen$. For two non-negative real numbers $α,
β\in\poskegel{\reellezahlen}$ this leads to the equation
\begin{align*}
α\vektor x + β\vektor y = (α\vektor x' + β\vektor y') +(αx+βy)\vektor a.
\end{align*}
Since the values $α,β,x,y\geq 0$ are non-negative, the inequality
$αx+βy\geq 0$ is true. Equally, the sum 
$α\vektor x'+β\vektor y'\in\hyperebene$ belongs to the
hyperplane. Consequently, for a non-negative sum $αx+βy\neq 0$ also the
sum of the vectors
$α\vektor x + β\vektor y\in P$ is a member of the set $P$. On the other hand,
 the equation $αx+βy=0$ induces
$αx=βy=0$ as both values are non-negative. This means that $α\vektor
x\in\poskegel\hyperebene$, as well as $β\vektor
y\in\poskegel\hyperebene$. From the semigroup properties of the
positive cone of the hyperplane
$\hyperebene$ follows that also the linear combination
$α\vektor x + β\vektor y \in\poskegel\hyperebene$ is an element of
this subspace.
So we have shown for all vectors $\vektor x,\vektor y \in P$ from the
half space $P$ and all non-negative real numbers $α,β\in\reellezahlen$ that
the linear combination $α\vektor x+β\vektor y
\in P$ belongs to $P$, too. Thus, $P$ is a subsemigroup of
$\direktsumme{E}\reellezahlen$, so the Condition
\eqref{eq:poskegelbedingung_halbgruppe} holds, and --
since $\direktsumme{E}\reellezahlen$ is commutative -- $P$ is also an
invariant subsemigroup of the vector space, which proves condition
\eqref{eq:poskegelbedingung_invariant}. 
From Theorem \ref{satz:poskegelbedingungen} follows that $P$ defines
a linear order on the vector
space $\direktsumme{E}\reellezahlen$ considered as a group. 
Using the unique decomposition $\vektor{x} = \vektor{x}' +x\vektor{a}$
for any vector $\vektor{x}\in P$ with $\vektor{x}'\in\hyperebene$ we
infer that $α\vektor{x}'\in\hyperebene$ and $α\vektor{x}' +
αx\vektor{a}\in P$.
Thus the vector space $\direktsumme{E}\reellezahlen$ is a linearly
ordered vector space.
\end{proof}%

\noindent In the opposite direction we can only limit the codimention of the set
of infinitesimal elements.

\begin{lemma}\label{lemma:rand-ist-lin-unterraum}
  Let $\direktsumme{E}\reellezahlen = \OrdGruppe{\direktsumme
    M\reellezahlen}+-\nullvektor\leq$ be an ordered vector space. Then
  the boundary
  $\rand{\bigl(\poskegel{(\direktsumme{E}\reellezahlen)}\bigr)}$ of
  the positive cone is a linear subspace of
  $\direktsumme{E}\reellezahlen$.
\end{lemma}

\begin{proof}
  Let $\vektor x\in \rand{\bigl(\poskegel{(\direktsumme{E}\reellezahlen)}\bigr)}$ a boundary
  vector of the positive cone. Then for all
  $r\in\reellezahlen\setminus\menge\nullvektor$ the intersections
  $B_{r}(\vektor x)\cap
  \poskegel{(\direktsumme{E}\reellezahlen)}$ and $B_{r}(\vektor x)\cap
  \negkegel{(\direktsumme{E}\reellezahlen)}$ are non-empty. Furthermore, for $\vektor x
  \neq \nullvektor$  either
  $\vektor x\in \poskegel{(\direktsumme{E}\reellezahlen)}$ or $-\vektor x\in
  \poskegel{(\direktsumme{E}\reellezahlen)}$ holds.

 In fact, $\vektor x\in\rand{\bigl(\poskegel{(\direktsumme{E}\reellezahlen)}\bigr)}$ is satisfied iff
 $-\vektor x\in\rand{\bigl(\poskegel{(\direktsumme{E}\reellezahlen)}\bigr)}$.
 This follows as $\vektor x'\in B_{r}(\vektor x)\cap
  \poskegel{(\direktsumme{E}\reellezahlen)}$ is a member of the subset of the
  positive elements of a neighbourhood of the vector \vektor x iff $-\vektor x'\in
  B_{r}(-\vektor x) \setminus\poskegel{(\direktsumme{E}\reellezahlen)}$ is located in
  a corresponding subset of negative elements. An analogue result is
  true for $\negkegel{(\direktsumme{E}\reellezahlen)}$. Thus, the neighbourhood
  $B_r(-\vektor x)$ contains as well positive as well as negative elements.

 In a next step we show that $\rand{\bigl(\poskegel{(\direktsumme{E}\reellezahlen)}\bigr)}$ is a
  linear subspace in $\direktsumme{E}\reellezahlen$. To do that, we
  choose two vectors $\vektor x,\vektor
  y\in
  \rand{\bigl(\poskegel{(\direktsumme{E}\reellezahlen)}\bigr)}$. Then
  for arbitrary positive real numbers $r,\alpha,\beta \in \poskegel\reellezahlen\setminus \menge{0}$
  four vectors $\vektor x'$, $\vektor x''$, $\vektor y'$ and
  $\vektor y''$ exist in the following way:
  \begin{align*}
    \vektor x' &\in B_{\frac r{2\alpha}}(\vektor x)\cap \poskegel{(\direktsumme{E}\reellezahlen)}&
    \vektor x'' &\in B_{\frac r{2\alpha}}(\vektor x)\setminus \poskegel{(\direktsumme{E}\reellezahlen)}\\
    \vektor y' &\in B_{\frac r{2β}}(\vektor y)\cap \poskegel{(\direktsumme{E}\reellezahlen)}&
    \vektor y'' &\in B_{\frac r{2β}}(\vektor y)\setminus \poskegel{(\direktsumme{E}\reellezahlen)}
  \end{align*}

 Depending on the values of  $\alpha$ and $\beta$ this leads to the
 following six cases:
  \begin{description}
  \item[$\alpha>0$ and $\beta>0$] Here, we get the equations
    \begin{align*}
      \alpha\vektor x'+\beta\vektor y'&=
      \underbrace{\alpha\vektor x +
        \underbrace{\alpha\underbrace{(\vektor x'-\vektor x)}_{{}\in B_{\frac
            r{2\alpha}}(\nullvektor)\cap
          \poskegel{(\direktsumme{E}\reellezahlen)}}}_{\in
        B_{\frac r2}(\nullvektor)\cap
        \poskegel{(\direktsumme{E}\reellezahlen)}}
      + \beta\vektor y +\underbrace{\beta\underbrace{(\vektor y'-\vektor y)}_{{}\in B_{\frac r{2\beta}}(\nullvektor)\cap \poskegel{(\direktsumme{E}\reellezahlen)}}}_{\in B_{\frac r2}(\nullvektor)\cap \poskegel{(\direktsumme{E}\reellezahlen)}}}_
      {\in B_{r}(\alpha\vektor x + \beta\vektor y)\cap
        \poskegel{(\direktsumme{E}\reellezahlen)}}\text{ and}\\
      \alpha\vektor x''+\beta\vektor y''&=
      \underbrace{\alpha\vektor x
        + \underbrace{\alpha\underbrace{(\vektor x''-\vektor x)}_{{}\in B_{\frac
            r{2\alpha}}(\nullvektor)\setminus
          \poskegel{(\direktsumme{E}\reellezahlen)}}}_{\in B_{\frac r2}(\nullvektor)\setminus\poskegel{(\direktsumme{E}\reellezahlen)}}
      +\beta\vektor y + \underbrace{\beta\underbrace{(\vektor y''-\vektor y)}_{{}\in B_{\frac r{2\beta}}(\nullvektor)\setminus\poskegel{(\direktsumme{E}\reellezahlen)}}}_{\in B_{\frac r2}(\nullvektor)\setminus\poskegel{(\direktsumme{E}\reellezahlen)}}}_
      {\in B_{r}(\alpha\vektor x + \beta\vektor y)\setminus\poskegel{(\direktsumme{E}\reellezahlen)}}.
\end{align*}
Thus, in $B_{r}(\alpha\vektor x + \beta\vektor y)$ are positive as
well as negative elements of the vector space
$\direktsumme{E}\reellezahlen$, so $\alpha\vektor x+\beta\vektor y$
is also an element of $\rand{\bigl(\poskegel{(\direktsumme{E}\reellezahlen)}\bigr)}$.
\item[$\alpha>0$ and $\beta<0$] This case an be deduced to the last
  case using the representation
\[\alpha\vektor x + \beta \vektor y = \alpha\vektor x +
(-\beta)(-\vektor y),\] since $\vektor y$ is an element of the
boundary iff $-\vektor y$ is within the boundary, too.
\item[$\alpha<0$ and $\beta>0$] As the addition is commutative this is identical to the predeceasing
  case.
\item[$\alpha<0$ and $\beta<0$] With the same idea as in the
  predeceasing cases, the result follows from the first statement.
\item[$\alpha\neq 0$ and $\beta = 0$] Here, we have $\beta\vektor y =
  \nullvektor$. So analogue to the first case follows:
\begin{align*}
	\alpha\vektor x'+\beta\vektor y'&=
        \alpha\vektor x + \alpha\underbrace{(\vektor x'-\vektor x)}_{{}\in B_{\frac r{2\alpha}}(\nullvektor)\cap \poskegel{(\direktsumme{E}\reellezahlen)}}
        \in B_{r}(\alpha\vektor x)\cap \poskegel{(\direktsumme{E}\reellezahlen)}\\
        \alpha\vektor x''+\beta\vektor y''&=
        \alpha\vektor x + \alpha\underbrace{(\vektor x''-\vektor x)}_{{}\in B_{\frac r{2\alpha}}(\nullvektor)\setminus \poskegel{(\direktsumme{E}\reellezahlen)}}
        \in B_{r}(\alpha\vektor x)\setminus\poskegel{(\direktsumme{E}\reellezahlen)}\\
\end{align*}
\item[$\alpha = 0$] For $\beta \neq 0$ this case is analogue to
  the predeceasing case. Otherwise we know $\alpha\vektor x+\beta\vektor y = \nullvektor \in \rand{\bigl(\poskegel{(\direktsumme{E}\reellezahlen)}\bigr)}$.
\end{description}
So far, we have shown that the relevant operations of the vector
space, addition of vectors, and multiplication with scalars do not
exit
$\rand{\bigl(\poskegel{(\direktsumme{E}\reellezahlen)}\bigr)}$. Thus,
this set is a linear subspace of $\direktsumme{E}\reellezahlen$.  
\end{proof}

\noindent The next step is to determine the dimension of the boundary

\begin{lemma}\label{lemma:rand-ist-lin-unterraum2}
  Let $\direktsumme{E}\reellezahlen = \OrdGruppe{\direktsumme
    M\reellezahlen}+-\nullvektor\leq$ be an ordered vector space. Then
  for the codimension of the boundary
  $\rand{\bigl(\poskegel{(\direktsumme{E}\reellezahlen)}\bigr)}$ the following equation holds:
  \begin{equation}
    \codim \rand{\bigl(\poskegel{(\direktsumme{E}\reellezahlen)}\bigr)} 
    =
    \begin{cases}
      0,&\text{ iff }{\inneres{\bigl(\poskegel{(\direktsumme{E}\reellezahlen)}\bigr)}} = ∅\\
      1,&\text{ else}
    \end{cases}
  \end{equation}
\end{lemma}

\begin{proof}
  Let us assume that the interior of
  $\poskegel{(\direktsumme{E}\reellezahlen)}$ is non-empty. By
  Lemma~\ref{lemma:rand-ist-lin-unterraum} the boundary
  $\rand{\bigl(\poskegel{(\direktsumme{E}\reellezahlen)}\bigr)}$ is a
  linear subspace of $\direktsumme{E}\reellezahlen$.  To be a
  hyperplane, this set has to be a subspace with codimension 1. To
  show this, we assume that we have at least codimension 2 and deduce
  a contradiction from that assumption. At the beginning we have to
  find a vector, which is not in
  $\rand{\bigl(\poskegel{(\direktsumme{E}\reellezahlen)}\bigr)}$. Let
  $\vektor x, \vektor y\in\poskegel{(\direktsumme{E}\reellezahlen)}$
  be two arbitrary positive elements of
  $\direktsumme{E}\reellezahlen$. Then either $\vektor x\leq \vektor
  y$ or $\vektor y\leq \vektor x$ follows from the linear order. Let
  w.l.o.g.\ $\vektor x\leq\vektor y$. This leads to the
  inequality $\nullvektor\leq\vektor y-\vektor x$ and we get for
  $0\leq \lambda\leq 1$:
  \begin{align*}
    \lambda\vektor x+(1-\lambda)\vektor y&=\lambda \vektor x +
    (1-\lambda)(\vektor y-\vektor x)+(1-\lambda)\vektor x\\
    &= \vektor x + (1-\lambda)(\vektor y-\vektor x)\geq \vektor x,
  \end{align*}
  since $0\leq 1-\lambda$ holds, and the positive cone
  $\poskegel{(\direktsumme{E}\reellezahlen)}$ is convex in the sense of
  ordered groups. Consequently,
  $\poskegel{(\direktsumme{E}\reellezahlen)}$ is convex in the sense of
  vector space and with the positive cone, the negative cone
  $\negkegel{(\direktsumme{E}\reellezahlen)}=-\poskegel{(\direktsumme{E}\reellezahlen)}$
  has this property, too. So we have the identity  $\linhuelle(\poskegel{(\direktsumme{E}\reellezahlen)}) =
  \linhuelle(\negkegel{(\direktsumme{E}\reellezahlen)})=\direktsumme{E}\reellezahlen$,
  since 
  $\direktsumme{E}\reellezahlen =
  \negkegel{(\direktsumme{E}\reellezahlen)}\cup\poskegel{(\direktsumme{E}\reellezahlen)}$.
  As for each base vector $\vektor e_a$ with $a\in E$ either $\vektor
  e_a$ or the vector $-\vektor{e_a}$ is positive, and we can replace in
  the basis
  each negative basis vector with its inverse we can assume
  w.l.o.g. that 
  $\Menge{\vektor e_a}{a\in
    E}\subseteq\poskegel{(\direktsumme{E}\reellezahlen)}$.
  This is also true for its convex hull, which leads to $\conv \Menge{\vektor
    e_a}{a\in E}\cap {\inneres{\bigl(\poskegel{(\direktsumme{E}\reellezahlen)}\bigr)}}\neq \emptyset$.

  Let $\vektor a\in {\conv \Menge{\vektor
      e_a}{a\in E}}\cap {\inneres{\bigl(\poskegel{(\direktsumme{E}\reellezahlen)}\bigr)}}$,
  and let us assume that the codimension of $\rand{\bigl(\poskegel{(\direktsumme{E}\reellezahlen)}\bigr)}$
  is greater than 1. Then besides $\vektor a$ another
  vector $\vektor x\in\poskegel{(\direktsumme{E}\reellezahlen)}$ exists which
  is linearly independent from $\menge {\vektor a}\cup {\rand{\bigl(\poskegel{(\direktsumme{E}\reellezahlen)}\bigr)}}$.
  Consequently, there exists an isomorphism $ξ$, which maps the
  linear closure of $\vektor a$ and $\vektor x$ onto $\reellezahlen^2$ such
  that $ξ(\vektor a)=\paar 10$ and $ξ(\vektor
  x)=\paar 01$ are true. Obviously, the natural order on
  $\reellezahlen^2$ is compatible to the order, which is transferred
  from $\direktsumme{E}\reellezahlen$ to $\reellezahlen^2$ by
  $ξ$ . So
  we can extend the order on  $\reellezahlen^2$ linearly in such a way,
  that $ξ$ is an o-isomorphism.

  Using this order extension, we know that  $\paar {-\frac 12}{-1}<\paar 00<\paar 1{\frac
    12}$ holds. Now, we will define a mapping  
  \[f:\reellezahlen\to\menge{-1,0,1}:\lambda\mapsto\begin{cases}1&\paar[\big] 1{\frac12}-\lambda\paar[\big]{\frac 32}{\frac 32} > \paar 00\\
    -1&\paar[\big] 1{\frac12}-\lambda\paar[\big]{\frac 32}{\frac 32} < \paar 00\\
    0&\text{else}\end{cases}.\]

  Obviously, $f(0) = 1$ and $f(1) = -1$. Thus, $f$ is discontinuous
  according the natural topology on $\reellezahlen$. Let $\lambda_{0}$
  be a place of discontinuity of $f$. Since $\reellezahlen^2$ is linearly
  ordered and 
  \[
  \paar 00 \not\in \paar[\left]1%\right[]
  {\frac12}-\reellezahlen\paar[\auto]{\frac 32}{\frac 32},
  \]
  one of the identities
  $f(\lambda_{0})=1$ or $f(\lambda_{0})=-1$ is met. Thus, for each
  $\varepsilon\in\reellezahlen$ there exists a real number $\lambda_{1}\in\reellezahlen$ such that
  $|λ_0-λ_1|<ε$ and $f(\lambda_{1}) = -f(\lambda_{0})$. This implies that
  in $B_{\varepsilon}\Bigl(\varphi^{-1}\bigl(\paar
  1{\frac12}-\lambda_{0}\paar{\frac 32}{\frac 32}\bigr)\Bigl)$ 
  a positive vector $\vektor
  y'\in\poskegel{(\direktsumme{E}\reellezahlen)}$ and a 
  negative vector $\vektor
  y''\not\in\poskegel{(\direktsumme{E}\reellezahlen)}$ exist. So
  $\varphi^{-1}\bigl(\paar 
  1{\frac12}-\lambda_{0}\paar{\frac 32}{\frac
    32}\bigr)\in\rand{\bigl(\poskegel{(\direktsumme{E}\reellezahlen)}\bigr)}\setminus\menge
  \nullvektor$, this is a contradiction to the assumption that this set
  has no non-zero representatives in the linear closure of $\vektor x$ and
  $\vektor a$. Thus, the boundary
  $\rand{\bigl(\poskegel{(\direktsumme{E}\reellezahlen)}\bigr)}$ is a hyperplane in
  $\direktsumme{E}\reellezahlen$ with
  $\rand{\bigl(\poskegel{(\direktsumme{E}\reellezahlen)}\bigr)}\vektorraumplus\reellezahlen\vektor a=\direktsumme{E}\reellezahlen$.

  As already shown ${\conv \Menge{\vektor e_a}{a\in E}}$ can be
  considered to be a subset of the positive cone. It also generates
  the vector space, so the set ${\conv \Menge{\vektor e_a}{a\in E}}\cap
  {\inneres{\bigl(\poskegel{(\direktsumme{E}\reellezahlen)}\bigr)}}$
  is empty iff
  $\inneres{\bigl(\poskegel{(\direktsumme{E}\reellezahlen)}\bigr)}$ is
  empty. In that case the codimension is zero.
\end{proof}

\noindent Consequently the boundary of the positive cone in any ordered vector
space is either a hyperplane or it is the space itself. With this
knowledge we can define some additional notations for a vector $\vektor a$ of an ordered vector space:
\begin{align}
  \vektorraum_{\vektor a} &:= \linhuelle
  \Menge{\vektor x∈\direktsumme{E}\reellezahlen}{\nullvektor≤\vektor x≤|\vektor a|}\text{ and}\label{eq:Va}\\
  \hyperebene{}_{\vektor a} &:= \linhuelle\Menge{\vektor
    x∈\direktsumme{E}\reellezahlen}{\nullvektor≤\vektor x\infinitesimal<|\vektor a|}\label{eq:Ha}
\end{align}

In general the subspace $\vektorraum_{\vektor a}$ exists only for such
vectors $\vektor a$ that have assigned a unique absolute element $|\vektor
a|≥\nullvektor$. For that it is sufficient that there exists a lattice ordered
subspace containing $\vektor a$. If the vector space is linearly
ordered $\hyperebene_{\vektor a}$ is a hyperplane and
Equation~(\ref{eq:Ha}) can be extended to 
\[
  \hyperebene{}_{\vektor a} = \rand{\bigl(\poskegel{(\vektorraum_{\vektor a})}\bigr)}
\]

After we have shown that all
linearly ordered hyperplanes define linear orders, we can prove that
all linear orders on the vector space $\direktsumme{E}\reellezahlen$
have such a representation.
 
\begin{theorem}\label{satz:poskegelrand=Hyperebene}
  Let $\direktsumme{E}\reellezahlen = \OrdGruppe{\direktsumme
    E\reellezahlen}+-\nullvektor\leq$ be an ordered vector space. The
  partial order $\leq$ is linear iff one of the following conditions is true:
  \begin{enumerate}
  \item \label{item:1}The positive cone 
    $\poskegel{(\direktsumme{E}\reellezahlen)}$ has the form
    \begin{equation}
      \poskegel{(\direktsumme{E}\reellezahlen)}=\poskegel\hyperebene \disjverein\bigl(\hyperebene +(\poskegel\reellezahlen \setminus\menge 0)\vektor a\bigr),
    \end{equation}
    where $\hyperebene$ is a linearly ordered hyperplane in
    $\direktsumme{E}\reellezahlen$ and $\vektor a\in\direktsumme{E}\reellezahlen\setminus\hyperebene$
    is a positive vector from the interior of the positive cone.
  \item\label{item:2} The vector space $\direktsumme{E}\reellezahlen$ is the union
    of an infinite chain of proper subspaces of the form given in \ref{item:1}.
  \item\label{item:3} $E=∅$.
  \end{enumerate}
\end{theorem}

\begin{proof}
  As \ref{item:3} is trivial we have to prove that for every non-trivial vector space either \ref{item:1}
  or \ref{item:2} is true. By
  Lemma~\ref{lemma:hyperebene-definiert-ordnung} the vector space is
  linearly ordered if the condition \ref{item:1} holds. By
  Lemma~\ref{lemma:rand-ist-lin-unterraum} and
  Lemma~\ref{lemma:rand-ist-lin-unterraum2} the remaining part consists of two questions:
  \begin{enumerate}
  \item Is it sufficient that a vector space is the union of an
    infinite chain of linearly ordered subspaces in order to get a linear
    order?
  \item Can every linear order with empty interior of the positive
    cone be described by such a chain?
  \end{enumerate}

  In order to answer the first question, suppose there are two
  incomparible elements $\vektor a,\vektor
  b∈\direktsumme{E}\reellezahlen$. Then there exist linearly ordered
  subspaces $\vektorraum[U]_a$ and $\vektorraum[U]_b$ in the chain of
  subspaces that generates the order such that $\vektor
  a∈\vektorraum[U]_a$ and $\vektor b∈\vektorraum[U]_b$. As both vector
  spaces belong to the same chain either
  $\vektorraum[U]_a≤\vektorraum[U]_b$ or
  $\vektorraum[U]_b≤\vektorraum[U]_a$ holds. As $\vektor a$ and
  $\vektor b$ are incomparible and the vector spaces are linearly
  orderded we get the additional conditons $\vektor
  a\not∈\vektorraum[U]_b$ and $\vektor b\not∈\vektorraum[U]_a$. Thus
  both vector spaces do not belong to the same chain in the subspace
  lattice, which is a contradiction to the assumption. Consequently,
  we can answer the first question with “Yes”.

  To answer the other question let $\vektor a,\vektor
  b∈\direktsumme{E}\reellezahlen$ be two vectors with $\vektor
  a<\vektor b$. Then either $\vektor a∈\hyperebene_b$ or
  $\hyperebene_a\vektorraumplus\reellezahlen\vektor
  a=\hyperebene_b\vektorraumplus\reellezahlen\vektor
  a=\hyperebene_b\vektorraumplus\reellezahlen\vektor b$ holds. Thus, the set 
  \[
  \mathfrak W_{\vektor b}:=\Menge{\hyperebene_{\vektor
      a}\vektorraumplus \reellezahlen\vektor a}{\vektor a≤\vektor b}
  \]
  is a chain in the subspace lattice of $\direktsumme{E}\reellezahlen$.
  Thus, each order ideal in the ordered set $\paar{\mathfrak W}{≤}$ with
  \[\mathfrak W = ⋃_{\vektor b∈\direktsumme{E}\reellezahlen}\mathfrak W_{\vektor b}\]
  is a chain. With a similar argumentation we can
  prove that each order filter in this ordered set is a chain, too.
  Thus, $\paar{\mathfrak W}{≤}$ is a chain itself. As it describes the
  linear order as discussed in the previous question, also this
  question can be answered with “Yes”.
\end{proof}

\begin{corollary}
  Considering the hyperplane $\hyperebene$ from
  Lemma~\ref{lemma:hyperebene-definiert-ordnung},
  Condition~\ref{item:1} and one of the linear orders, which arise
  from this hyperplane, in Theorem \ref{satz:poskegelrand=Hyperebene}
  the following equation holds:
\begin{equation}
\poskegel\hyperebene = \poskegel{(\direktsumme{E}\reellezahlen)}\cap\rand{\bigl(\poskegel{(\direktsumme{E}\reellezahlen)}\bigr)}.
\end{equation}
\end{corollary}
\begin{proof}
  The positive cone of the hyperplane is uniquely defined.
\end{proof}

\noindent The previous lemma shows that linear vector space orders on
$\direktsumme{E}\reellezahlen$ can be characterised by hyperplanes in
sufficiently chosen subspaces. Thus, each linear order on a given
abelian group can be represented by one or an infite chain of linearly
ordered hyperplanes in this space. Obviously, for each linearly
ordered hyperplane $\hyperebene$ two linear orders on
$\direktsumme{E}\reellezahlen$ exist: A first one, where the open
half-space marked by a given vector
$\vektor{a}\in\direktsumme{E}\reellezahlen\setminus\hyperebene$ is
positive and a second one where it is negative. Consequently, a given
basis of the vector space can be a subset of the positive cone only in
one of these two cases.

As the lattice of linear subspaces is inductively ordered and since all
subspaces of a given subspace are subspaces of hyperplanes in the
bigger subspace, we can apply Zorn's lemma and receive the result
that each linear order on a vector space \vektorraum{} is uniquely
defined by a given linearly ordered base $\Halbord{\mathfrak B}{\leq}$
according to the following construction: 

Let $\vektor{a}\in\mathfrak B$ be a base vector. Furthermore, let
$\vektorraum_{\vektor{a}}=
\linhuelle{\Menge{\vektor{a}'}{\vektor{a}'\leq\vektor{a}}}$ and
$\hyperebene_{\vektor{a}}:=
\linhuelle{\Menge{\vektor{a}'}{\vektor{a}'<\vektor{a}}}$ a linearly ordered
hyperplane. Then the subspace $\vektorraum_\vektor{a}$ can be
linearly ordered by Lemma \ref{lemma:hyperebene-definiert-ordnung}
such that $\vektor a>\nullvektor$ is positive and
$\vektorraum_{\vektor{a}}$ is itself ordered as described above for
the vectors $\vektor a''<\vektor{a}$ in $\mathfrak B$.

\iffalse
$\hyperebene_{\vektor{a}}$ be the linearly ordered hyperplane
$\hyperebene_{\vektor{a}}:=
\linhuelle{\Menge{\vektor{a}'}{\vektor{a}'<\vektor{a}}}$ which is
ordered by a relation $\leq$ according to $\vektor{x}\leq\vektor{y}$,
iff either $x\in\hyperebene_{\vektor{a''}}$ for some
$\vektor{a''}\in\mathfrak B$ and
$\vektor{y}\not\in\hyperebene_{\vektor{a''}}$ or
$\menge{\vektor{x},\vektor{y}}\subseteq\vektorraum_{\vektor{a''}}\setminus\hyperebene_{\vektor{a''}}$
for some $\vektor{a}''$ and with respect to the basis $\mathfrak B$
the scalar coefficient of $\vektor x$, which belongs to $\vektor{a}''$
is less or equal to the coefficient of $\vektor{y}$ in the basis
representation of these vectors. Then the subspace
$\vektorraum_\vektor{a}$ can be linearly ordered by Lemma
\ref{lemma:hyperebene-definiert-ordnung} such that $\vektor
a>\nullvektor$ is positive and $\vektorraum_{\vektor{a}}$ is itself
ordered as described above for the vectors $\vektor a''<\vektor{a}$ in
$\mathfrak B$.
\fi

Moreover, the vector spaces
$\vektorraum_{\vektor{a}}$ form a chain in the lattice of subspaces of
$\direktsumme{E}\reellezahlen$ that fulfils the following condition:
\[
\direktsumme{E}\reellezahlen=\bigcup_{\vektor{a}\in\mathfrak
  B}\vektorraum_{\vektor{a}}.
\]
Thus, $\direktsumme{E}\reellezahlen$ can be linearly ordered, since
the order on the hyperplane above could be chosen
arbitrarily. Obviously, for a given linear order the basis is not
unique as the subspaces and their order is not influenced by the
multiplication of basis vectors with scalars.

We will discuss this construction not in further details,\footnote{The
complete discussion will be made available via \cite{dd/MATH-AL-02-2011}} as some
aspects of this topic have 
been used already in different flavours in other literature e.g., the works of The~\cite{Teh:1961} and
Zajceva \cite{M.I.Zajceva:1953}.

\section{Linear orders on groups}

The characterisations of linear orders on groups in The's~\cite{Teh:1961} and
Zajceva's~\cite{M.I.Zajceva:1953} articles both have drawbacks.  The first description
needs different vector spaces to describe all linear orders on a given
group, while the inductive nature of the latter restricts it to mainly
countable dimensional vector spaces. Nevertheless, several authors have
discussed the existence of linear orders on infinitely generated
groups by considering finitely generated subgroups
(cf. \cite{Fuchs:1958,springerlink:10.1007/BF02036756,Ohnishi:1952}). So
it is convenient to expect also characterisations of linear orders
on groups based on those which have been found for finitely generated
groups.

The approach we use here has been developed independently from the
mentioned articles. With the results of the preceding sections, we are
able to describe the linear orders on an abelian group with the help
of linearly ordered hyperplanes. Given an abelian group $\gruppe
=\Gruppe{G}{+}{-}{0}$ and a maximal independent set $E\subseteq G$,
each linearly ordered base in $\direktsumme{E}\reellezahlen$ defines
an order on \gruppe, and each linear order on the vector space
$\direktsumme{E}\reellezahlen$ can be reached that way. Linearly
ordered bases define systems of hyperplanes that are linearly ordered
by set inclusion. Each of this systems defines a different order on
$\direktsumme{E}\reellezahlen$. This is assured by Lemma
\ref{lemma:hyperebene-definiert-ordnung}, Theorem
\ref{satz:poskegelrand=Hyperebene}, and the monomorphisms in the
Diagram \ref{subfig:gzqr-monomorphismen2}. Theorem
\ref{lemma:Z-vektorraum} and Lemma \ref{lemma:ordungsuebertragung-G-Q}
assure that each linear order on the group can be described as linear
order on the real vector space $\direktsumme{E}\reellezahlen$.

So there is (among others) one question remaining: Which linearly ordered bases
define the same linear order on the Group $\gruppe$? We will address
this question in the current section.

A subspace $\vektorraum[U]$ of $\direktsumme{E}\reellezahlen$ may be
multidimensional, though it contains no rational vectors
($\vektorraum[U]∩\direktsumme{E}\rationalezahlen =\menge 0$). On the
other hand we know that the linear closure of the rational vectors is
the complete vector space. Thus, omitting a vector from a basis leads
to loss of rational vectors. This means, that every vector in a
linearly ordered base is necessary in order to define a linear order
on the rational „subspace“. In order to describe this we will use the
following notions.

\begin{definition}
  Let $\paar\basis{≤}$ be a linearly ordered basis of a vector space
  $\direktsumme E\reellezahlen$. A vector $\vektor x$ is called
  \defindex{(rationally) active} in a subset $U⊆\basis$ iff the linear
  closures of $U$ and $U\setminus\menge{\vektor x}$ differ in the
  rational vectors, i.e.
  \begin{equation}
    \linhuelle(U\setminus\menge{\vektor x})∩\direktsumme E\rationalezahlen ≠ 
    \linhuelle U∩\direktsumme{E}\rationalezahlen.
  \end{equation}
  Otherwise $\vektor x$ is called \defindex{(rationally) passive}.

  An element $\vektor a∈\basis$ is called \defindex{activator} of
  another element $\vektor x∈\basis{}$, iff both vectors $\vektor a$
  and $\vektor x$ are active in the principal ideal $\ordnungsideal
  {\vektor a}$. If $\vektor x$ is its own activator it is called
  \defindex{self-activator}. The set of activators of an element
  $\vektor x$ is denoted by $\Act \vektor x$.

  Furthermore let $\vektorraum$ be a subspace of
  $\direktsumme{E}\reellezahlen$. Then the set
  $\vAct \vektorraum := \linhuelle(\direktsumme E\rationalezahlen∩\vektorraum)$ is called
  the \defindex{active subspace} of $\vektorraum$.
\end{definition}

\begin{corollary}
  Let $\Halbord\basis{≤}$ be a linearly ordered basis of
  $\direktsumme{E}\reellezahlen$ and $U⊆\basis$ a corresponding
  subset. Then the active subspace of $\linhuelle U$ is the linear
  hull of the active vectors from $U$ in $\linhuelle U$.
\end{corollary}

\begin{lemma}
  Let $\Halbord\basis{≤}$ the linearly ordered basis of the linearly
  ordered vector space $\direktsumme{E}\reellezahlen$. Then the set of
  activators of an element $\vektor x$ is an order filter in the set
  of self-activators of $\basis$.
\end{lemma}

\begin{proof}
  By definition, a vector $\vektor a$ is an activator of $\vektor x$ iff it is a
  self-activator and the condition 
  \[
    \linhuelle(\ordnungsideal[\Halbord\basis{≤}]{\vektor a}\setminus\menge{\vektor x})∩\direktsumme E\rationalezahlen ≠ 
    \linhuelle(\ordnungsideal[\Halbord\basis{≤}]{\vektor a})∩\direktsumme{E}\rationalezahlen.
  \]
  holds. The latter condition is equivalent to the existence of a
  vector $\vektor b∈\vektorraum_{\vektor
    a}∩\direktsumme{E}\rationalezahlen$ whose projection along
  $\basis$ into $\linhuelle \vektor x$ is non-zero. Such a vector
  exists, iff there exists a basis vector $\vektor c$ such that
  $\vektorraum_{\vektor c} = \vektorraum_{\vektor b}$ and $\vektor c$
  is an activator of $\vektor x$. Since in linearly ordered vector
  spaces $\vektor b∈\vektorraum_{\vektor a}$ is the same as
  $\vektorraum_{\vektor b}⊆\vektorraum_{\vektor a}$, we know that
  $\vektor c≤\vektor a$. The vector $\vektor b$ also proves that this
  basis vector $\vektor c$ is also an activator of $\vektor x$ as
  $\vektor b\not∈\hyperebene_{\vektor b}=\hyperebene_{\vektor c}$. On
  the other hand for each activator $\vektor c≤\vektor a$ of $\vektor
  x$ we can find a vector $\vektor b$ whose projection into
  $\linhuelle{\vektor x}$ is non-zero. Thus we have proved so far: A
  self-activator $\vektor a∈\basis$ is an activator of $\vektor x$ iff
  there exists a self-activator $\vektor c≤\vektor a$ that is also an
  activator of $\vektor x$. This implies that for each activator
  $\vektor c∈\Act \vektor x$ of $\vektor x$ each self-activator $\vektor
  a≥\vektor c$ is also an activator of $\vektor x$.
\end{proof}

\begin{corollary}
  For every basis vector $\vektor x$ and every activator $\vektor y∈\Act
  \vektor x$ the relation \mbox{$\Act \vektor y ⊆ \Act \vektor x$} holds.
\end{corollary}
\iffalse
\begin{proof}
  Let $\vektor y∈\Act\vektor x$. Then there exists one vector $\vektor
  a∈(\linhuelle(\ordnungsideal{\vektor y})\setminus
  \linhuelle(\ordnungsideal{\vektor y}\setminus\menge{\vektor
  x}))∩\direktsumme{E}\rationalezahlen$. This vector has a unique
  coordinate description corresponding to $\basis$. By definiton the
  coordinate that corresponds to $\vektor x$ is nonzero. As the
  description is unique, this vector is not contained in any subspace
  that is spanned by a set of basis vectors that do not contain
  $\vektor x$. On the other hand for $\vektor b∈\Act \vektor y$ the inequality $\linhuelle(\ordnungsideal{\vektor b})≠
  \linhuelle(\ordnungsideal{\vektor b}\setminus\menge{\vektor
  y}))$ implies $\vektor
  y ≤ \vektor b$. Thus, the principal ideal
  $\ordnungsideal{\vektor b}$ is a superset of $\ordnungsideal{\vektor y}$, which implies $\vektor
  x∈\linhuelle{\ordnungsideal{\vektor b}}$ but $\vektor
  x\not∈\linhuelle(\ordnungsideal{\vektor b}\setminus\menge{\vektor x})$,
  which implies: $\vektor b∈\Act{\vektor x}$ iff $\vektor b∈\Act{\vektor b}$.
\end{proof}
\fi

\noindent We consider two linearly ordered bases equivalent iff they
generate the same order, i.e. their primary ideals generate the same
subspaces $\direktsumme{E}\rationalezahlen$.

\begin{lemma}
  For every linearly ordered basis $\basis$ there exists an equivalent
  basis $\basis'$ such that for each basis vector $\vektor x∈\basis'$
  the condition
  \begin{equation}
    \vektor x∈\Act' \vektor x⇔\vektor x∈\direktsumme E\rationalezahlen
  \end{equation}
  holds.
\end{lemma}
\begin{proof}
  Let $\vektor x∈\basis$ be a self-activator ($\vektor x∈\Act'\vektor
  x$). We know by definition that
  $\linhuelle\bigl(\linhuelle(\ordnungsideal{\vektor
    x})∩\direktsumme{E}\rationalezahlen\bigr)=\linhuelle(\vektorraum_{\vektor
    x}∩\direktsumme{E}\rationalezahlen)$ and
  $\linhuelle(\ordnungsideal{\vektor x}\setminus\vektor
  x)∩\direktsumme{E}\rationalezahlen = \hyperebene_{\vektor
    x}∩\direktsumme{E}\rationalezahlen$ hold. By the presumption we
  know $(\vektorraum_{\vektor x}\setminus\hyperebene_{\vektor
    x})∩\direktsumme{E}\rationalezahlen ≠∅$. Let us choose a vector
  $f(\vektor x) ∈ (\vektorraum_{\vektor
    x}\setminus\hyperebene_{\vektor
    x})∩\direktsumme{E}\rationalezahlen$ for every self-activator
  $\vektor x∈\basis$. In that case we know $\vektorraum_{f(\vektor x)
  } = \vektorraum_{\vektor x}$. Otherwise choose $f(\vektor x) =
  \vektor x$. Then $f$ is a mapping from $\basis$ into the vector
  space $\direktsumme{E}\reellezahlen$. And $\basis':=\Menge{f(\vektor
    x)}{\vektor x∈\basis}$ fulfils the condition if $f(\vektor
  x)≤'f(\vektor y) :⇔ \vektor x≤\vektor y$.
\end{proof}

\noindent Thus, we can develop a standard form for linearly ordered bases.

\begin{definition}
  A linearly ordered base $\basis$ is called \defindex{clarified}, iff it fulfils the condition 
  \begin{equation}
    \vektor x∈\Act' \vektor x⇔\vektor x∈\direktsumme E\rationalezahlen.\qedhere
  \end{equation}
\end{definition}

\noindent As we can describe every linear order by a clarified base, in the ramaining part
of this section, we will consider all bases to be clarified if it is
not stated otherwise.

\begin{lemma}\label{lemma:projektion-erhaelt-dimension}
  Let $\basis$ be a base of $\direktsumme{E}\reellezahlen$ and
  $\vektor a,\vektor b∈\basis$ be two different basis vectors. Then the
  projection of $\direktsumme{E}\rationalezahlen$ along
  $\basis\setminus\menge{\vektor a,\vektor b}$ into
  $\linhuelle({\vektor a,\vektor b})$ is 2-dimensional.
\end{lemma}

\begin{proof}
  Since both vectors $\vektor a$ and $\vektor b$ are linear
  combinations of the standard basis, the projection of the standard
  base vectors contains a nontrivial vector $\vektor z$. Let $P$ be
  the projection along $\basis\setminus\menge{\vektor a,\vektor
    b}$. Then from $\vektor z∈P[\direktsumme{E}\rationalezahlen]$
  follows the set inclusion $\linhuelle\menge{\vektor
    z}⊆\linhuelle(P[\direktsumme{E}\rationalezahlen])⊆\linhuelle\menge{\vektor
    a,\vektor b}$. 

  Suppose $\linhuelle\menge{\vektor
    z}=\linhuelle(P[\direktsumme{E}\rationalezahlen])$. Then we can
  express every rational vector as a linear combination of $\vektor z$ and
  $\basis\setminus \menge{\vektor a,\vektor b}$. This leads to the set inclusion 
  \[
  \direktsumme{E}\rationalezahlen⊆
  \linhuelle\bigl((\basis\setminus\menge{\vektor a,\vektor b})∪
  \menge{\vektor z}\bigr)⊆\linhuelle\basis.
  \]
  As the linear hull is a closure operator, we get immediately the set inclusion
  \[
  \linhuelle(\direktsumme{E}\rationalezahlen)⊆
  \linhuelle\bigl((\basis\setminus\menge{\vektor a,\vektor b})∪\menge{\vektor z}\bigr)⊆
  \linhuelle\basis=\linhuelle(\direktsumme{E}\rationalezahlen).
  \]

  As $\vektor z$ is
  linear dependend from the vectors $\vektor a$ and $\vektor b$, at
  least one of the equations $\linhuelle \menge{\vektor a,\vektor b} =
  \linhuelle\menge{\vektor a,\vektor z}$ and $\linhuelle
  \menge{\vektor a,\vektor b} = \linhuelle\menge{\vektor z,\vektor b}$
  holds. W.l.o.g. we can assume $\linhuelle \menge{\vektor a,\vektor
    b} = \linhuelle\menge{\vektor a,\vektor z}$. Then we can exchange
  $\vektor b$ in $\basis$ with $\vektor z$, too. We get
  \[
  \linhuelle(\direktsumme{E}\rationalezahlen)=\linhuelle(\direktsumme{E}\reellezahlen)
  =\linhuelle{\basis}=\linhuelle\bigl((\basis\setminus\menge{\vektor b})∪\menge{\vektor z}\bigr)
  \]
  As we have seen above, from this follows that $\vektor a$ is linear
  dependent from the set $(\basis \setminus\menge{\vektor a,\vektor
    b})∪\menge{\vektor z}$. Thus, we can express $\vektor a$ as a
  linear combination of the vectors from $(\basis
  \setminus\menge{\vektor a,\vektor b})∪\menge{\vektor z}$. As
  $\basis$ is a basis this representation has a nonzero $\vektor z$
  component and can express $\vektor z$ as a linear combination of the
  vectors from $\basis\setminus\menge{\vektor b}$. Consequently we get
  the following equations
  \[
  \linhuelle\bigl((\basis\setminus\menge{\vektor a,\vektor b})∪\menge{\vektor z}\bigr)
  =\linhuelle(\basis\setminus \menge{\vektor b})
  =\linhuelle(\direktsumme{E}\rationalezahlen)
  =\linhuelle(\basis).
  \]
  This shows that either $\basis$ is not a basis, that $\vektor a=\vektor b$, or that the assumption
  $\linhuelle\menge{\vektor
    z}=\linhuelle(P[\direktsumme{E}\rationalezahlen])$ is wrong.
\end{proof}

Using this lemma, we can find a first condition, which explains when
two different orders on $\direktsumme{E}\reellezahlen$ imply the same
order on $\direktsumme{E}\rationalezahlen$.

\begin{lemma}
  Let the vector space $\direktsumme{E}\reellezahlen$ be orderded
  twice with the linear order relations $≤$ and $≤'$. Furthermore let
  $\basis$ a Basis that describes both $≤$ and $≤'$, and $\Act$
  respectively $\Act'$ the corresponding operators that provide the
  set of activators of $≤$ respectively $≤'$. Then the induced orders
  on the rational vector space $\direktsumme{E}\rationalezahlen$
  coincide iff for evey element $\vektor x∈\basis$ the sets of
  Activators $\Act\vektor x$ and $\Act'\vektor x$ are equal.
\end{lemma}

\begin{proof}
  We have to consider four cases: 
  \begin{enumerate}
  \item\label{item:4}The two bases are equally ordered. In that case
    the orders also share the same sets of activators.
  \item\label{item:5}There exist two self-activators $\vektor x$ and $\vektor y$
    with $\vektor x≤\vektor y$ and $\vektor y≤'\vektor x$.
  \item\label{item:6}The sets of self-activators are different and the orders are
    the same on those elements which are self-activators according to
    both orders.
  \item\label{item:7}The sets of self-activators are the same and the orders differ only
    on non-self-activators.
  \end{enumerate}

  Let us consider point~\ref{item:5}. Then we can choose two vectors 
  \begin{align*}
    \vektor a&∈(\vektorraum_{\vektor x}\setminus\hyperebene_{\vektor x})∩\direktsumme{E}\rationalezahlen\\
    \vektor b&∈(\vektorraum'_{\vektor y}\setminus\hyperebene'_{\vektor y})∩\direktsumme{E}\rationalezahlen
  \end{align*}
  such that the component of $\vektor a$ in $\vektor x$ direction and the
  component of $\vektor b$ in $\vektor y$ direction are non-zero. Then
  $\vektor a<\vektor b<'\vektor a$ holds.

  In case~\ref{item:6} a vector $\vektor x∈\basis$ exists such that $\vektor
  x∈\Act \vektor x$ and $\vektor x\not∈\Act' \vektor x$ hold. Then we
  choose a rational vector $\vektor a∈(\vektorraum_{\vektor
    x}\setminus\hyperebene_{\vektor
    x})∩\direktsumme{E}\rationalezahlen$ such that the component in
  $\vektor x$ direction is different from zero. Let us choose $\vektor
  y∈\basis$ such that $\vektorraum'_{\vektor y}=\vektorraum'_{\vektor
    a}$. As $\vektor a$ has a non-zero $\vektor y$ component,
  according to $\basis$, we have the inequality $\vektorraum_{\vektor
    y}⊆\vektorraum_{\vektor x}$ leading to $\Act\vektor x⊆\Act\vektor
  y$. As it has a non-zero $\vektor x$ component, too, we get
  $\Act'\vektor y⊆\Act'\vektor x$. By
  Lemma~\ref{lemma:projektion-erhaelt-dimension} we can find a vector
  $\vektor b ∈\vektorraum'_{\vektor
    y}∩\direktsumme{E}\rationalezahlen$ whose $\vektor x$ component is
  larger than the one of $\vektor a$ and whose $\vektor y$ component
  is smaller than the one of $\vektor a$. Then $\vektor a<\vektor
  b<'\vektor a$ proves the assertion.

  For the remaining part \ref{item:7} we prove that two inactive basis
  vectors $\vektor x,\vektor y∈\basis$ do not influence the
  result. Suppose the relations $\vektor x\not∈\Act'\vektor x$,
  $\vektor x\not∈\Act\vektor x$, and $\vektor x<\vektor y<'\vektor
  x$. Suppose there are two vectors $\vektor a,\vektor
  b∈\direktsumme{E}\vektorraum$ with $\vektor a<\vektor b$. As all
  vectors have finite support with respect to the standard base, their
  support with respect to $\basis$ is finite, too. Thus, there exists
  a largest basis vector $\vektor z∈\basis$ such that the $\vektor z$
  component of $\vektor a$ is smaller than the one of $\vektor
  b$. Then both vectors $\vektor a,\vektor b∈\vektorraum_{\vektor z}$
  belong to its vector space, while at least one of them is not in the
  corresponding hyperplane $\hyperebene_{\vektor z}$ leading to $\vektor
  z∈\Act\vektor z$. Thus, $\vektor x≠\vektor z$. Suppose $\vektor
  x<\vektor z <'\vektor x$ and $\vektor a$ has a non-zero $\vektor x$
  component. Then there would be a self-activator $\vektor z'∈\Act' \vektor
  z'∩\Act' \vektor x$ such that $\vektor a$ has a non-zero $\vektor
  z'$ component. As $\vektor a∈\vektorraum_{\vektor z}$ this implies
  $\vektor z'<\vektor z<'\vektor z'$ which has been excluded by the
  presumption. Consequently, either $\vektor a$ has a zero $\vektor x$
  component or $\vektor x<\vektor z ⇔ \vektor x<'\vektor z$ holds. The
  same is true for any other non-self-activator. As the orders on the
  self-activators coincide between the relations $≤$ and $≤'$, the
  $\vektor z$ component remains the main component to determine the
  order of $\vektor a$ and $\vektor b$ leading to $\vektor a≤'\vektor
  b$.
\end{proof}

\noindent In Theorem \ref{satz:poskegelrand=Hyperebene} we have shown that in
$\direktsumme{E}\reellezahlen$ each linearly ordered hyperplane
defines a different linear order. On the other hand it is not obvious,
which linearly ordered hyperplanes define different linear orders on
$\direktsumme{E}\rationalezahlen$.

\begin{lemma}\label{lemma:R-Ordnung-induziert-Z-ordnung-eindeutig}
  Let $\leq$ and $\leq'$ be two linear orders on\/
  $\direktsumme{E}\reellezahlen$ and $\vektor
  a∈\direktsumme{E}\rationalezahlen$ a rational vector. Furthermore we
  denote the vector spaces from Equation~(\ref{eq:Va}) with
  $\vektorraum_{\vektor a}$ and $\vektorraum_{\vektor a}'$ and the
  corresponding linearly ordered hyperplanes according to
  Equation~(\ref{eq:Ha}) with $\hyperebene_{\vektor a}$ and
  $\hyperebene_{\vektor a}'$. If the active subspaces
  $\vAct\vektorraum_{\vektor a}$ and $\vAct\vektorraum'_{\vektor a}$
  are not equal or the restrictions of the hyperplanes
  $\hyperebene_{\vektor a}∩\vAct\vektorraum_{\vektor a}$ and
  $\hyperebene_{\vektor a}'∩\vAct\vektorraum_{\vektor a}'$ to the
  active subspaces $\vektorraum_{\vektor a}$ and $\vektorraum_{\vektor
    a}'$ differ, then the induced orders w.r.t. Lemma
  \ref{lemma:monomorphismusordnung} differ on
  $\direktsumme{E}\rationalezahlen$.
\end{lemma}
\begin{proof}
  Suppose at first $\vAct\vektorraum_{\vektor a}=\vAct\vektorraum_{\vektor a}'$.
  According to the presumption the set 
  \[U:=\poskegel{(\vAct\vektorraum_{\vektor
      a})}\setminus\poskegel{\vektorraum_{\vektor
      a}'}=\poskegel{(\vAct\vektorraum_{\vektor a})}\cap
  \negkegel{\vektorraum_{\vektor
      a}'}\setminus\poskegel{\hyperebene_{\vektor a}'}
  =\poskegel{(\vAct \vektorraum_{\vektor a})}∩\negkegel{(\vektorraum'_{\vektor a})}\setminus\menge{\nullvektor}
  \]
  
  is non-empty. Let $\basis$ the linearly ordered base that describes
  the order relation $≤$ and $\basis'$ the one of $≤'$. As the
  hyperplanes $\hyperebene_{\vektor a}$ and $\hyperebene'_{\vektor a}$
  differ, there exists a base vector $\vektor
  b∈\basis∩\hyperebene_{\vektor a}∩\vAct\vektorraum_{\vektor
    a}\setminus\hyperebene_{\vektor a}'$ and a base vector $\vektor
  b'∈\basis'∩\hyperebene'_{\vektor a}∩\vAct\vektorraum'_{\vektor
    a}\setminus\hyperebene_{\vektor a}$. Both vectors are linear
  independent. As $\vektor b'\not∈\hyperebene_{\vektor a}$, for the
  base vector $\vektor y∈\basis$ with $\vektorraum_{\vektor
    y}=\vektorraum_{\vektor a}$ the set
  $\hat\basis=\basis\setminus\menge{\vektor y}∪\menge{\vektor b'}$ is
  another basis that describes $≤$. As $\vektor a∈\direktsumme
  E\rationalezahlen$ the vector $\vektor b'$ is active in
  $\vektorraum_{\vektor a}$, by
  Lemma~\ref{lemma:projektion-erhaelt-dimension} the projection of
  $\direktsumme{E}\rationalezahlen∩\vektorraum_{\vektor a}$ along
  $\hat\basis$ into $\linhuelle\menge{\vektor b,\vektor b'}$ is
  two-dimensional. Thus it contains for all 4 combinations of positive
  or negative (non-zero) scalars with respect to $\vektor b$ and
  $\vektor b'$ at least one element. This implies that at least one
  element exists that is positive with respect to both orders, one that is
  positive with respect to $≤$ and negative with respect to $≤'$, one
  that is negative with respect to both orders and at least one that
  is negative with respect to $≤$ and positive with respect to
  $≤'$. The latter one is sufficient to prove that the orders differ.

  The remaining part of the proof considers $\vAct\vektorraum_{\vektor
    a}≠\vAct\vektorraum_{\vektor a}'$. W.l.o.g. considering the
  definition of $\vAct$ there exists a vector $\vektor
  b∈\vAct\vektorraum_{\vektor
    a}∩\direktsumme{E}\rationalezahlen\setminus\vAct\vektorraum_{\vektor
    a}'$. Furthermore we can choose an integer $n∈\ganzzahlen$ such
  that $\vektor b<\vektor a^n$ and as $\vektor
  b\not\in\vektorraum_{\vektor a}'$ we get the other inequality
  $\vektor a \infinitesimal{<}'\vektor b$.
\end{proof}

\noindent After having shown that different hyperplanes imply different orders,
the next question we have to address, considers the effect of having different orders
on $\direktsumme{E}\reellezahlen$, which share a common dividing hyperplane.
 \begin{lemma}\label{lemma:Ordnungen-definieren-hyperebene-eindeutig}
   Let $\leq$ and $\leq'$ be two differing linear orders on
   $\direktsumme{E}\reellezahlen$, which have the same hyperplane
   $\hyperebene_{\vektor a}=\hyperebene'_{\vektor a}$ for some rational vector
   $\vektor a∈\direktsumme{E}\rationalezahlen$. If the vector
   $\vektor a$ is positive w.r.t.\ each of the orders
   ($\nullvektor \leq\vektor a$ and $\nullvektor \leq'\vektor a$),
   then the induced orders on $\vektorraum_{\vektor
     a}∩\direktsumme{E}\rationalezahlen$ are identical iff the induced
   orders on $\hyperebene_{\vektor
     a}\cap\direktsumme{E}\rationalezahlen$ and $\hyperebene'_{\vektor
     a}\cap\direktsumme{E}\rationalezahlen$ are equal.
 \end{lemma}
 \begin{proof}
   According to Lemma  \ref{lemma:monomorphismusordnung} all elements
   of $\direktsumme{E}\rationalezahlen\setminus\hyperebene$ and all
   elements of $\hyperebene\cap\direktsumme{E}\rationalezahlen$ are either
   positive with respect to both orders on
   $\direktsumme{E}\rationalezahlen$ or negative. Thus, the order on
   $\direktsumme{E}\reellezahlen\setminus\direktsumme{E}\rationalezahlen$
   does not influence the order on $\direktsumme{E}\rationalezahlen$. 
 \end{proof}

 \begin{lemma}\label{folg:lin-ordnung-auf-Q-Z-aequivalent}
   Let $≤$ and $≤'$ be two linear orders on
   $\direktsumme{E}\reellezahlen$. They are equal on
   $\direktsumme{E}\ganzzahlen$ iff the coincide on
   $\direktsumme{E}\rationalezahlen$.
 \end{lemma}

 \begin{proof}
   As $\direktsumme{E}\ganzzahlen⊆\direktsumme{E}\rationalezahlen$ the
   orders are equal on $\direktsumme{E}\ganzzahlen$ if they are equal
   on $\direktsumme{E}\rationalezahlen$. For the other direction let
   $\vektor a∈\direktsumme{E}\rationalezahlen$ such that $\vektor
   a<\nullvektor <'\vektor a$, i.e. it is positive with respect to one
   and negative with respect to the other order. Then it has finitely
   many rational coordinates with respect to the standard base of
   $\direktsumme{E}\rationalezahlen$. Let $a∈\natzahlen$ be the least
   common denominator of the coordinates of $\vektor
   a$. Then $a\vektor a∈\direktsumme{E}\ganzzahlen$ is an integer
   vector. Thus we get the inequality $a\vektor a<\nullvektor
   <'a\vektor a$.
 \end{proof}

\noindent With Corollary \ref{folg:lin-ordnung-auf-Q-Z-aequivalent} and the
monomorphisms $φ$ and $ψ$ as defined in Theorems \ref{satz:Monomorphismus-G-Q}
and \ref{lemma:Monomorphismus-Z-G} we also have a characterisation of
linear orders on an arbitrary abelian group:

\begin{theorem}[Characterisation of linear orders on abelian groups]\label{satz:charakterisierung-linearer-ordnungen}\hfill\\
  Let $\gruppe=\Gruppe{G}{+}{-}{0}$ be an abelian group, $E\subseteq
  G$ a maximal independent set in \gruppe{}, and $φ$ the canonical
  embedding of $\gruppe$ into $\direktsumme{E}\reellezahlen$ as
  described in Lemma \ref{satz:Monomorphismus-G-Q}. Let further
  $U⊆\direktsumme{E}\rationalezahlen$ be a set that generates
  $\direktsumme E\reellezahlen$ in the following way:
  \begin{equation}
    \direktsumme{E}\reellezahlen = ⋃_{\vektor a∈U}\vektorraum_{\vektor a}.
  \end{equation}
  Two linear orders $\leq$ and $\leq'$ on
  \gruppe{} are different iff some vector $\vektor a∈U$ fulfils one of the following
  conditions:
  \begin{enumerate}
  \item $\hyperebene_{\vektor a}∩\vAct\vektorraum_{\vektor a}
    \neq\hyperebene'_{\vektor a}∩\vAct\vektorraum'_{\vektor a}$\label{item:8}
  \item $\vAct\vektorraum_{\vektor a}\neq\vAct\vektorraum'_{\vektor a}$\label{item:9}
  \item $\poskegel{(\hyperebene_{\vektor a})}\cap\direktsumme{E}\rationalezahlen\neq\poskegel{(\hyperebene'_{\vektor a})}\cap\direktsumme{E}\rationalezahlen$.\label{item:10}
  \end{enumerate}
\end{theorem}
\begin{proof}
  As shown in the first part
  of this article Diagram \ref{subfig:gzqr-monomorphismen2} commutes
  and with Corollary~\ref{folg:lin-ordnung-auf-Q-Z-aequivalent} each
  linear order on $\direktsumme{E}\rationalezahlen$ uniquely defines
  an order on the group \gruppe{}. Lemma
  \ref{lemma:ordungsuebertragung-G-Q} finally assures that there are
  no more linear orders on \gruppe{}, since each linear order is
  semiclosed.
  
  Lemma \ref{lemma:monomorphismusordnung} tells us that we can
  transfer the order from $\direktsumme{E}\reellezahlen$ onto
  $\direktsumme{E}\rationalezahlen$ and
  $\direktsumme{E}\ganzzahlen$. The Theorems
  \ref{satz:Monomorphismus-G-Q} and \ref{lemma:Monomorphismus-Z-G}
  assure that these orders can also be considered as orders on the
  group \gruppe. Let $\vektor a∈U$ one of the selected vectors. As
  $\vektorraum_{\vektor a}=\hyperebene_{\vektor a}\vektorraumplus
  \reellezahlen\vektor a$ and similar for $≤'$,
  Theorem~\ref{satz:poskegelrand=Hyperebene} describes the positive
  cone of $\vektorraum_{\vektor a}$ using a
  hyperplane. Lemma~\ref{lemma:R-Ordnung-induziert-Z-ordnung-eindeutig}
  shows that Conditions~\ref{item:8} and \ref{item:9} indicate a
  difference between the orders $≤$ and $≤'$. If these conditions do
  not hold,
  Lemma~\ref{lemma:Ordnungen-definieren-hyperebene-eindeutig} proves
  Condition~\ref{item:10}. Thus the restrictions of the orders $≤$ and
  $≤'$ to $\vektorraum_{\vektor a}∪\vektorraum'_{\vektor a}$ are equal
  iff none of these three conditions holds.
  
  From $\direktsumme{E}\reellezahlen = ⋃_{\vektor
    a∈U}\vektorraum_{\vektor a}= \vAct\direktsumme{E}\reellezahlen$ also the restricted equation
  $\direktsumme E\rationalezahlen =
  ⋃_{\vektor a∈U}\vektorraum_{\vektor a}∩\direktsumme{E}\rationalezahlen$ follows as well as the
  same for the other order relation $≤$. Thus, the orders $≤$ and $≤'$
  are equal on $\direktsumme{E}\rationalezahlen$ iff their
  restrictions to $\vAct \vektorraum_{\vektor a}$ are equal for each
  vektor $\vektor a∈U$.
\end{proof}
\begin{corollary}
  Linear orders on groups can be completely described by linearly
  ordered hyperplanes and their intersections in certain subspaces of a real vector
  space $\direktsumme{E}\reellezahlen$.
\end{corollary}

\noindent Theorem \ref{satz:charakterisierung-linearer-ordnungen}
provides a characterisation of linearly ordered groups independent of
the ideas described in The's~\cite{Teh:1961} and
Zajceva's~\cite{M.I.Zajceva:1953} papers. Teh describes linear orders
by embedding them into vector spaces of the form $\direktsumme{E}\reellezahlen$
where the cardinality $|E|$ is smaller or equal to that of the continuum $c$.
%\footnote{Recall,
 % $\direktsumme{M}\reellezahlen\subseteq\reellezahlen^M$ where the
%  inclusion can be a proper inclusion if the set $M$ is infinite.}
For the cardinality of these subsets he defines a prototype of an
order and considers all other orderings to be homomorphic images of
the prototypic vector spaces. The dimensions of these spaces
correspond to the archimedian rank of the generated order in
\gruppe{}. Our characterisation does not touch the notion of
archimedian orders. The description using linearly ordered bases
provides additional details about in the case of infinite archimedian
rank. The work of Zajceva describes archimedian orders by certain
equations. It is not evident that finite linear equations can be used
in general for infinitely generated groups. Nevertheless, these
equations can be used in several cases since they are descriptions of
the hyperplanes, which we used in Theorem
\ref{satz:charakterisierung-linearer-ordnungen}. So far, we didn't use
the properties of $\direktsumme{E}\reellezahlen$ as a scalar product
space. This will be necessary for discussing the plane equations.

Beforehand we head over to linear order extensions, we can shortly
bridge to Zajceva's results. Two elements $g,h\in G$ of an abelian
linearly ordered group $\gruppe=\OrdGruppe{G}{+}{-}{0}{\leq}$ with
$g\leq h$ are considered \defindex{archimedian equivalent} iff there
exist two integers $m,n\in\ganzzahlen$ such that $h^m\leq g\leq h\leq
g^n$ holds. The group $\gruppe$ is archimedian iff all of its elements
are archimedian equivalent to each other.  In other words: For each
two elements $g,h\in G$ with $0\leq h-g$ there exists a number
$n\in\ganzzahlen$ such that $h-ng\leq 0$. When we describe the linear
order on $\gruppe$ by a hyperplane as in Theorem
\ref{satz:charakterisierung-linearer-ordnungen}, this is equivalent
to: for all vectors $\vektor{x},\vektor{y}\in φ[G]$ with
$\vektor{y}-\vektor{x}\in
\poskegel{\hyperebene}\cup\hyperebene+\bigl(\poskegel{\reellezahlen}\setminus\menge{0}\bigr)\vektor{a}$
there exists a number $n\in\ganzzahlen$ such that
$\vektor{y}-n\vektor{x}\in
\negkegel{\hyperebene}\cup\hyperebene+\bigl(\negkegel{\reellezahlen}\setminus\menge{0}\bigr)\vektor{a}$. If
the intersection $\hyperebene\cap φ[G]$ of the dividing hyperplane of
the order and the image of the group is empty, then the vector
$\vektor{x}$ is linearly independent from \hyperebene{} and, thus, the
linear closure $\linhuelle{(\hyperebene\cup\menge{\vektor{x}})} =
\direktsumme{E}\reellezahlen$ is the whole vector space. Thus, there
exists an integer $n$ such that $\vektor{y}-n\vektor{x}$ is an element
of the negative half space. On the other hand if
$\vektor{x}\in\hyperebene$ is an element of the hyperplane and
$\vektor{y}\not \in\hyperebene$ is not, then the straight line
$\vektor{y}+\reellezahlen\vektor{x}$ is parallel to \hyperebene{} an
thus $\vektor{y}-n\vektor{x}$ will never be negative. Thus, the group
is archimedian iff there is no element $g\in G$ such that
$φ(g)\in\hyperebene$. According to Corollary
\ref{folg:lin-ordnung-auf-Q-Z-aequivalent} this is equivalent to the
condition that no vector in \hyperebene{} has only rational
coordinates with respect to the standard base defined by the
independent set $E$ as described in the preceding sections. This
implies that the quotient of any two coordinates of the normal vector
of $\hyperebene$ must be irrational if it exists. Otherwise we can
scale the normal vector by a real number such that it has two rational
coordinates.  On the other hand if all such quotients are irrational,
all vectors containing only rational coordinates are never orthogonal
to the normal vector.

In \cite{M.I.Zajceva:1953} Zajceva discusses linear orders on abelian
groups with a maximal independent set $E$ by assigning a set of
rationally independent numbers $E'$ with a bijective mapping
$\wabbildung{\abbildung{α}{E}{E'}}{e}{α(e)}$ to them. So each element
$g\in G$ with a rational representation $g=\sum_{a\in E}ξ_a(g)a$
($\abbildung{ξ_a}{G}{\rationalezahlen}$) is positive iff the real
number $\sum_{a\in E}ξ_a(g)α(a)$ is positive.  These equations and
inequalities can be considered as the scalar product of the normal
vector to the dividing hyperplane $\hyperebene$ in the case of abelian
groups with a finite maximal independent set. However if the maximal
independent set $E$ is infinite, the normal vector must have infinite
support and thus is not an element of the vector space. In that case
using our approach to find such a representation is less suitable than
Zajceva's method. It is again the hierarchy of subspaces generated by
the order, that extends this description.

\section{Linear order extensions}
Before we discuss linear order extensions on arbitrary partially
ordered torsion free abelian groups, we start with the special case of
ordered vector spaces. So far we have discussed
linearly ordered hyperplanes as representations of linear orders on
vector spaces. Now we can relate them to other orders.
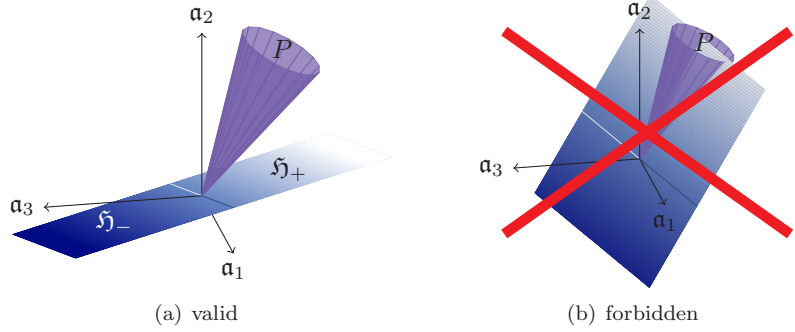
\begin{figure}[tbp]%
  \hbox to \linewidth{\hfill
  \subfigure[][valid]{% Sketch output, version 0.3 (build 7d, Fri Nov 22 14:14:09 2013)
% Output language: PGF/TikZ,LaTeX
\begin{tikzpicture}[line join=round]
\tikzstyle{ann} = [fill=white,font=\footnotesize,inner sep=1pt]
\providecommand\vektor{}\tikzstyle{plane}=[shading=axis,top color=white,bottom
color=HKS41K80,shading angle=-50]\tikzstyle{plane2}=[shading=axis,top color=white,bottom
color=HKS41K80,shading angle=-30]\draw[arrows=<-,line width=.4pt](.42,-.767)--(0,0);
\filldraw[plane,draw=black,draw opacity=0.01,line width=0.1,fill opacity=1](-1.68,-.833)--(2.52,.513)--(1.68,.833)--(-2.52,-.513)--cycle;
\draw[draw opacity=0.8,draw=white](-.42,.16)--(0,0);
\draw[arrows=->,line width=.4pt](0,0)--(0,2.167);
\draw[draw opacity=0.8,draw=HKS41K80](.42,-.16)--(0,0);
\draw[arrows=->,line width=.4pt](0,0)--(-2.1,-.153);
\filldraw[draw=HKS36K80,draw opacity=0.7,line width=0.0001,fill=HKS36K80,fill opacity=0.7](0,0)--(.548,2.207)--(.624,2.247)--(0,0)--cycle;
\filldraw[draw=HKS36K80,draw opacity=0.7,line width=0.0001,fill=HKS36K80,fill opacity=0.7](0,0)--(.624,2.247)--(.765,2.239)--(0,0)--cycle;
\filldraw[draw=HKS36K80,draw opacity=0.7,line width=0.0001,fill=HKS36K80,fill opacity=0.7](0,0)--(.765,2.239)--(.951,2.183)--(0,0)--cycle;
\filldraw[draw=HKS36K80,draw opacity=0.7,line width=0.0001,fill=HKS36K80,fill opacity=0.7](0,0)--(.951,2.183)--(1.152,2.087)--(0,0)--cycle;
\filldraw[draw=HKS36K80,draw opacity=0.7,line width=0.0001,fill=HKS36K80,fill opacity=0.7](0,0)--(1.152,2.087)--(1.339,1.968)--(0,0)--cycle;
\filldraw[draw=HKS36K80,draw opacity=0.7,line width=0.0001,fill=HKS36K80,fill opacity=0.7](0,0)--(1.339,1.968)--(1.482,1.842)--(0,0)--cycle;
\filldraw[draw=HKS36K80,draw opacity=0.7,line width=0.0001,fill=HKS36K80,fill opacity=0.7](0,0)--(1.482,1.842)--(1.561,1.728)--(0,0)--cycle;
\filldraw[draw=HKS36K80,draw opacity=0.7,line width=0.0001,fill=HKS36K80,fill opacity=0.7](0,0)--(1.561,1.728)--(1.562,1.645)--(0,0)--cycle;
\filldraw[draw=HKS36K80,draw opacity=0.7,line width=0.0001,fill=HKS36K80,fill opacity=0.7](0,0)--(1.562,1.645)--(1.486,1.605)--(0,0)--cycle;
\filldraw[draw=HKS36K80,draw opacity=0.7,line width=0.0001,fill=HKS36K80,fill opacity=0.7](0,0)--(1.486,1.605)--(1.344,1.613)--(0,0)--cycle;
\filldraw[draw=HKS36K80,draw opacity=0.7,line width=0.0001,fill=HKS36K80,fill opacity=0.7](0,0)--(1.344,1.613)--(1.159,1.669)--(0,0)--cycle;
\filldraw[draw=HKS36K80,draw opacity=0.7,line width=0.0001,fill=HKS36K80,fill opacity=0.7](0,0)--(1.159,1.669)--(.957,1.764)--(0,0)--cycle;
\filldraw[draw=HKS36K80,draw opacity=0.7,line width=0.0001,fill=HKS36K80,fill opacity=0.7](0,0)--(.957,1.764)--(.771,1.884)--(0,0)--cycle;
\filldraw[draw=HKS36K80,draw opacity=0.7,line width=0.0001,fill=HKS36K80,fill opacity=0.7](0,0)--(.771,1.884)--(.627,2.01)--(0,0)--cycle;
\filldraw[draw=HKS36K80,draw opacity=0.7,line width=0.0001,fill=HKS36K80,fill opacity=0.7](0,0)--(.627,2.01)--(.549,2.123)--(0,0)--cycle;
\filldraw[draw=HKS36K80,draw opacity=0.7,line width=0.0001,fill=HKS36K80,fill opacity=0.7](0,0)--(.549,2.123)--(.548,2.207)--(0,0)--cycle;
\path (1.068,1.95) node {$P$};\path (.42,-.767) node[below] {$\vektor a_1$}
     (0,2.167) node[above] {$\vektor a_2$}
     (-2.1,-.153) node[left] {$\vektor a_3$};\path (1.159,.353) node {$\poskegel\hyperebene$};\path[white] (-1.159,-.353) node {$\negkegel\hyperebene$};\end{tikzpicture}% End sketch output
}\hfill
  \subfigure[][forbidden]{% Sketch output, version 0.3 (build 7d, Fri Nov 22 14:14:09 2013)
% Output language: PGF/TikZ,LaTeX
\begin{tikzpicture}[line join=round]
\tikzstyle{ann} = [fill=white,font=\footnotesize,inner sep=1pt]
\providecommand\vektor{}\tikzstyle{plane}=[shading=axis,top color=white,bottom
color=HKS41K80,shading angle=-50]\tikzstyle{plane2}=[shading=axis,top color=white,bottom
color=HKS41K80,shading angle=-30]\filldraw[draw=HKS36K80,draw opacity=0.7,line width=0.0001,fill=HKS36K80,fill opacity=0.7](0,0)--(1.249,1.383)--(1.25,1.316)--(0,0)--cycle;
\filldraw[draw=HKS36K80,draw opacity=0.7,line width=0.0001,fill=HKS36K80,fill opacity=0.7](0,0)--(1.186,1.473)--(1.249,1.383)--(0,0)--cycle;
\filldraw[fill=HKS36K80,fill opacity=0.7,draw=none](0,0)--(1.189,1.284)--(1.128,1.287)--(0,0)--cycle;
\draw[draw=HKS36K80,draw opacity=0.7,line width=0.0001](0,0)--(0,0)--(1.189,1.284)--(1.128,1.287);
\filldraw[draw=HKS36K80,draw opacity=0.7,line width=0.0001,fill=HKS36K80,fill opacity=0.7](0,0)--(1.25,1.316)--(1.189,1.284)--(0,0)--cycle;
\filldraw[draw=HKS36K80,draw opacity=0.7,line width=0.0001,fill=HKS36K80,fill opacity=0.7](0,0)--(1.071,1.574)--(1.186,1.473)--(0,0)--cycle;
\filldraw[draw=HKS36K80,draw opacity=0.7,line width=0.0001,fill=HKS36K80,fill opacity=0.7](0,0)--(.922,1.67)--(1.071,1.574)--(0,0)--cycle;
\filldraw[fill=HKS36K80,fill opacity=0.7,draw=none](0,0)--(.87,1.695)--(.922,1.67)--(0,0)--cycle;
\draw[draw=HKS36K80,draw opacity=0.7,line width=0.0001](.87,1.695)--(.922,1.67)--(0,0)--(0,0);
\filldraw[plane2,draw=black,draw opacity=0.01,line width=0.1,fill opacity=0.8](.13,-1.717)--(1.668,.921)--(.142,2.182)--(-1.397,-.455)--cycle;
\draw[draw opacity=0.8,draw=white](-.763,.631)--(0,0);
\filldraw[fill=HKS36K80,fill opacity=0.7,draw=none](0,0)--(.761,1.746)--(.87,1.695)--(0,0)--cycle;
\draw[draw=HKS36K80,draw opacity=0.7,line width=0.0001](0,0)--(0,0)--(.761,1.746)--(.87,1.695);
\filldraw[fill=HKS36K80,fill opacity=0.7,draw=none](0,0)--(1.128,1.287)--(1.076,1.291)--(0,0)--cycle;
\draw[draw=HKS36K80,draw opacity=0.7,line width=0.0001](1.128,1.287)--(1.076,1.291)--(0,0)--(0,0);
\draw[arrows=->,line width=.4pt](0,0)--(-1.68,-.123);
\filldraw[draw=HKS36K80,draw opacity=0.7,line width=0.0001,fill=HKS36K80,fill opacity=0.7](0,0)--(.438,1.765)--(.499,1.798)--(0,0)--cycle;
\filldraw[draw=HKS36K80,draw opacity=0.7,line width=0.0001,fill=HKS36K80,fill opacity=0.7](0,0)--(.439,1.699)--(.438,1.765)--(0,0)--cycle;
\filldraw[draw=HKS36K80,draw opacity=0.7,line width=0.0001,fill=HKS36K80,fill opacity=0.7](0,0)--(.499,1.798)--(.612,1.791)--(0,0)--cycle;
\filldraw[draw=HKS36K80,draw opacity=0.7,line width=0.0001,fill=HKS36K80,fill opacity=0.7](0,0)--(.502,1.608)--(.439,1.699)--(0,0)--cycle;
\filldraw[draw=HKS36K80,draw opacity=0.7,line width=0.0001,fill=HKS36K80,fill opacity=0.7](0,0)--(.612,1.791)--(.761,1.746)--(0,0)--cycle;
\filldraw[draw=HKS36K80,draw opacity=0.7,line width=0.0001,fill=HKS36K80,fill opacity=0.7](0,0)--(.617,1.507)--(.502,1.608)--(0,0)--cycle;
\filldraw[draw=HKS36K80,draw opacity=0.7,line width=0.0001,fill=HKS36K80,fill opacity=0.7](0,0)--(.766,1.412)--(.617,1.507)--(0,0)--cycle;
\filldraw[draw=HKS36K80,draw opacity=0.7,line width=0.0001,fill=HKS36K80,fill opacity=0.7](0,0)--(.927,1.335)--(.766,1.412)--(0,0)--cycle;
\filldraw[draw=HKS36K80,draw opacity=0.7,line width=0.0001,fill=HKS36K80,fill opacity=0.7](0,0)--(1.076,1.291)--(.927,1.335)--(0,0)--cycle;
\draw[arrows=<->,line width=.4pt](.336,-.613)--(0,0)--(0,1.733);
\draw[draw opacity=0.8,draw=HKS41K80](.763,-.631)--(0,0);
\draw[draw=red,line width=4](-1.8,-1)--(2,1.7);
\draw[draw=red,line width=4](-1.8,1.7)--(2,-1);
\path (.854,1.56) node {$P$};\path (.336,-.613) node[below] {$\vektor a_1$}
     (0,1.733) node[above] {$\vektor a_2$}
     (-1.68,-.123) node[left] {$\vektor a_3$};\end{tikzpicture}% End sketch output
}\hfill}
  \caption{Positive cone ($P$) and hyperplane $\hyperebene$
    according to the characterisation in the Theorems 
    \ref{satz:Hyperebenencharakterisierung} and
    \ref{satz:ordnungserweiterungscharakterisierung}}
  \label{fig:ordnungserweiterungscharakterisierung}
\end{figure}
Figure \ref{fig:ordnungserweiterungscharakterisierung} illustrates
of the following theorem.

\begin{lemma}\label{satz:Hyperebenencharakterisierung}
  Let $\paar {\vektorraum[R]}\leq$ be an ordered real vector space,
  which has the positive cone $\poskegel{\vektorraum[R]}$ with
  non-empty interior. A linearly ordered hyperplane \hyperebene\
  defines a linear order extension $\leq'$ of $\leq$ iff for its
  positive cone $\poskegel\hyperebene$ the following condition is
  true:
 \begin{equation}
   \poskegel{\vektorraum[R]}\cap\hyperebene \subseteq \poskegel 
   \hyperebene\cap \rand\poskegel{\vektorraum[R]}\label{eq:hyperebene-fuer-Ordnungserweiterung}
  \end{equation}
\end{lemma}

\begin{proof}
  Let $\leq'$ be a linear order on the vector space 
  $\vektorraum[R]$ defined by a hyperplane $\hyperebene$, which is an order extension of  the order
  $\leq$. Moreover, $\leq'$ is also a linear order on \hyperebene. Then for
  the positive cone $\poskegel{\vektorraum[R]}':=\Menge{\vektor x\in
    \vektorraum[R]}{\nullvektor \leq' \vektor x}$ the following
  relation is satisfied:
  \begin{equation}
    \label{eq:bw:ordnung-genau-dann-wenn-Hyperebene}
    \poskegel{\vektorraum[R]}\subseteq \poskegel{\vektorraum[R]}'
  \end{equation}
  Suppose there exists an element $\vektor x \in \inneres (
  \poskegel{\vektorraum[R]})\cap\hyperebene$. Hence, neighbourhood
  $U(\vektor x)$ exists such
  that $U(\vektor x)\subseteq \poskegel{\vektorraum[R]}$. On the
  other hand with $\hyperebene =\rand \poskegel{\vektorraum[R]}'$ the
  inequality $U(\vektor x)\setminus
  \poskegel{\vektorraum[R]}'\neq\emptyset$ holds.  So we deduce
  $\bigl(U(\vektor x)\cap\poskegel{\vektorraum[R]}\bigr)\setminus
  \poskegel{\vektorraum[R]}'\neq\emptyset$, which implies
  $\poskegel{\vektorraum[R]}\setminus
  \poskegel{\vektorraum[R]}'\neq\emptyset$ (\widerspruch{}). This is a
  contradiction to
  \eqref{eq:bw:ordnung-genau-dann-wenn-Hyperebene}. Consequently, the
  intersection
  $\inneres\poskegel{\vektorraum[R]}\cap\hyperebene=\emptyset$ is
  empty. Hence, $\poskegel{\vektorraum[R]}\cap\hyperebene\subseteq
  \rand \poskegel{\vektorraum[R]}$, which implies
  $\poskegel{\vektorraum[R]}\cap\hyperebene\subseteq \rand
  \poskegel{\vektorraum[R]}\cap\hyperebene$.

  Let now $\vektor
  x\in\poskegel{\vektorraum[R]}\cap\negkegel\hyperebene\setminus\menge{\nullvektor}$. In
  this case the vector $\vektor x$ would be positive with respect to
  the given order $\leq$ and with respect to the order $\leq'$ that is
  generated by the hyperplane \hyperebene\ it would be
  negative. Hence, the orders would not be
  compatible(\widerspruch). Since \hyperebene\ is linearly ordered,
  the remaining case is $\vektor x\in\poskegel\hyperebene$ and thus,
  \begin{equation*} 
    \poskegel{\vektorraum[R]}\cap\hyperebene \subseteq
    \poskegel\hyperebene. 
  \end{equation*}
  So far we have shown that
  each linear order extension of $\leq$ fulfils condition
  \eqref{eq:hyperebene-fuer-Ordnungserweiterung}.

  If the hyperplane \hyperebene\ fulfils the Condition
  \eqref{eq:hyperebene-fuer-Ordnungserweiterung}, on one side of the
  hyperplane there are no elements of $\poskegel{\vektorraum[R]}$.  In
  the contrary case two vectors $\vektor x
  \in\inneres\poskegel{\vektorraum[R]}$ and $\vektor y\in\hyperebene$
  and a positive real number
  $α\in\poskegel\reellezahlen\setminus\menge0$ exist such that
  $\vektor x$ is on one side of \hyperebene{} and $\vektor y' :=
  \vektor y -α(\vektor x-\vektor y)$
%
  \iffalse a neighbourhood $B_r(\vektor y-α\vektor
  x)\subseteq\poskegel{\vektorraum[R]}\setminus\hyperebene$ \fi
%
  is on the other side of \hyperebene. Then $\vektor y'-\vektor x =
  (1+α)(\vektor y - \vektor x)$. Let $U$ be an open environment of
  $\vektor x$.  Then the set
  \[
  U':=\Menge{\vektor z'\in\direktsumme E\reellezahlen}
  {\vektor z'=\vektor z + \frac{1}{1+α}(\vektor z - \vektor y'), \vektor z∈U}
  \]
  is neighbourhood of $\vektor y$ and as $\poskegel{\vektorraum[R]}$
  is convex $U'⊆\poskegel{\vektorraum[R]}$ is one of its subsets.
  
  Thus, we have
  $\vektor y\in \inneres\poskegel{\vektorraum[R]}\cap\hyperebene$
  (\widerspruch) in contradiction to 
  \eqref{eq:hyperebene-fuer-Ordnungserweiterung}. Thus only on one
  side of $\hyperebene$ can be positive elements with respect to the
  order $\leq$.

  If on both sides of $\hyperebene$ there are no elements of
  $\poskegel{\vektorraum[R]}$ i.e., 
  $\poskegel{\vektorraum[R]}\subseteq\hyperebene$, we can freely choose
  the side of the positive cone of $\leq'$. In this case we choose an
  arbitrary vector $\vektor
  a\in\vektorraum[R]\setminus\hyperebene$ and receive in dependence
  from $\vektor{a}$ one of the two possible orders.

  If the set  $\poskegel{\vektorraum[R]}\setminus\hyperebene$ is
  non-empty, we choose a vector 
  $\vektor a\in \poskegel{\vektorraum[R]}\setminus\hyperebene$. In
  this case 
 $P:=\poskegel\hyperebene\cup\bigl(\hyperebene+(\poskegel\reellezahlen\setminus\menge
  0)\vektor a\bigr)$ is the set of positive elements of a
  linearly ordered group. The independence from the chosen
  representative $\vektor a$ has been shown for 
 $\poskegel{\vektorraum[R]}\setminus\hyperebene\neq\emptyset$ in the
 proof of Theorem \ref{satz:poskegelrand=Hyperebene}.
\end{proof}

\noindent As for each vector $\vektor a∈\direktsumme{E}\reellezahlen$ the
subspace $\vektorraum'_{\vektor a}$ has a hyperplane that divides its
positive cone according to the order relation $≤'$ from the negative
one, we can immediately conclude the following corollary:

\begin{corollary}\label{folg:Hyperebenencharakterisierung}
  Let $U⊆\direktsumme{E}\rationalezahlen$ be a set that generates
  $\direktsumme E\reellezahlen$ in the following way:
  \begin{equation}
    \direktsumme{E}\reellezahlen = ⋃_{\vektor a∈U}\vektorraum_{\vektor a}.
  \end{equation}
  Then $\leq'$ is a linear extension of $\leq$ iff for each element
  $\vektor a∈U$ the following relation holds:
 \begin{equation}
   \poskegel{\vektorraum[R]}\cap\hyperebene'_{\vektor a} \subseteq 
   \poskegel{(\hyperebene'_{\vektor a})}\cap \rand\poskegel{\vektorraum[R]}
   \label{eq:hyperebene-fuer-Ordnungserweiterung2}
  \end{equation}
\end{corollary}

\begin{figure}\noindent\parbox[t]{0.4\linewidth}{\caption[Arrangement of dividing hyperplane and positive cone]
{\linebreak[2]\raggedright Arrangement of dividing hyperplane and positive cone}\label{fig:hyperebene-und-kegel}}%
\parbox[t]{0.6\linewidth}{\centerline{\raisebox{-\height}{% Sketch output, version 0.3 (build 7d, Fri Nov 22 14:14:09 2013)
% Output language: PGF/TikZ,LaTeX
\begin{tikzpicture}[line join=round]
\tikzstyle{ann} = [fill=white,font=\footnotesize,inner sep=1pt]
\providecommand\vektor{}\tikzstyle{plane}=[shading=axis,top color=white,bottom
color=HKS41K80,shading angle=-50]\tikzstyle{plane2}=[shading=axis,top color=white,bottom
color=HKS41K80,shading angle=-30]\draw[arrows=<-,line width=.4pt](.42,-.767)--(0,0);
\filldraw[plane,draw=black,draw opacity=0.01,line width=0.1,fill opacity=1](-1.68,-.833)--(2.52,.513)--(1.68,.833)--(-2.52,-.513)--cycle;
\draw[draw opacity=0.8,draw=white](-.42,.16)--(0,0);
\draw[arrows=->,line width=.4pt](0,0)--(0,2.167);
\draw[draw opacity=0.8,draw=HKS41K80](.42,-.16)--(0,0);
\draw[arrows=->,line width=.4pt](0,0)--(-2.1,-.153);
\filldraw[draw=HKS36K80,draw opacity=0.7,line width=0.0001,fill=HKS36K80,fill opacity=0.7](0,0)--(.548,2.207)--(.624,2.247)--(0,0)--cycle;
\filldraw[draw=HKS36K80,draw opacity=0.7,line width=0.0001,fill=HKS36K80,fill opacity=0.7](0,0)--(.624,2.247)--(.765,2.239)--(0,0)--cycle;
\filldraw[draw=HKS36K80,draw opacity=0.7,line width=0.0001,fill=HKS36K80,fill opacity=0.7](0,0)--(.765,2.239)--(.951,2.183)--(0,0)--cycle;
\filldraw[draw=HKS36K80,draw opacity=0.7,line width=0.0001,fill=HKS36K80,fill opacity=0.7](0,0)--(.951,2.183)--(1.152,2.087)--(0,0)--cycle;
\filldraw[draw=HKS36K80,draw opacity=0.7,line width=0.0001,fill=HKS36K80,fill opacity=0.7](0,0)--(1.152,2.087)--(1.339,1.968)--(0,0)--cycle;
\filldraw[draw=HKS36K80,draw opacity=0.7,line width=0.0001,fill=HKS36K80,fill opacity=0.7](0,0)--(1.339,1.968)--(1.482,1.842)--(0,0)--cycle;
\filldraw[draw=HKS36K80,draw opacity=0.7,line width=0.0001,fill=HKS36K80,fill opacity=0.7](0,0)--(1.482,1.842)--(1.561,1.728)--(0,0)--cycle;
\filldraw[draw=HKS36K80,draw opacity=0.7,line width=0.0001,fill=HKS36K80,fill opacity=0.7](0,0)--(1.561,1.728)--(1.562,1.645)--(0,0)--cycle;
\filldraw[draw=HKS36K80,draw opacity=0.7,line width=0.0001,fill=HKS36K80,fill opacity=0.7](0,0)--(1.562,1.645)--(1.486,1.605)--(0,0)--cycle;
\filldraw[draw=HKS36K80,draw opacity=0.7,line width=0.0001,fill=HKS36K80,fill opacity=0.7](0,0)--(1.486,1.605)--(1.344,1.613)--(0,0)--cycle;
\filldraw[draw=HKS36K80,draw opacity=0.7,line width=0.0001,fill=HKS36K80,fill opacity=0.7](0,0)--(1.344,1.613)--(1.159,1.669)--(0,0)--cycle;
\filldraw[draw=HKS36K80,draw opacity=0.7,line width=0.0001,fill=HKS36K80,fill opacity=0.7](0,0)--(1.159,1.669)--(.957,1.764)--(0,0)--cycle;
\filldraw[draw=HKS36K80,draw opacity=0.7,line width=0.0001,fill=HKS36K80,fill opacity=0.7](0,0)--(.957,1.764)--(.771,1.884)--(0,0)--cycle;
\filldraw[draw=HKS36K80,draw opacity=0.7,line width=0.0001,fill=HKS36K80,fill opacity=0.7](0,0)--(.771,1.884)--(.627,2.01)--(0,0)--cycle;
\filldraw[draw=HKS36K80,draw opacity=0.7,line width=0.0001,fill=HKS36K80,fill opacity=0.7](0,0)--(.627,2.01)--(.549,2.123)--(0,0)--cycle;
\filldraw[draw=HKS36K80,draw opacity=0.7,line width=0.0001,fill=HKS36K80,fill opacity=0.7](0,0)--(.549,2.123)--(.548,2.207)--(0,0)--cycle;
\path (1.068,1.95) node {$P$};\path (.42,-.767) node[below] {$\vektor a_1$}
     (0,2.167) node[above] {$\vektor a_2$}
     (-2.1,-.153) node[left] {$\vektor a_3$};\path (1.159,.353) node {$\poskegel\hyperebene$};\path[white] (-1.159,-.353) node {$\negkegel\hyperebene$};\end{tikzpicture}% End sketch output
}}}%
\end{figure}
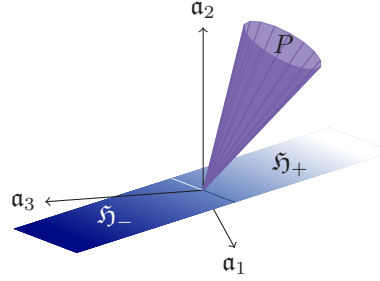%

\noindent The main theorem of this article is illustrated in Figure
\ref{fig:hyperebene-und-kegel}. It describes the set of linear order
extensions of a given order on an abelian partially ordered group:

\begin{theorem}[Characterisation of linear order extensions on
  abelian po-groups]\label{satz:ordnungserweiterungscharakterisierung}
  Let $\gruppe=\OrdGruppe G+-0\leq$ a po-group, $E$ a maximal
  independent set in $\gruppe$, $\direktsumme E\reellezahlen$ the
  corresponding vector space with the embeddings $\varphi:G\to
  \direktsumme E\reellezahlen$ and $\psi:\direktsumme E\ganzzahlen\to
  G$ as described in \eqref{eq:homomorphismus-G-Q} and
  \eqref{eq:z-g-mornomorphismus}. Then each linear order extension of
  $\leq$ can be characterised by a set
  $U⊆\direktsumme{E}\rationalezahlen$ and for each element $\vektor
  a∈U$ a linearly ordered hyperplane $\hyperebene_{\vektor a}$ which fulfils the
  condition
 \begin{equation}\label{eq:ordnungserweiterungscharakterisierung}
   \conv(\varphi[\poskegel{\gruppe}])\cap\hyperebene_{\vektor a}
   \subseteq \poskegel{(\hyperebene_{\vektor a})}\cap 
   \rand\Bigl(\conv\bigl(\varphi[\poskegel{\gruppe}]\bigr)\Bigr).
  \end{equation}
  For each Element $\vektor a∈ U$ the corresponding hyperplane
  $\hyperebene_{\vektor a}$ defines a unique linear order extension in
  the subspace $\vektorraum_{\vektor a}$ if $\vektorraum_{\vektor
    a}∩\varphi[\poskegel{\gruppe}]\setminus\hyperebene_{\vektor a}$ is
  non-empty, otherwise it can be used to define exactly 2 linear order
  extensions in this subspace.
\end{theorem}

\begin{proof}
  The existence of the monomorphisms has bee proved in the Lemmata
  \ref{satz:Monomorphismus-G-Q} and
  \ref{lemma:Monomorphismus-Z-G}. With Theorem
  \ref{lemma:Z-vektorraum} and
  Lemma~\ref{lemma:ordungsuebertragung-G-Q} we get the description of
  the positive cone $\poskegel{\gruppe}$ as
  $\conv(\varphi[\poskegel{\gruppe}])$ in the vector space
  $\direktsumme{E}\reellezahlen$. For non-semiclosed orders this
  construction adds only these elements to the positive cone of
  $\gruppe$, which are positive in every linear order extension on
  \gruppe{}, as Lemma~\ref{lemma:halbabschluss} shows. Thus this
  restriction can be loosened. Theorem
  \ref{satz:charakterisierung-linearer-ordnungen} gives us a
  description of linear orders on $\gruppe$ in the vector space
  $\direktsumme{E}\reellezahlen$.
  Corollary~\ref{folg:Hyperebenencharakterisierung} provides the
  Condition \eqref{eq:ordnungserweiterungscharakterisierung} for the
  characterisations of the linear orders on $\direktsumme
  E\reellezahlen$ as well as the number of orders defined by the
  hyperplane. Thus all parts of this theorem have been proved
  elsewhere in this article.
 \end{proof}

\noindent With Theorem
\ref{satz:ordnungserweiterungscharakterisierung} we have described a
geometrical characterisation of linear order extensions of partially
ordered torsion free abelian groups. This allows us, to use methods
from the convex geometry for the construction and investigation of such
linear order extensions. The following section illustrates this.

\section{Further results}\label{sec:further-results}
As an illustration of the method provided above, in this section some
other proofs and additional conclusions are added.

So far we have only described existing order extensions. The question
has not been touched whether such an order extension or hyperplane
does exist. The following theorem fills this gap. It has been proofed
several times
(cf. \cite{Everett:1950,Lorenzen:1939,fuchs:1963,Simbireva:1947}). We
add an additional proof:

\begin{theorem}\label{satz:existenz-lineare-ordnungserweiterung}
  Each partially ordered torsion free group \gruppe{} has a compatible linear
  order extension.
\end{theorem}
\begin{proof}
  Consider the complete sublattice of the subspace lattice of
  $\vektorraum$ that is generated by the set
  $W:=\Menge{\vektorraum_{\vektor a}}{\vektor a∈\vektorraum}$. This
  sublattice has a maximal chain $\mathfrak K$. As each vector
  $\vektor a∈\vektorraum$ is contained in an element of $\mathfrak K$,
  the union $⋃\mathfrak K=\vektorraum$ is the complete vector space
  $\vektorraum$. We define
  \[
  \vektorraum[S]_{\vektor a}:=
  ⋃\Menge{\vektorraum[W]∈\mathfrak K}{\vektor a\not ∈\vektorraum[W]}.
  \]
  As the subspace
  $\vektorraum[W]_a:=\vektorraum_a+\vektorraum[S]_{\vektor a}$ is the
  join of $\vektorraum_{\vektor a}$ and $\vektorraum[S]_{\vektor a}$
  in the subspace lattice, $\vektorraum[W]_{\vektor a}$ is the upper
  neighbour of $\vektorraum[S]_{\vektor a}$ in $\mathfrak K$. Then
  there exists a maximal linear independent set
  $B⊆\vektorraum[W]_{\vektor a}$ that contains $\vektor a$ and is
  liner independent from $\vektorraum[S]_{\vektor a}$. Then for each
  vector $\vektor b∈B$ the vector space $\vektorraum_{\vektor b}$ is
  the same as $\vektorraum_{\vektor a}$ (otherwise $\vektor
  b∈\vektorraum[S]_{\vektor a}$).

  For each vector $\vektor c∈\vektorraum_{\vektor
    a}∩\poskegel\vektorraum$ there exists a real number
  $c∈\reellezahlen$ such that $c\vektor c ≤ \vektor a$ which implies
  $\nullvektor ≤\vektor a-c\vektor c$. Thus there exists a basis
  $\basis'_{\vektor a}$ containing $\vektor a$ and a basis of
  $\vektorraum[S]_{\vektor a}$ such that for all Elements $\vektor
  x∈\negkegel{\vektorraum}∩\vektorraum_{\vektor a}$ the distance
  $\|\vektor a - \vektor x\|_{\infty,\basis}≥
  d∈\poskegel\reellezahlen\setminus\menge 0$ is strictly positive
  (choose all base vectors to be positive between $\vektor a$ and
  $\nullvektor$) and such that $\vektor a$ is the largest element in
  the basis $\basis'_{\vektor a}$. Then
  $\conv\bigl(B_{\frac{d}2,\infty,\basis'}(\vektor
  a)∪\vektorraum[W]_{\vektor a}\bigr)$ has non-empty interior and does not hit $\negkegel\vektorraum$.
%
  \iffalse
  With respect to the given norm we
  can find an open ball $B_r(\vektor
  a)⊆\poskegel\vektorraum∩\vektorraum[W]_{\vektor a}$ around $\vektor
  a$. As by definition for every element $\vektor
  x∈\poskegel\vektorraum∩\vektorraum[W]_{\vektor a}\setminus
  \vektorraum[S]_{\vektor a}$ the generated vector space
  $\vektorraum[W]_{\vektor x}=\vektorraum[W]_{\vektor a}$ is the whole
  containing space, we can find such an environment around every
  non-zero element of the positive cone outside of
  $\vektorraum[S]_{\vektor a}$. Thus $\vektorraum[W]_{\vektor
    a}∩\poskegel\vektorraum\setminus\vektorraum[S]_{\vektor a}$ is an
  open set with respect to the norm $\|\cdot\|_{1,\basis}$.
  \fi
  % 
  Then there exists a hyperplane $\hyperebene_{\vektor a}$ such that
  $\vektorraum[S]_{\vektor a}⊆\hyperebene_{\vektor a}$ and
  $\hyperebene_{\vektor a}$ does not hit $\vektorraum[W]_{\vektor
    a}∩\poskegel\vektorraum\setminus\vektorraum[S]_{\vektor a}$.
  
  Let us chose a base $\hat\basis_{\vektor a}$ such that it contains a base of
  $\vektorraum[S]_{\vektor a}$ together with $\vektor a$ with the
  additional condition that
  $\hat\basis_{\vektor a}\setminus(\menge{\vektor a}∪\vektorraum[S]_{\vektor
    a})⊆\hyperebene_{\vektor a}$. Then we define a linear order
  relation $≤'$ on $\hat\basis_{\vektor a}$ that has $\vektor a$ as largest
  element and is defined in a way such that for any two elements
  $\vektor x∈\hat\basis_{\vektor a}∩\hyperebene\setminus\vektorraum[S]_{\vektor a}$
  and $\vektor y∈\vektorraum[S]_{\vektor a}$ the inequality $\vektor
  y≤'\vektor x$ holds. Then $≤'$ defines a linear order on
  $\vektorraum[W]_{\vektor a}$ such that the positive cone of
  $\vektorraum[W]_{\vektor a}$ with respect to the original order $≤$
  is a subset of the positive cone that is defined by
  $≤'$. Furthermore we can use
  $\hat\basis_{\vektor a}\setminus\vektorraum[S]_{\vektor a}$ to extend any
  linearly ordered base of $\vektorraum[S]_{\vektor a}$ to a linearly
  ordered base of $\vektorraum[W]_{\vektor a}$ that fulfils the same
  condition. 

  As the vector $\vektor a$ has been chosen arbitrarily, we can find
  such a subbase $\hat\basis_{\vektor a}$ for every subspace
  $\vektorraum[W]_{\vektor a}∈\mathfrak K$, and the union $\hat\basis$
  of all these subbases is a basis of the original vector space
  $\vektorraum$. The order on the chain $\mathfrak K$ defines a linear
  preorder on $\hat\basis$, which can be turned into a linear order by
  using the orders of the subbases $\hat\basis_{\vektor a}$ as
  developed above. Thus we can define a linear order whose positive
  cone is a superset of the one of the original order $≤$, which means
  that it is a linear order extension.

  As every torsion free ordered abelian group $\gruppe$ has a closed
  order extension, we can use this construction to construct a linear
  order extension on $\gruppe$.
\end{proof}

\noindent Looking at the proof and the preceding sections it is
obvious that every linear order extension can be found with this
construction: Every linear order can be expressed by a maximal chain
in the subspace lattice, and it is an order extension iff it provides
an extension of the positive cone.

\noindent In general, an order relation on a given set is the
intersection of its linear order extensions
(cf. \cite{Szpilrajn:1930}). If the order must be compatible to a given
algebraic structure, the set of choices for the linear extensions is
much more restricted. For abelian groups the following theorem has
been proved a long time ago (cf.\ \cite{Fuchs:1958}, Corollary
6). Also in this case our method provides an alternative proof.

\begin{theorem}
  Let $\gruppe=\OrdGruppe G+-0\leq$ be a partially ordered torsion free
  abelian group. Then the order relation $\leq$ is the intersection of
  all of its linear order extensions iff it is semiclosed.
\end{theorem}

\begin{proof}
  One direction is obvious: Each lattice ordered group is semiclosed,
  thus also each linearly ordered group. So, the positive cone of an
  arbitrary order on a given group \gruppe{} contains for each 
  $n\in\natzahlen\setminus\menge 0$ together with $na$  also the element
  $a$ and, thus, the intersection is also semiclosed 
  (cf. Lemma \ref{lemma:halbabschluss}).

  Considering the other direction, let $\leq$ be
  semiclosed. Furthermore, let $\mathfrak O$ the set of all linear
  order extensions of $\leq$ which are compatible with \gruppe\ and $E$
  a maximal independent set in \gruppe. Then for each  monomorphism
$\varphi:\gruppe\to \direktsumme
  E\reellezahlen$ with $a\mapsto \vektor e_a$ for all  $a\in E$ the
  set $P:=\conv\bigl[\varphi[\poskegel{\gruppe}]\bigr]$ is the unique
  representation of the positive cone in $\direktsumme E\reellezahlen$
  according to Theorem \ref{lemma:Z-vektorraum} and Lemma
  \ref{lemma:ordungsuebertragung-G-Q}. Thus it is sufficient to show
  that $P$ is the intersection of the corresponding positive half
  spaces of elements of $\mathfrak O$. Let $P'$ the intersection of all of these half spaces of the
  linear order extensions. Then 
  $P'$ is non-empty because of Theorem
  \ref{satz:existenz-lineare-ordnungserweiterung}. Furthermore,
  $P\subseteq P'$. Now, we will show the inclusion 
  $\direktsumme E\reellezahlen\setminus P\subseteq \direktsumme
  E\reellezahlen\setminus P'$.

Let $\vektor x\in\direktsumme E\reellezahlen\setminus
  P$. We infer from the convexity of $P$ the equality
  $\reellezahlen\vektor x\cap P=\menge\nullvektor$. Let $\vektor
  y\in \conv(\poskegel\reellezahlen \vektor x\cup
  -P)\cap \conv(\negkegel\reellezahlen \vektor x\cup
  P)=\linhuelle^+(\poskegel\reellezahlen \vektor x\cup -P)\cap
  \linhuelle^+(\negkegel\reellezahlen \vektor x\cup P)$. The two sets
  are equal as both are convex and $\nullvektor \in
  \poskegel\reellezahlen\vektor x\cap P$. Then there exist
  $\alpha,β\in\poskegel\reellezahlen$ and vectors  $\vektor p,\vektor q\in P$
 such that $α\vektor x-\vektor p=\vektor y=-β\vektor x+\vektor
 q$. Consequently,
  $(α+β)\vektor x=\vektor p+\vektor q\in P$, which implies with the
  properties of being a cone for $α\neq 0$ or
  $β\neq 0$ that also $\vektor x$ is in $P$
  (\widerspruch). Consequently, $α=β=0$ and it follows $\vektor y=\vektor
  p=\vektor q=\nullvektor$.
  This implies $\conv (\negkegel \reellezahlen {\vektor
    x}\cup P)\cap \conv(\poskegel\reellezahlen{\vektor x}\cup
  -P)=\menge \nullvektor$.

  As proven in Theorem
  \ref{satz:existenz-lineare-ordnungserweiterung}, a linear order
  extension of $\leq$ exists, for which $\vektor
  x$ is negative. Consequently, $\vektor x \not\in P'$. Thus, $\direktsumme E\reellezahlen\setminus
  P\subseteq \direktsumme E\reellezahlen\setminus
  P'$  follows and so $P'\subseteq P$, from which we get  $P=P'$. This means
  that a semiclosed order $\leq$ is the intersection of all of its
  order extensions.
\end{proof}

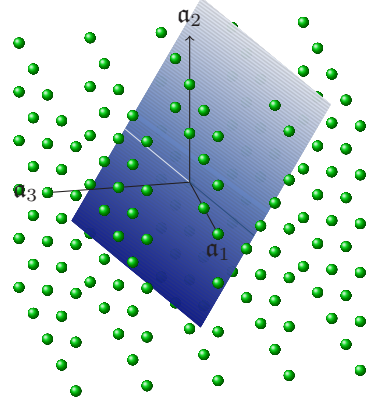
\begin{figure}\noindent\parbox[t]{0.5\linewidth}{\caption[Hyperplane and the image of  a partially ordered group]
{\linebreak[2]\raggedright Hyperplane and the image of  a partially
  ordered group. The hyperplane is arranged in a way that it does not intersect
with the images of the group elements.}\label{fig:ebene-im-raster}}%
\parbox[t]{0.5\linewidth}{\centerline{\raisebox{-\height}{% Sketch output, version 0.3 (build 7d, Fri Nov 22 14:14:09 2013)
% Output language: PGF/TikZ,LaTeX
\begin{tikzpicture}[line join=round]
\tikzstyle{ann} = [fill=white,font=\footnotesize,inner sep=1pt]
\providecommand\vektor{}\tikzstyle{plane}=[shading=axis,top color=white,bottom
color=HKS41K80,shading angle=-50]\tikzstyle{plane2}=[shading=axis,top color=white,bottom
color=HKS41K80,shading angle=-30]\tikzstyle{rasterpunkt}=[shading=ball,ball color=HKS57K100,draw=none]\filldraw[rasterpunkt](6.512,-1.122) circle (2pt);
\filldraw[rasterpunkt](5.567,-1.191) circle (2pt);
\filldraw[rasterpunkt](6.512,-.472) circle (2pt);
\filldraw[rasterpunkt](4.622,-1.26) circle (2pt);
\filldraw[rasterpunkt](5.567,-.541) circle (2pt);
\filldraw[rasterpunkt](6.512,.178) circle (2pt);
\filldraw[rasterpunkt](3.677,-1.329) circle (2pt);
\filldraw[rasterpunkt](4.622,-.61) circle (2pt);
\filldraw[rasterpunkt](5.567,.109) circle (2pt);
\filldraw[rasterpunkt](2.732,-1.398) circle (2pt);
\filldraw[rasterpunkt](6.512,.828) circle (2pt);
\filldraw[rasterpunkt](3.677,-.679) circle (2pt);
\filldraw[rasterpunkt](4.622,.04) circle (2pt);
\filldraw[rasterpunkt](5.567,.759) circle (2pt);
\filldraw[rasterpunkt](4.622,.69) circle (2pt);
\filldraw[rasterpunkt](5.567,1.409) circle (2pt);
\filldraw[rasterpunkt](4.811,-1.605) circle (2pt);
\filldraw[rasterpunkt](4.622,1.34) circle (2pt);
\filldraw[rasterpunkt](5.567,2.059) circle (2pt);
\filldraw[rasterpunkt](4.811,-.955) circle (2pt);
\filldraw[rasterpunkt](5.756,-.236) circle (2pt);
\filldraw[rasterpunkt](4.811,-.305) circle (2pt);
\filldraw[rasterpunkt](5.756,.414) circle (2pt);
\filldraw[rasterpunkt](6.701,1.133) circle (2pt);
\filldraw[rasterpunkt](4.811,.345) circle (2pt);
\filldraw[rasterpunkt](5.756,1.064) circle (2pt);
\filldraw[rasterpunkt](5.756,1.714) circle (2pt);
\filldraw[rasterpunkt](5,-1.3) circle (2pt);
\filldraw[rasterpunkt](5,-.65) circle (2pt);
\filldraw[rasterpunkt](5.945,.069) circle (2pt);
\filldraw[rasterpunkt](5,0) circle (2pt);
\filldraw[rasterpunkt](5.945,.719) circle (2pt);
\filldraw[rasterpunkt](5.945,1.369) circle (2pt);
\filldraw[rasterpunkt](5.189,-1.645) circle (2pt);
\filldraw[rasterpunkt](5.189,-.995) circle (2pt);
\filldraw[rasterpunkt](6.134,.374) circle (2pt);
\filldraw[plane2,draw=black,draw opacity=0.01,line width=0.1,fill opacity=0.8](5.146,-1.931)--(6.876,1.036)--(5.159,2.455)--(3.429,-.512)--cycle;
\filldraw[rasterpunkt](6.701,-1.467) circle (2pt);
\filldraw[rasterpunkt](2.732,-.748) circle (2pt);
\filldraw[rasterpunkt](6.512,1.478) circle (2pt);
\filldraw[rasterpunkt](3.677,-.029) circle (2pt);
\filldraw[rasterpunkt](5.756,-1.536) circle (2pt);
\draw[draw opacity=0.8,draw=white](4.141,.709)--(5,0);
\filldraw[rasterpunkt](6.701,-.817) circle (2pt);
\filldraw[rasterpunkt](2.732,-.098) circle (2pt);
\filldraw[rasterpunkt](6.512,2.128) circle (2pt);
\filldraw[rasterpunkt](3.677,.621) circle (2pt);
\filldraw[rasterpunkt](5.756,-.886) circle (2pt);
\filldraw[rasterpunkt](6.701,-.167) circle (2pt);
\filldraw[rasterpunkt](3.866,-1.674) circle (2pt);
\filldraw[rasterpunkt](2.732,.552) circle (2pt);
\filldraw[rasterpunkt](3.677,1.271) circle (2pt);
\filldraw[rasterpunkt](2.921,-1.743) circle (2pt);
\filldraw[rasterpunkt](4.622,1.99) circle (2pt);
\filldraw[rasterpunkt](6.701,.483) circle (2pt);
\filldraw[rasterpunkt](3.866,-1.024) circle (2pt);
\filldraw[rasterpunkt](2.732,1.202) circle (2pt);
\filldraw[rasterpunkt](6.89,-1.812) circle (2pt);
\filldraw[rasterpunkt](3.677,1.921) circle (2pt);
\filldraw[rasterpunkt](2.921,-1.093) circle (2pt);
\filldraw[rasterpunkt](3.866,-.374) circle (2pt);
\filldraw[rasterpunkt](5.945,-1.881) circle (2pt);
\filldraw[rasterpunkt](2.732,1.852) circle (2pt);
\filldraw[rasterpunkt](6.89,-1.162) circle (2pt);
\filldraw[rasterpunkt](2.921,-.443) circle (2pt);
\filldraw[rasterpunkt](5,-1.95) circle (2pt);
\filldraw[rasterpunkt](6.701,1.783) circle (2pt);
\filldraw[rasterpunkt](3.866,.276) circle (2pt);
\filldraw[rasterpunkt](5.945,-1.231) circle (2pt);
\filldraw[rasterpunkt](4.811,.995) circle (2pt);
\filldraw[rasterpunkt](6.89,-.512) circle (2pt);
\filldraw[rasterpunkt](4.055,-2.019) circle (2pt);
\filldraw[rasterpunkt](2.921,.207) circle (2pt);
\filldraw[rasterpunkt](3.866,.926) circle (2pt);
\filldraw[rasterpunkt](5.945,-.581) circle (2pt);
\filldraw[rasterpunkt](3.11,-2.088) circle (2pt);
\filldraw[rasterpunkt](4.811,1.645) circle (2pt);
\filldraw[rasterpunkt](6.89,.138) circle (2pt);
\filldraw[rasterpunkt](4.055,-1.369) circle (2pt);
\filldraw[rasterpunkt](2.921,.857) circle (2pt);
\filldraw[rasterpunkt](7.079,-2.157) circle (2pt);
\filldraw[rasterpunkt](3.866,1.576) circle (2pt);
\filldraw[rasterpunkt](3.11,-1.438) circle (2pt);
\filldraw[rasterpunkt](6.89,.788) circle (2pt);
\filldraw[rasterpunkt](4.055,-.719) circle (2pt);
\filldraw[rasterpunkt](6.134,-2.226) circle (2pt);
\filldraw[rasterpunkt](2.921,1.507) circle (2pt);
\draw[arrows=->,line width=.4pt](5,0)--(3.11,-.138);
\draw[arrows=<->,line width=.4pt](5.378,-.69)--(5,0)--(5,1.95);
\draw[draw opacity=0.8,draw=HKS41K80](5.858,-.71)--(5,0);
\filldraw[rasterpunkt](7.079,-1.507) circle (2pt);
\filldraw[rasterpunkt](3.11,-.788) circle (2pt);
\filldraw[rasterpunkt](5.189,-2.295) circle (2pt);
\filldraw[rasterpunkt](6.89,1.438) circle (2pt);
\filldraw[rasterpunkt](4.055,-.069) circle (2pt);
\filldraw[rasterpunkt](6.134,-1.576) circle (2pt);
\filldraw[rasterpunkt](5,.65) circle (2pt);
\filldraw[rasterpunkt](7.079,-.857) circle (2pt);
\filldraw[rasterpunkt](4.244,-2.364) circle (2pt);
\filldraw[rasterpunkt](3.11,-.138) circle (2pt);
\filldraw[rasterpunkt](4.055,.581) circle (2pt);
\filldraw[rasterpunkt](6.134,-.926) circle (2pt);
\filldraw[rasterpunkt](3.299,-2.433) circle (2pt);
\filldraw[rasterpunkt](5,1.3) circle (2pt);
\filldraw[rasterpunkt](7.079,-.207) circle (2pt);
\filldraw[rasterpunkt](4.244,-1.714) circle (2pt);
\filldraw[rasterpunkt](3.11,.512) circle (2pt);
\filldraw[rasterpunkt](7.268,-2.502) circle (2pt);
\filldraw[rasterpunkt](4.055,1.231) circle (2pt);
\filldraw[rasterpunkt](6.134,-.276) circle (2pt);
\filldraw[rasterpunkt](3.299,-1.783) circle (2pt);
\filldraw[rasterpunkt](7.079,.443) circle (2pt);
\filldraw[rasterpunkt](4.244,-1.064) circle (2pt);
\filldraw[rasterpunkt](6.323,-2.571) circle (2pt);
\filldraw[rasterpunkt](3.11,1.162) circle (2pt);
\filldraw[rasterpunkt](5.189,-.345) circle (2pt);
\filldraw[rasterpunkt](7.268,-1.852) circle (2pt);
\filldraw[rasterpunkt](3.299,-1.133) circle (2pt);
\filldraw[rasterpunkt](5.378,-2.64) circle (2pt);
\filldraw[rasterpunkt](7.079,1.093) circle (2pt);
\filldraw[rasterpunkt](4.244,-.414) circle (2pt);
\filldraw[rasterpunkt](6.323,-1.921) circle (2pt);
\filldraw[rasterpunkt](5.189,.305) circle (2pt);
\filldraw[rasterpunkt](7.268,-1.202) circle (2pt);
\filldraw[rasterpunkt](4.433,-2.709) circle (2pt);
\filldraw[rasterpunkt](6.134,1.024) circle (2pt);
\filldraw[rasterpunkt](3.299,-.483) circle (2pt);
\filldraw[rasterpunkt](5.378,-1.99) circle (2pt);
\filldraw[rasterpunkt](4.244,.236) circle (2pt);
\filldraw[rasterpunkt](6.323,-1.271) circle (2pt);
\filldraw[rasterpunkt](3.488,-2.778) circle (2pt);
\filldraw[rasterpunkt](5.189,.955) circle (2pt);
\filldraw[rasterpunkt](7.268,-.552) circle (2pt);
\filldraw[rasterpunkt](4.433,-2.059) circle (2pt);
\filldraw[rasterpunkt](3.299,.167) circle (2pt);
\filldraw[rasterpunkt](5.378,-1.34) circle (2pt);
\filldraw[rasterpunkt](4.244,.886) circle (2pt);
\filldraw[rasterpunkt](6.323,-.621) circle (2pt);
\filldraw[rasterpunkt](3.488,-2.128) circle (2pt);
\filldraw[rasterpunkt](7.268,.098) circle (2pt);
\filldraw[rasterpunkt](4.433,-1.409) circle (2pt);
\filldraw[rasterpunkt](3.299,.817) circle (2pt);
\filldraw[rasterpunkt](5.378,-.69) circle (2pt);
\filldraw[rasterpunkt](6.323,.029) circle (2pt);
\filldraw[rasterpunkt](3.488,-1.478) circle (2pt);
\filldraw[rasterpunkt](7.268,.748) circle (2pt);
\filldraw[rasterpunkt](4.433,-.759) circle (2pt);
\filldraw[rasterpunkt](5.378,-.04) circle (2pt);
\filldraw[rasterpunkt](6.323,.679) circle (2pt);
\filldraw[rasterpunkt](3.488,-.828) circle (2pt);
\filldraw[rasterpunkt](4.433,-.109) circle (2pt);
\filldraw[rasterpunkt](5.378,.61) circle (2pt);
\filldraw[rasterpunkt](3.488,-.178) circle (2pt);
\filldraw[rasterpunkt](4.433,.541) circle (2pt);
\filldraw[rasterpunkt](3.488,.472) circle (2pt);
\path (5.378,-.69) node[below] {$\vektor a_1$}
     (5,1.95) node[above] {$\vektor a_2$}
     (3.11,-.138) node[left] {$\vektor a_3$};\end{tikzpicture}% End sketch output
}}}%
\end{figure}%

Hölder's theorem\footnote{In fact, Hölder
  proved the theorem not for archimedian orders, but for such ordered semigroups that
  provide a Dedekind cut. For a link to the current theorem
  cf. \cite{Birkhoff93}, XIII.12.} has not been used so far, but can be proved easily,
now. The main idea will be formulated as a separate lemma:

\begin{lemma}
  A linearly ordered base of $\direktsumme{E}\reellezahlen$ defines an
  archimedian order on $\direktsumme{E}\rationalezahlen$ iff it
  fulfils the following conditions:
  \begin{enumerate}
  \item There exists a largest basis vector $\vektor b$.
  \item The largest vector $\vektor b$ is the only self-activator in $\basis$.
  \end{enumerate}
\end{lemma}
\begin{proof}
  At first we show that an archimedian order has only one
  self-activator. Suppose there are two self-activators $\vektor a$
  and $\vektor b$. W.l.o.g. we can assume $\vektor a,\vektor
  b∈\direktsumme{E}\rationalezahlen$. As the basis is linearly ordered
  either $\vektor a<\vektor b$ or $\vektor b <\vektor a$ is
  true. Let's assume the first one holds. This implies $\vektor
  a∈\hyperebene_{\vektor b}$, which means $\vektor a$ is infinitesimal
  smaller than $\vektor b$ (\widerspruch). Thus there is at most one
  self-activator in $\basis$. As
  $\direktsumme{E}\rationalezahlen⊆\linhuelle\basis$ holds there is at
  least one self-activator in $\basis$.

  The projection along $\basis$ into the subspace that is generated by
  two basis vectors is 2-dimensional as proofed in
  Lemma~\ref{lemma:projektion-erhaelt-dimension}. This implies that
  every non-self-activator has a self-activator above it. Thus the
  unique self-activator of an archimedian linear order must be the
  maximal element.

  On the other hand if these two conditions are fulfilled, we can
  choose the maximal element $\vektor b$. Then the corresponding
  hyperplane $\hyperebene_{\vektor b}$ is disjoint from
  $\direktsumme{E}\rationalezahlen$. Thus the order is completely
  defined by the projection along \basis\ into $\reellezahlen\vektor
  b$. Suppose there there are two vectors $\vektor a_1$ and $\vektor
  a_2$ in $\direktsumme{E}\rationalezahlen$ with the same projection
  into $\reellezahlen\vektor b$. Then the difference $\vektor
  a_1-\vektor a_2∈\hyperebene_{\vektor b}∩
  \direktsumme{E}\rationalezahlen$ is both infinitesimal smaller than
  $\vektor b$ and rational (\widerspruch).
\end{proof}

\noindent Geometrically, this implies that every archimedian linear
order can be represented by a linearly ordered hyperplane in the
complete vector space $\direktsumme{E}\reellezahlen$ that does not
intersect with the rational vectors from
$\direktsumme{E}\rationalezahlen$ as shown in
Figure~\ref{fig:ebene-im-raster}.

It has been proved many times that an archimedian $o$-group is
abelian. This shall not be repeated, here. Thus, we prove the
following theorem only for abelian groups.

\begin{theorem}[\cite{Birkhoff93}, XIII.12; Hölder, \cite{Holder:1901}; ]
  Every archimedian linearly ordered abelian group is isomorphic to a subgroup
  of the real numbers.
\end{theorem}

\begin{proof}
  Let $E$ be a maximal independent set in the group and $\basis$ a
  linearly ordered base of $\direktsumme{E}\reellezahlen$. Then, by
  the preceding lemma the order on the group is defined by the
  projection along $\basis$ into the subspace $\reellezahlen\vektor b$
  of the maximal basis vector.
\end{proof}

\begin{corollary}
  The abelian group $\ganzzahlen^\reellezahlen$  with respect to
  elementwise addition and subtraction does not permit an archimedian
  linear order.
\end{corollary}

\begin{proof}
  For the cardinalities we get the inequality
  \[
  |\ganzzahlen^\reellezahlen|≥2^{|\reellezahlen|}=|\potmenge{\reellezahlen}|>|\reellezahlen|.
  \]
  Thus, there is no isomorphism from $\ganzzahlen^\reellezahlen$ into a
  subset of the real numbers.
\end{proof}

\noindent The following theorem touches a question that has been
published in \cite{Bludov:2009} as Problem~1.7. As it is weaker, it
can be considered only as a step into the right direction.

\begin{theorem}
  Every torsion free abelian group permits an archimedian directed order.
\end{theorem}
\begin{proof}
  Let $\gruppe$ be an abelian group. And $E$ a maximal independent
  set. Then \gruppe\ can be embedded into
  $\direktsumme{E}\rationalezahlen$. Let
  $\abbildung{φ}{\gruppe}{\direktsumme{E}\rationalezahlen}$ be the
  corresponding embedding. The product order ($\vektor a ≤\vektor b ⇔
  ∀e∈E:\vektor a(e)≤\vektor b(e)$) on
  $\direktsumme{E}\rationalezahlen$ is an archimedian order: For any
  two vectors $\vektor a$ and $\vektor b$ the set $\supp{\vektor
    a}∪\supp{\vektor b}$ is finite. If both $\supp{\vektor a}\setminus
  \supp{\vektor b}≠∅$ and $\supp{\vektor b}\setminus \supp{\vektor
    a}≠∅$ are non-empty either $\vektor a\parallel\vektor b$ or
  $\vektor a\parallel -\vektor b$ holds, thus the two vectors fulfil
  the archimedian property. 

  For the remaining case we may assume $\supp{\vektor b}⊆
  \supp{\vektor a}$. Then there exist an element $e∈E$ and integers
  $a,b∈\ganzzahlen$ such that $a\vektor a(e)=b\vektor b(e)≠0$. If
  $\vektor a(a)$ is strictly positive, then $a$ can be chosen strictly
  positive, too.  This leads to $(a+1)\vektor a(a)>b\vektor
  b(e)>\vektor b(e)$ and $(a+1)\vektor a\not<\vektor b$ on the one
  hand, and $2b\vektor b(e) = 2a\vektor a(e) > \vektor a(e)$ and
  $2b\vektor b\not < \vektor a$ on the other hand, the archimedian
  property. For a negative coordinate $\vektor a(e)$ and a negative
  scalar $a$ the inequalities read $(a-1)\vektor a(e)>b\vektor
  b(e)>\vektor b(e)$, $(a-1)\vektor a\not<\vektor b$, $2b\vektor b(e)
  = 2a\vektor a(e) > \vektor a(e)$ and $2b\vektor b\not < \vektor
  a$. As $φ$ is an embedding, this proves also that the given order is
  also archimedian on $\gruppe$.

  Finally the order is a directed order. Let $\vektor a$ and $\vektor
  b$ the vectors that correspond to two arbitrary group elements. Then
  we define a mapping $\abbildung{\vektor c}{E}{\ganzzahlen}$ such
  that $\vektor c(e)$ is the least integer that fulfils the conditions
  \[
  \vektor a(e)≤\vektor c(e)\qquad\text{and}\qquad\vektor b(e)≤\vektor c(e).
  \]
  With the relation $\supp {\vektor c} ⊆ \supp{\vektor
    a}∪\supp{\vektor b}$ the vector $\vektor c$ is a well-defined
  element of $\direktsumme{E}\ganzzahlen$ and thus it represents an
  element of $\gruppe$ such that $\vektor a ≤ \vektor c$ and $\vektor
  b≤\vektor c$.
\end{proof}

\section{Conclusion and further topics}

In Theorem \ref{satz:charakterisierung-linearer-ordnungen} we have
provided a method to describe linearly ordered abelian groups by means
of linearly ordered hyperplanes in a the vector space, arising from
the set of mappings with finite support from a maximal independent set
into the real numbers. The vector space has been defined as the direct
sum $\direktsumme{E}\reellezahlen$, where $E$ denotes a maximal
independent set in the given abelian torsion free group
\gruppe{}. This description enables us to investigate all linear
orders on $\gruppe$ in one common vector space. Thus, we can use it to
describe all linear order extensions of a partially ordered abelian
group. We concentrated on torsion free groups, as only those can be
ordered linearly. The description used for linear orders can also be
used for semiclosed partially ordered groups, which have been shown to
be compatible with the attempt to catalogue all compatible linear
order extensions of a given partially ordered group. Finally such a
characterisation has been given in Theorem
\ref{satz:ordnungserweiterungscharakterisierung}. Some additional
examples provided an insight, how this method can be applied to
mathematical problems arising from the work with partially ordered
abelian groups.

Though abelian groups have interesting applications, a general
description for arbitrary groups would be interesting. The groups
considered here are torsion free abelian groups and have been embedded
into some vector space $\direktsumme{E}\rationalezahlen$ which is
always possible as the latter is an injective group.  Another idea
would be to consider those groups, which can be embedded into vector
spaces over skew fields.

Linear order extensions play an important role in the modelling of
tone systems (cf. \cite{NeumaierWille1990}) and ordered generalised
interval systems \cite{davi2007generalized} in mathematical
music theory. In this topic it is necessary to describe
factorisations of $\ell$-groups by non-convex subgroups based
on linear order extensions. This leads to further interesting questions
in the direction of cylindrically and cyclically ordered groups. This
theory has its application in software development e.\,g., for
understanding just intonation logics based on the Tonnetz as described
by Martin Vogel \cite{Vogel-1975} and provided
by \textsc{Mutabor} (\cite{Mutabor-home}).

\section{Acknowledgements}
Special thanks to Charles W. Holland for providing  his illustration of Zajceva's
theorem and some further references. This has been the initial inspiration for the current work.

\bibliography{Literatur}
\end{document}